\documentclass[12pt,a4paper,twoside]{amsart}
\usepackage[utf8]{inputenc}

\usepackage[english]{babel}  
\usepackage[T1]{fontenc}     
\usepackage{amsmath,amssymb,amsfonts}
\usepackage{amsthm} 
\usepackage{graphicx}     
\usepackage{ae}
\usepackage{aecompl}
\usepackage{mathrsfs,titletoc}
\usepackage{xspace}

\usepackage[
  top=3cm,
  bottom=3cm,
  inner=2.5cm,      
  outer=2.5cm,
  footskip=40pt     
]{geometry}

\usepackage{multicol}
\usepackage{framed}
\usepackage[all]{xy}
\usepackage{float}
\usepackage{bm}
\usepackage[shortlabels]{enumitem}
\usepackage[usenames,dvipsnames]{xcolor}
\usepackage[colorlinks,linkcolor=BrickRed,citecolor=CadetBlue,urlcolor=black,hypertexnames=true]{hyperref}

\usepackage{tikz}
\usetikzlibrary{shapes,patterns,calc,snakes}

\makeatletter
\newtheorem*{rep@thm}{\rep@title}
\newcommand{\newreptheorem}[2]{%
	\newenvironment{rep#1}[1]{%
		\def\rep@title{#2 \ref{##1}}%
		\begin{rep@thm}}%
		{\end{rep@thm}}}
\makeatother

\newreptheorem{thm}{Theorem}
\newreptheorem{prop}{Proposition}
\newreptheorem{cor}{Corollary}

\theoremstyle{definition} 
\newtheorem{definicija}{Definition}[section]
\newtheorem{primer*}[definicija]{Example}
\newtheorem{opomba}[definicija]{Remark}

\newenvironment{dokaz}[1][Proof]{\begin{proof}[#1]}{\end{proof}}
\newenvironment{primer}[1][]{\begin{primer*}[#1]\pushQED{\qed}}{\popQED\end{primer*}}

\theoremstyle{plain} 
\newtheorem{lema}[definicija]{Lemma}
\newtheorem{izrek}[definicija]{Theorem}
\newtheorem{trditev}[definicija]{Proposition}
\newtheorem{posledica}[definicija]{Corollary}

\makeatletter
\theoremstyle{definition}
\newtheorem*{rep@def}{\rep@title}
\newcommand{\newrepdefinition}[2]{%
	\newenvironment{rep#1}[1]{%
		\def\rep@title{#2 \ref{##1}}%
		\begin{rep@def}}%
		{\end{rep@def}}}
\makeatother

\newrepdefinition{def}{Definition}

\numberwithin{equation}{section}  
\newcommand{\R}{\mathbb R}

\newcommand{\ti}{matrix face\xspace}
\newcommand{\tii}{$C^\ast$-face\xspace}
\newcommand{\tiii}{weak matrix face\xspace}
\newcommand{\tmi}{matrix multiface\xspace}
\newcommand{\tmii}{matrix convex multiface\xspace}

\newcommand{\eti}{matrix exposed face\xspace}
\newcommand{\etii}{$C^\ast$-exposed face\xspace}
\newcommand{\etiii}{weak matrix exposed face\xspace}
\newcommand{\etmi}{matrix exposed multiface\xspace}
\newcommand{\etmii}{matrix convex exposed multiface\xspace}

\makeatletter \g@addto@macro\bfseries{\boldmath} \makeatother

\newtheorem{thmA}{Theorem}

\begin{document}

\setcounter{page}{1}    
\pagenumbering{arabic}  

\setcounter{tocdepth}{3}
\contentsmargin{2.55em} 
\dottedcontents{section}[3.8em]{}{2.3em}{.4pc} 
\dottedcontents{subsection}[6.1em]{}{3.2em}{.4pc}
\dottedcontents{subsubsection}[8.4em]{}{4.1em}{.4pc}

\title{Facial structure of matrix convex sets}
\author{Igor Klep}
\address{Igor Klep: Faculty of Mathematics and Physics, Department of Mathematics,  University of Ljubljana, Slovenia}
\email{igor.klep@fmf.uni-lj.si}
\thanks{IK was supported by the 
Slovenian Research Agency grants 
J1-2453, N1-0217, J1-3004 and P1-0222.}
\author{Tea Štrekelj}
\address{Tea Štrekelj: Institute of Mathematics, Physics and Mechanics, Ljubljana, Slovenia}
\email{tea.strekelj@fmf.uni-lj.si}
\thanks{This is a part of the PhD thesis written by the second author under the supervision of the first author at the University of Ljubljana, Faculty of Mathematics and Physics.}



\maketitle

\begin{abstract} This article investigates the notions of exposed points and (exposed) faces in the matrix convex setting.
	Matrix exposed points in finite dimensions were first defined by Kriel in 2019. Here this notion is extended to matrix convex sets in infinite-dimensional vector spaces. 
	Then a connection between matrix exposed points and matrix extreme points is established:~a matrix extreme point is ordinary exposed if and only if it is matrix exposed. This leads to a Krein-Milman type result for matrix exposed points that is due to Straszewicz-Klee in classical convexity: a compact matrix convex set is the closed matrix convex hull of its matrix exposed points.
	
	Several notions of a fixed-level as well as a multicomponent matrix face and matrix exposed face are introduced to extend the concepts of a matrix extreme point and a matrix exposed point, respectively. Their properties resemble those of (exposed) faces in the classical sense, e.g., it is shown that the $C^\ast$-extreme (matrix extreme) points of a matrix face (matrix multiface) of a matrix convex set \textbf{$K$} are matrix extreme in \textbf{$K$}. As in the case of extreme points, any fixed-level matrix face is ordinary exposed if and only if it is a matrix exposed face.  From this it follows that every fixed-level matrix face of a free spectrahedron is matrix exposed. On the other hand, matrix multifaces give rise to the noncommutative counterpart of the classical theory connecting (archimedean) faces of compact convex sets and (archimedean) order ideals of the corresponding function systems.\looseness=-1
\end{abstract}

\tableofcontents

\section{Introduction}
In the classical theory of convexity an important role is played by distinguished points and subsets of the relative boundary of a convex set. The extreme points ext\,$C$ of a convex set $C$ in a locally convex vector space are those points $c \in C$ that cannot be expressed as a nontrivial convex combination of the elements of $C.$ Equivalently, the set $C \backslash \{c\}$ is convex. Geometrically, any line with an extreme point in its relative interior has at least one of its endpoints outside $C.$ A compact convex set $C$ is the closed convex hull of ext\,$C$ by the Krein-Milman theorem \cite[Section III.4]{Ba}, so in this case the extreme points generate $C.$\looseness=-1

Another type of boundary points of interest are the exposed points exp\,$C.$ These are the ones that can be weakly separated from $C$ by an affine hyperplane and they are easily seen to be extreme. In (finite-dimensional) polyhedra the exposed and extreme points coincide, but in general the exposed points form a proper subset of the extreme boundary. If $C$ is a compact convex set in a normed vector space, the exposed points also recover $C$ as they are dense in ext\,$C$ by the Straszewicz theorem (see \cite[Section II.2]{Ba} and \cite{K}). Hence a compact convex set in a normed vector space is the closed convex hull of its exposed points.\looseness=-1

Extending the concept of extremeness from points to sets leads to the notion of faces. A face $F \subseteq C$ is a convex subset that contains the endpoints of all the lines in $C$ which intersect the relative interior of $F$. Equivalently, faces are those subsets $F$ of $C$ for which $C \backslash F$ is convex. Trivially, $\emptyset$ and $C$ are faces, but any other proper face $F \subsetneq C$ is contained in the relative boundary of $C$. If the face $F$ is a singleton $\{x\}$, then $x$ is an extreme point, meaning that faces extend the concept of extreme points to sets. A key feature of faces is that their extreme points are extreme in the convex set they are contained in. This property features prominently in the proof of the Krein-Milman theorem. In this article these notions are explored in the matrix convex setting.

\subsection{Matrix convex sets} Matrix convex sets are the noncommutative counterpart to classical convex sets, and were first introduced by Wittstock \cite{Wit}. As they are categorically dual to operator systems \cite[Proposition 3.5]{WW}, they introduce convex-geometric ideas and tools to the study operator systems. A Hahn-Banach separation theorem for matrix convex sets was proved by Effros-Winkler \cite{EW}, the matricial Krein-Milman theorem is due to Webster-Winkler \cite{WW}, and further fundamental results were developed recently by a plethora of authors: representations of convex sets by linear matrix inequalities \cite{HV, HM12, FNT}, further results on free spectrahedra including the convex Positivstellensatz \cite{HKM12}, inclusion problems and dilation theory \cite{HKM13, DDOSS17, HL}, minimal and maximal matrix convex sets \cite{Pas} and matrix convex hull approximation \cite{HKM16}, (absolute) extreme points of matrix convex sets and free spectrahedra \cite{EHKM, EH} and matrix exposed points of free spectrahedra \cite{Kr}, the theory of $C^\ast$-convexity, i.e., fixed-dimension matrix convexity, via operator systems and the correspondence between $C^\ast$-extreme points and pure states \cite{FM, F, BM}, noncommutative Choquet theory \cite{DK} and the connection between nonunital operator systems and noncommutative (nc) convex sets \cite{KKM}, the correspondence between compact rectangular matrix convex sets and operator spaces \cite{FHL}, etc.\looseness=-1

Let $V$ be a complex vector space with predual space $V^\prime.$ We will usually refer to $V$ as a dual vector space without explicitly mentioning the predual $V^\prime.$ Denote by $M_{m,n}(V)$ the space of $m \times n$ matrices over $V$ and use the abbreviation $M_n(V) = M_{n,n}(V).$ To simplify notation we write $\mathbb{M}_{m,n} = M_{m,n}(\mathbb{C}),$ $\mathbb{M}_n = M_{n,n}(\mathbb{C}),$ and $\mathbb{I}_n \in \mathbb{M}_n$ for the identity matrix. Unless mentioned otherwise,
we endow $V$ with the weak topology and the corresponding matrix spaces with the product topology. We say a set $\textbf{$S$}= (S_n)_{n \in \mathbb{N}} \subseteq (M_n(V))_{n \in \mathbb{N}}$ is closed (compact) if it is levelwise closed (compact), i.e., each component $S_n$ is closed (compact).

\begin{definicija}
	Suppose for each $n \in \mathbb{N}$ the set $K_n$ is a subset of $M_n(V)$ and denote by \textbf{$K$} the graded family $(K_n)_{n \in \mathbb{N}}$.
	
	(a) Let $A_1,\ldots, A_k \in \textbf{$K$}$ with $A_i \in K_{n_i}$. An expression of the form
	\begin{equation}\label{eq-11}
		\sum_{i=1}^k \gamma_i^\ast A_i \gamma_i,
	\end{equation}
	where $\gamma_i \in \mathbb{M}_{n_i,n}$ are complex matrices with $\sum_{i=1}^k \gamma_i^\ast \gamma_i = \mathbb{I}_n$, is a \textbf{matrix convex combination} of the  points  $A_1,\ldots, A_k$. 
	
	(b) We call \textbf{$K$} a \textbf{matrix convex set} in $V$ if it is closed under 
	matrix convex combinations.
	
	(c) A set $C \subseteq M_n(V)$ is a \textbf{$C^\ast$-convex set} if it is closed under formation of \textbf{$C^\ast$-convex combinations}, i.e., matrix convex combinations \eqref{eq-11} with $A_i \in C$ and $\gamma_i \in \mathbb{M}_n$.
\end{definicija}

Equivalently, a set \textbf{$K$} is matrix convex if and only if it is closed under formation of direct sums and conjugations by isometries. If $0 \in K_1,$ then \textbf{$K$} is closed under conjugations by arbitrary contractions (see Proposition \ref{pr15} below for a proof of this simple observation).
Note that for a matrix convex set \textbf{$K$}, each $K_n$ is a convex set in the classical sense.

For any graded set $\textbf{$S$} = (S_n)_{n \in \mathbb{N}}$ with $S_n \subseteq M_n(V),$ the intersection of all matrix convex sets containing $\textbf{$S$}$ is called the \textbf{matrix convex hull} of $\textbf{$S$}$ and is denoted by mconv\,\textbf{$S$}. Its closure is denoted by $\overline{\text{mconv}}\,\textbf{$S$}$. 

\begin{definicija}
A morphism between matrix convex sets \textbf{$K$} and \textbf{$L$} over spaces $V$ and $W$ respectively, a \textbf{matrix affine map}, is a continuous linear map $\Phi : V \to W$ that satisfies $\Phi_r(K_r) \subseteq L_r$ for all $r \in \mathbb{N}$ and
\begin{align}\label{eq: mtx aff} 
\Phi_r\bigg(\sum_{i=1}^k \gamma_i^\ast A_i \gamma_i\bigg) =
\sum_{i=1}^k\gamma_i^\ast \, \Phi_{r_i}(A_i)  \gamma_i 
\end{align}
for all $k$-tuples $(A_i)_{i=1}^k$ and $(\gamma_i)_{i=1}^k$ such that $A_i \in K_{r_i}$ and $\gamma_i \in \mathbb{M}_{r_i,r}$ for $i = 1,\ldots, k$ with the property $\sum_{i=1}^k \gamma_i^\ast \gamma_i = \mathbb{I}_r$. Here for any positive integer $r$ and $B = (B_{i,j}) \in M_r(V)$ we denote by 
$$\Phi_r(B)  = \big(\Phi(B_{i,j})\big)$$ 
the $r$-th \textbf{ampliation} of $\Phi.$ We call $\Phi$ a \textbf{matrix affine homeomorphism} if each of the $\Phi_r$ is a homeomorphism.
\end{definicija}

\begin{opomba}
	Note that when a matrix affine map $\Phi$ on a matrix convex set \textbf{$K$} maps into  $(M_n(\mathbb{M}_r))_n$ for some $r \in \mathbb{N},$ then by definition, the multiplication in \eqref{eq: mtx aff} means multiplying blocks of $\Phi_{r_i}(A_i)$ with elements (scalars) of $\gamma_i$ or $\gamma_i^*.$
	In this case, by identifying each $M_n(\mathbb{M}_r)$ with $\mathbb{M}_{nr},$ Equation \eqref{eq: mtx aff} can be replaced with
	$$
	\Phi_r\bigg(\sum_{i=1}^k \gamma_i^\ast A_i \gamma_i\bigg) =
	\sum_{i=1}^k (\gamma_i^\ast \otimes \mathbb{I}_r) \, \Phi_{r_i}(A_i)  (\gamma_i \otimes \mathbb{I}_r),
	$$
	where multiplication is now the usual matrix multiplication. We use this equivalent interpretation henceforth.
\end{opomba}

\subsubsection{Matrix extreme points} We recall the definition of the matrix counterparts of classical extreme points: matrix extreme points were introduced by Webster-Winkler \cite{WW}, while $C^\ast$-extreme points already appeared earlier in \cite{FM}.

\begin{definicija}\label{mext}
	Let \textbf{$K$} be a matrix convex set.
	
	(a) A matrix convex combination \eqref{eq-11} is \textbf{proper} if all of the matrices $\gamma_i$ are onto.
	
	(b) A point $A \in K_n$ is \textbf{matrix extreme} if from any expression
	of $A$ as a proper matrix convex combination of elements $A_i \in K_{n_i}$ it follows that $n_i = n$ and each of the $A_i$ is unitarily equivalent to $A$.
	
	(c) A point $A \in K_n$ is a \textbf{$C^\ast$-extreme point} if any expression
	of $A$ as a proper $C^\ast$-convex combination of elements $A_i \in K_{n}$ implies each of the $A_i$ is unitarily equivalent to $A$.
\end{definicija}

Any matrix extreme point of a compact matrix convex set is extreme in the classical sense by \protect{\cite[Corollary 3.6]{WW}}. This also holds for non-compact sets; it is e.g.~an easy corollary of the characterization \protect{\cite[Proposition 4.6]{EHKM}} (see also \cite{HL}).

\subsubsection{Matrix exposed points}
In this paper we study the notions of exposed points and (exposed) faces in the matrix convex setting. While matrix exposed points in the finite dimensional setting were first introduced by Kriel \cite{Kr}, in Section \ref{sec2} we generalize the notion to arbitrary, infinite-dimensional vector spaces, and then investigate their properties and streamline some of the arguments appearing in \cite{Kr}. 

\begin{repdef}{def511}
	Let $\textbf{$K$} = (K_n)_{n \in \mathbb{N}}$ be a matrix convex set in a dual vector space $V.$ An element $A \in K_n$ is called a \textbf{matrix exposed point} of \textbf{$K$} if there exist a continuous linear map $\Phi : V \to \mathbb{M}_n$ and a self-adjoint matrix $\alpha \in \mathbb{M}_n$ such that the following conditions hold:\looseness=-1
	\begin{enumerate}[(a)]
		\item for all positive integers $r$ and $B \in K_r$ we have $\Phi_r(B) \preceq \alpha \otimes \mathbb{I}_r;$ 
		\item $\{ B \in K_n \ | \  \alpha \otimes \mathbb{I}_n - \Phi_n(B) \succeq 0 \text{ singular}\} = \{U^\ast A U \ | \ U \in \mathbb{M}_n \text{ unitary}\}.$
	\end{enumerate}
\end{repdef}

It is straightforward that any exposed point in the classical sense is extreme, although the proof of the matricial analogue, stated as Proposition \ref{tr15}, is more involved, and needs some careful preliminary observations on the exposing map $\Phi$ given in Proposition \ref{lema14}. As expected, matrix exposed points form a proper subset of the matrix extreme points in general (Example \ref{ex36}). \looseness=-1

\subsection{Main results}
Inspired by the Effros-Winkler matricial Hahn-Banach separation techniques developed in \cite{EW}, we establish our first main result, Theorem \ref{thmA}, giving the precise connection between matrix extreme points and matrix exposed points via classical exposed points.

\begin{thmA}\label{thmA}
	Let $\textbf{$K$} = (K_n)_{n \in \mathbb{N}}$ be a matrix convex set. Then a point $A \in K_n$ is matrix exposed if and only if it is a matrix extreme point, which is ordinary exposed in $K_n.$
\end{thmA}


The part of Theorem \ref{thmA} asserting that matrix exposed points are matrix extreme is the above mentioned Proposition \ref{tr15}, while the remaining claims are stated  and proved separately as Theorem \ref{tr111}.
Here the idea in \cite{Kr} of introducing non-archimedean real closed fields is key to allow the Effros-Winkler separation techniques to apply in the context of a general (not necessarily closed) matrix convex set. The idea of separating over a real closed field extension of $\R$ also appears in the theory of convexity over arbitrary ordered fields, developed in \cite{SMR}, the related separation techniques in \cite{Rob} and the real closed separation theorem for convex sets in \cite{NT}. The methods used in the proof of Theorem \ref{tr111} also yield an Effros-Winkler type weak separation theorem for (non-closed) matrix convex sets.

\begin{repcor}{th310}[Weak separation theorem for matrix convex sets]
	Let \textbf{$K$} be a matrix convex set in a dual space $V$ with $0 \in K_1$ and $A \notin K_n.$ Suppose there is a continuous linear functional $\varphi: M_n(V) \to \mathbb{C}$ and real number $a > 0$ such that $\text{\textnormal{Re}}\,\varphi|_{K_n} < a$ and $\varphi(A) = a.$ Then there exists a continuous linear map $\Phi : V \to \mathbb{M}_n$ such that
	$$
	\mathbb{I}_n \otimes \mathbb{I}_r  - \text{\textnormal{Re}}\, \Phi_r(B) \succ 0
	$$
	for every positive integer $r$ and $B \in K_r,$ but
	$$
	\text{\textnormal{ker}}\,\big(\mathbb{I}_n \otimes \mathbb{I}_n  - \text{\textnormal{Re}}\, \Phi_n(A)\big) \neq \{0\}.
	$$
\end{repcor}

The second important result in Section \ref{sec2} is Theorem \ref{th213} giving a matrix analogue of the classical Straszewicz theorem \cite[Section II.2]{Ba}, more precisely, its generalisation to compact convex sets in normed spaces due to Klee \cite{K}. 

\begin{repthm}{th213}[\textbf{The Straszewicz-Klee theorem for matrix convex sets}]
	Let \textbf{$K$} be a compact matrix convex set in a normed vector space $V.$ Then $\text{\textnormal{mexp}}\,\textbf{$K$} \neq \emptyset$ and
	$$
	\textbf{$K$} = \overline{\text{\textnormal{mconv}}}\,(\text{\textnormal{mexp}}\,\textbf{$K$}).
	$$
\end{repthm}
The proof goes along the lines of the Webster-Winkler matricial Krein-Milman theorem \cite{WW} in combination with the techniques in \cite{HL} of assigning to a matrix convex set \textbf{$K$} a family of convex sets, whose exposed points are shown to be in correspondence with the matrix exposed points of \textbf{$K$}.

\subsubsection{Exposed points of state spaces} 
A (concrete) \textbf{operator system} $\mathcal{R}$ is a closed self-adjoint subspace of the operators on a Hilbert space that contains the identity. As before, for positive integers $r, n$ and  a linear map $\varphi: \mathcal{R} \to \mathbb{M}_n,$ the $r$-th ampliation $\varphi_r: M_r(\mathcal{R}) \to M_r(\mathbb{M}_n)$ is defined by applying $\varphi$ entrywise, i.e.,
$$
\varphi_r\big((A_{i,j})\big) = \big(\varphi(A_{i,j})\big)
$$ 
for $(A_{i,j}) \in M_r(\mathcal{R}).$ The map $\varphi$ is \textbf{completely positive} (cp) if for all $r \in \mathbb{N},$ the $r$-th ampliation $\varphi_r$ is positive, meaning that if $A \in M_r(\mathcal{R})$ is positive semidefinite, so is $\varphi_r(A) \in M_r(\mathbb{M}_n) \cong \mathbb{M}_{rn}.$ The collection of all unital completely positive (ucp) maps from an operator system $\mathcal{R}$ to the matrix spaces $\mathbb{M}_n$ for $n \in \mathbb{N}$ is referred to as the \textbf{matrix state space} of $\mathcal{R}$ and is easily seen to be matrix convex. In fact, it is the core example of a compact matrix convex set by the categorical duality established in \cite[Proposition 3.5]{WW}.
Further, the matrix extreme points of the state space of an operator system $\mathcal{R}$ are precisely the pure states on $\mathcal{R}$ (see \cite{F}). 

In Subsection \ref{subsec34} we give some insight into the matrix exposed points of the matrix state space of an operator system $\mathcal{R}$, while keeping in mind that by \cite[Theorem 2.2]{F}, the extreme rays in the space of cp maps on $\mathcal{R}$ are determined by the matrix extreme points of the state space of $\mathcal{R}$. Proposition \ref{tr116} presents the analogous connection between the exposed rays and matrix exposed points of the respective sets. If $\mathcal R=\mathcal A$ is a separable $C^*$-algebra, then every 
matrix extreme point of the matrix state space of $\mathcal A$ is
matrix exposed (Example \ref{ex:cstarexposed}).

\subsubsection{Matrix faces and matrix exposed faces}
To investigate the facial structure of a matrix convex set \textbf{$K$} we discuss several possible notions of a face and an exposed face of \textbf{$K$}. Here the main distinction is whether one considers subsets of a single component $K_n$ for some $n \in \mathbb{N}$ or multicomponent subsets of a matrix convex set \textbf{$K$}.
Section \ref{sec3} first introduces three concurrent definitions of a fixed-level matrix face, which aim to extend the concepts of a matrix extreme point or a matrix exposed point, and demonstrates their suitableness to the theory of matrix convexity.

\begin{repdef}{def41}
	Let $\textbf{$K$} = (K_r)_{r \in \mathbb{N}}$ be a matrix convex set in the space $V$ and $F$ a convex subset of $K_n$ for some $n \in \mathbb{N}$.
	
	(a) Then $F$ is a \textbf{\ti} if for every tuple of points $A_1, \ldots, A_k$ from \textbf{$K$} with $A_i \in K_{n_i}$ and every tuple of surjective matrices $\gamma_i \in \mathbb{M}_{n_i, n}$ satisfying $\sum_{i=1}^k \gamma_i^\ast \gamma_i = \mathbb{I}_n,$ the condition
	\begin{equation} 
		\sum_{i=1}^k \gamma_i^\ast A_i \gamma_i \in F,
	\end{equation}
	implies $n_i=n$ and $A_i \in F$ for $i=1, \ldots, k.$
	
	(b) If $F$ is a $C^\ast$-convex \ti, then it is a \textbf{\tii}.
	
	(c) The set $F$ is a \textbf{\tiii} if for every tuple of points $A_1, \ldots, A_k$ from \textbf{$K$} with $A_i \in K_{n_i}$ and every tuple of surjective matrices $\gamma_i \in \mathbb{M}_{n_i, n}$ satisfying $\sum_{i=1}^k \gamma_i^\ast \gamma_i = \mathbb{I}_n,$ the condition
	\begin{equation*}
		\sum_{i=1}^k \gamma_i^\ast A_i \gamma_i \in F,
	\end{equation*}
	implies $n_i=n$ and each $A_i$ is unitarily equivalent to some element in $F.$ We will denote by $\mathcal{U}(F) = \{U^\ast A U \ | \ A \in F, \ U \in \mathbb{M}_n \text{ unitary}\}$
	the unitary orbit of $F.$
\end{repdef}

We show that a key hereditary property of extreme points of classical faces has its matrix counterpart. 

\begin{repprop}{tr34}
	Let \textbf{$K$} be a matrix convex set and $F \subseteq K_n$ a matrix face of any type. Every $C^\ast$-extreme point of $F$ is a matrix extreme point of \textbf{$K$}.
\end{repprop}

In Subsection \ref{subsec42}, the corresponding three types of matrix exposed faces are introduced.

\begin{repdef}{def21}
	Let $\textbf{$K$} = (K_r)_{r \in \mathbb{N}}$ be a matrix convex set in a dual vector space $V$ and $F$ a convex subset of $K_n$. 
	
	(a) Then $F$ is a \textbf{\eti} if there exists a continuous linear map $\Phi:V \to \mathbb{M}_n$ and a self-adjoint matrix $\alpha \in \mathbb{M}_n$ satisfying the following conditions:
	\begin{enumerate}[(i), leftmargin=2cm]
		\item 
		for every positive integer $m$ and $B \in K_m$ we have $\Phi_m(B) \preceq \alpha \otimes \mathbb{I}_m;$ 
		
		\item for any $m < n$ and $B \in K_m$ we have $\Phi_m(B) \prec \alpha \otimes \mathbb{I}_m;$
		
		\item 
		$\{ B \in K_n \ | \  \alpha \otimes \mathbb{I}_n - \Phi_n(B) \succeq 0 \text{ is singular}\} = F.$
	\end{enumerate}
	
	(b) If $F$ is a $C^\ast$-convex \eti, then it is a \textbf{\etii}.
	
	(c) We call $F$ a \textbf{\etiii} if there exists a continuous linear map $\Phi:V \to \mathbb{M}_n$ and a self-adjoint matrix $\alpha \in \mathbb{M}_n$ satisfying the following conditions:
	\begin{enumerate}[(i), leftmargin=2cm]
		\item  for every positive integer $m$ and $B \in K_m$ we have $\Phi_m(B) \preceq \alpha \otimes \mathbb{I}_m;$
		
		\item for any $m < n$ and $B \in K_m$ we have $\Phi_m(B) \prec \alpha \otimes \mathbb{I}_m;$ 
		
		\item $\{ B \in K_n \ | \  \alpha \otimes \mathbb{I}_n - \Phi_n(B) \succeq 0 \text{ is singular}\} = \mathcal{U}(F).$
	\end{enumerate}
\end{repdef}
 After adapting the observations in Proposition \ref{lema14} on the exposing maps for matrix exposed points to the matrix face setting, we prove the following expected, yet not entirely obvious result. 
 
 \begin{repprop}{tr412}
 	Let \textbf{$K$} be a closed matrix convex and $F \subsetneq K_n$ a matrix exposed face of any type. Then $F$ is a matrix face of the corresponding type.
 \end{repprop}
 
 We proceed by giving a generalisation of Theorem \ref{tr111} for faces, namely Theorem \ref{tr211}, presenting an interplay between matrix faces and matrix exposed faces. 
 
 \begin{repthm}{tr211} Let \textbf{$K$} be a matrix convex set and $F \subseteq K_n$ a matrix face of any type that is also an exposed face.  Then $F$ is a matrix exposed face of the corresponding type.
 \end{repthm}
 
 Finally, Proposition \ref{tr416} and Corollary \ref{tr315} give a sufficient condition for a point to lie in a \tiii and a \ti, respectively. This leads to a family of examples of weak matrix faces presented in Example \ref{ex416}.
 We then observe in Subsection \ref{sub331} that as a corollary of Theorem \ref{tr211}, every \ti of a free spectrahedron is matrix exposed. 
 
\subsubsection{Matrix multifaces and matrix exposed multifaces}

In Section \ref{sec5} we discuss two aspirant notions of a multicomponent face of a matrix convex set. The main attribute of multilevel faces is their role in the noncommutative counterpart to the classical theory connecting (archimedean) faces of compact convex sets and (archimedean) order ideals of the corresponding function systems presented in \cite[Section II.5]{Al}.

While our aim is to extend the properties of a matrix extreme point, a notion similar to that of a matrix multiface, but mimicking absolute extreme points (see \cite{EHKM}), was recently explored under the name nc face in \cite{KKM}.

\begin{repdef}{def51}
	Let $\textbf{$K$} = (K_r)_{r \in \mathbb{N}}$ be a matrix convex set in the space $V$ and $\textbf{$F$} = (F_r)_{r \in \mathbb{N}} \subseteq \textbf{$K$}$ a levelwise convex subset of \textbf{$K$}.  
	
	(a) Then $\textbf{$F$}$ is a \textbf{\tmi} if for every tuple of points $A_1, \ldots, A_k$ from \textbf{$K$} and every tuple of surjective matrices $\gamma_i \in \mathbb{M}_{n_i, n}$ satisfying $\sum_{i=1}^k \gamma_i^\ast \gamma_i = \mathbb{I}_n,$ the condition
	\begin{equation} 
		\sum_{i=1}^k \gamma_i^\ast A_i \gamma_i \in \textbf{$F$},
	\end{equation}
	implies $A_i \in \textbf{$F$}$ for $i=1, \ldots, k.$
	
	(b) If $\textbf{$F$}$ is a matrix convex \tmi, then it is a \textbf{\tmii}.
\end{repdef}

Inspired by \cite[Section II.5]{Al} we give a family of examples of matrix convex multifaces. For a compact matrix convex set \textbf{$K$} we denote by \looseness=-1
$$
A(\textbf{$K$}) := \{\theta = (\theta_n: K_n \to \mathbb{M}_n)_{n \in \mathbb{N}} \ | \ \theta \text{ continuous matrix affine}\}
$$ 
its dual operator system, which is abstractly characterized by the Choi-Effros axioms as a matrix-ordered $\ast$-vector space with an Archimedean matrix order unit (see \cite[Chapter 13]{Pa}). Recall that by \cite[Proposition 3.5]{WW}, \textbf{$K$} is matrix affinely homeomorphic to UCP$(A(\textbf{$K$}))$ by the evaluation map sending $X \in \textbf{$K$}$ to the ucp map $\Phi_X \in \text{UCP}(A(\textbf{$K$})),$ where
$$
\Phi_X(\theta) = \theta(X)
$$
for any $\theta \in A(\textbf{$K$}).$ Now suppose $\Phi$ is a ucp map on $A(\textbf{$K$})$ with kernel $J$ and let 
\begin{equation} 
	J_n^{\perp} := \{A \in K_n \ | \ \theta_n(A) = 0 \ \  \forall \theta \in J\}
\end{equation}
for $n \in \mathbb{N}.$ We show in Example \ref{ex418} that if $J$ is spanned by its positive elements, i.e., $J= J^+ - J^+,$ then $\textbf{$J$}^{\perp}:= (J_n^{\perp})_{n \in \mathbb{N}}$ is a \tmii of  \textbf{$K$}. In Remark \ref{op422} we then observe how this construction gives a sufficient condition for a point to be contained in some matrix multiface. 

The matrix convex version of the interplay between classical faces and order ideals is stated as Proposition \ref{tr423}.
Let us say that a multicomponent subset $\textbf{$F$} \subseteq \textbf{$K$}$ of a compact matrix convex set \textbf{$K$} satisfies condition $(\ast)$ if for each $n \in \mathbb{N}$ and $\theta \in M_n(A(\textbf{$K$}))$ with $\theta|_{\textbf{$F$}} \succeq 0$ there is a positive element $\psi \in M_n(A(\textbf{$K$}))^+$ such that 
$$\psi \succeq \theta \quad \text{ and } \quad \psi|_{\textbf{$F$}} = \theta|_{\textbf{$F$}}.
$$ 
Further, a ucp map $\Phi: A(\textbf{$K$}) \to \mathbb{M}_n$ is called \textbf{partially order reflecting} if it satisfies 
\begin{equation*}
	\Phi_m\big(M_m(A(\textbf{$K$}))^+\big) = \Phi_m\big(M_m(A(\textbf{$K$}))\big)^+
\end{equation*}
for all $m \in \mathbb{N}$, i.e., for every $m \in \mathbb{N}$ and $A \in M_m(A(\textbf{$K$}))$ with $\Phi_m(A) \succeq 0$ there exists a $B \succeq 0$ such that $\Phi_m(A) = \Phi_m(B).$

\begin{repprop}{tr423} Let \textbf{$K$} be a compact matrix convex set.
	\begin{enumerate}[\text\textnormal{{(a)}}]
		\item Let $n \in \mathbb{N}$ and $\Phi: A(\textbf{$K$}) \to \mathbb{M}_n$ be a partially order reflecting ucp map with kernel $J$ spanned by its positive elements. Then $\textbf{$J$}^{\perp} \subseteq\textbf{$K$}$ is a closed \tmii that satisfies condition $(\ast)$.
	\end{enumerate}
	\begin{enumerate}[\text\textnormal{{(b)}}]
		\item Suppose $\textbf{$F$}  \subseteq \textbf{$K$}$ is a closed \tmii that satisfies $(\ast)$.
			Then
		$$
		J := \{\theta \in A(\textbf{$K$}) \ | \ \theta|_{\textbf{$F$}} = 0\}
		$$
		is spanned by its positive elements and is the kernel of a ucp map $\Phi: A(\textbf{$K$}) \to \mathcal{R}$ for some operator system $\mathcal{R},$ where $\Phi$ satisfies the partially order reflecting property $\Phi_n\big(M_n(A(\textbf{$K$}))^+\big) = \Phi_n\big(M_n(A(\textbf{$K$}))\big)^+$ for all $n \in \mathbb{N}.$
	\end{enumerate}
\end{repprop}

Example \ref{ex519} then explains how
every vertex of a simplex $S$ in an Euclidean space $\mathbb{R}^n$ lies in a \tmii of mconv$(S).$ More precisely, each vertex defines a partially order reflecting evaluation map whose kernel is spanned by its positive elements. 

Next, a hereditary property of matrix extreme points is established.

\begin{repprop}{tr58}
	Let \textbf{$K$} be a matrix convex set and $\textbf{$F$} \subseteq \textbf{$K$}$ a matrix \text{\textnormal{(}}convex\,\text{\textnormal{)}} multiface of any type. Every matrix extreme point of \textbf{$F$} is a matrix extreme point of \textbf{$K$}.
\end{repprop}

Section \ref{subsec52} introduces the exposed counterparts of the multilevel matrix faces and investigates their properties. For instance, in Proposition \ref{pr428} every component of a matrix exposed multiface is shown to be an ordinary exposed face.

\begin{repdef}{def59}
	Let $\textbf{$K$} = (K_r)_{r \in \mathbb{N}}$ be a matrix convex set in a dual vector space $V$ and \textbf{$F$} a levelwise convex subset of \textbf{$K$}. 
	
	(a) Then \textbf{$F$} is a \textbf{\etmi} if there exists a positive integer $r,$ a continuous linear map $\Phi:V \to \mathbb{M}_r$ and a self-adjoint matrix $\alpha \in \mathbb{M}_r$ satisfying the following conditions:
	\begin{enumerate}[(i), leftmargin=2cm]
		\item for every positive integer $n$ and $B \in K_n$ we have $\Phi_n(B) \preceq \alpha \otimes \mathbb{I}_n;$ 
		
		\item for each $n \in \mathbb{N}$ we have $\{ B \in K_n \ | \  \alpha \otimes \mathbb{I}_n - \Phi_n(B) \succeq 0 \text{ is singular}\} = F_n.$
	\end{enumerate}
	
	(b) If \textbf{$F$} is a matrix convex \etmi, then it is a \textbf{\etmii.}
\end{repdef}

Subsection \ref{subsec53} is an extension of Subsection \ref{subsec43} and explores the connection between matrix multifaces and matrix exposed multifaces.

\subsection{Reader's guide} This paper is organized as follows.
Section \ref{sec1} contains basic definitions and preliminaries in both classical and matrix convexity. Section \ref{sec2} then deals with generalising the notion of an exposed point and its properties to the matrix setting. It includes the proofs of Proposition \ref{tr15} and Theorem \ref{tr111} giving the interplay between matrix extreme and matrix exposed points. Subsection \ref{subsec33} is dedicated to the proof of the Straszewicz-Klee theorem for matrix convex sets (Theorem \ref{th213}), while Subsection \ref{subsec34} deals with exposed points of state spaces. Section \ref{sec3} introduces fixed-level matrix faces and matrix exposed faces and establishes their connection in Proposition \ref{tr412} and Theorem \ref{tr211}. The extreme points hereditary property of matrix faces is proved as Proposition \ref{tr34}, and in Subsection \ref{sub331} the correspondence between matrix faces and matrix exposed faces for free spectrahedra is deduced. Section \ref{sec5} covers  matrix multifaces and matrix exposed multifaces. Their correspondence with the kernels of partially order reflecting ucp maps is established in Proposition \ref{tr423}, while the hereditary property of matrix extreme points is stated as Proposition \ref{tr58}. 

\subsection*{Acknowledgments} We express our special thanks to Eric Evert for his insightful comments and valuable suggestions and are appreciative of the helpful comments on the earlier versions of the manuscript provided by Jurij Volčič and Scott McCullough. 
We thank Raphaël Clouâtre for notifying us of an issue in an earlier version of the paper and we thank the anonymous referee for their detailed reading and thoughtful suggestions.

\section{Preliminaries}\label{sec1}

We recall the formal definitions of extreme and exposed points from the classical theory, as well as of their set analogues, faces and exposed faces (see \cite{Ba}). We then present a convenient translation argument and a property of interior points of matrix convex sets. 
Lastly, we give some background on free spectrahedra in Subsection \ref{subsec21}. 

\begin{definicija} \label{def14}
	Let $K \subseteq V$ be a convex set.
	
	(a) A point $x \in K$ is called an \textbf{extreme point} of $K$ if any expression $x = t y + (1-t) z$ for some $y, z \in V$ and $0 < t < 1$ forces $x = y = z.$ Equivalently, the set $K\backslash\{x\}$ is convex.
	
	(b) A point $x \in K$ is an \textbf{exposed point} of $K$ if there exists a continuous functional $\varphi: V \to \mathbb{C}$ and a real number $a$ such that $\varphi(x) = a$ and $\varphi(y) < a$ for all $y \in K \backslash \{x\}.$
\end{definicija}

Every exposed point is extreme, while the converse holds, e.g., for (finite-dimensional) polyhedra, but not in general (see Figure \ref{izpneext} below). 
The next definition extends the concepts of extreme and exposed points to sets.

\begin{definicija}\label{def17}
	Let $K \subseteq V$ be a convex set.
	
	(a) A convex subset $F \subseteq K$ is called a \textbf{face} of $K$ if $t x + (1-t) y \in F$ for some $x, y \in K$ and $0 < t < 1$ forces $x, y \in F.$ Equivalently, the set $K\backslash F$ is convex.
	
	(b) A convex subset $F \subseteq K$ is an \textbf{exposed face} of $K$ if there exists a continuous functional $\varphi: V \to \mathbb{C}$ and a real number $a$ such that $\varphi(x) = a$ for all $x \in F$ and $\varphi(y) < a$ for all $y \in K \backslash F.$
\end{definicija}

Every exposed face is indeed a face and it is straightforward that for a singleton (exposed) face $F = \{x\},$ the point $x$ is extreme (exposed). We now explain some technical assumptions on the matrix convex set \textbf{$K$} we are considering that will appear throughout the paper.
If convenient, we may assume $0 \in K_1$ since one can instead consider the matrix convex set $-\lambda \,+\, \textbf{$K$}$ for some $\lambda \in K_1.$ This assumption is usually made without loss of generality as translations preserve matrix extreme points, etc.

\begin{trditev} \label{pr15}
	A matrix convex set \textbf{$K$} in a dual space $V$ with $0 \in K_1$ is closed under conjugation by contractions.
\end{trditev}

\begin{dokaz}
	Let $A \in K_r$ and let $\alpha \in \mathbb{M}_{r, n}$ be a contraction. Since $0 \in K_1,$ we have $0_n = \oplus_n 0 \in K_n.$ Letting $\beta = (\mathbb{I}_n - \alpha^\ast \alpha)^\frac{1}{2},$ we have $\alpha^\ast \alpha + \beta^\ast \beta = \mathbb{I}_n$ and 
	\[
	\alpha^\ast A \alpha = (\alpha^\ast\  \beta^\ast)(A \oplus 0_{n-r}) 
	\begin{pmatrix} 
		\alpha  \\
		\beta 
	\end{pmatrix} \in K_n. \qedhere
    \]
\end{dokaz}

\begin{trditev} \label{tr18}
	If \textbf{$K$} is a matrix convex set and $v \in \text{\textnormal{int}}\,K_1,$ then $\oplus_n v \in \text{\textnormal{int}}\,K_n.$
\end{trditev}

\begin{dokaz}
We may without loss of generality assume that 
$v = 0 \in \text{int}\,K_1.$ Otherwise \textbf{$K$} can be replaced by $-v + \textbf{$K$}$ so that $0 \in \text{\textnormal{int}}\, (-v + K_1).$ Then $\oplus_n 0 \in \text{\textnormal{int}}\,(-\oplus_n v + K_n)$ implies $\oplus_n v \in \text{\textnormal{int}}\,K_n.$

	Recall that open neighbourhoods of $0$ in the weak topology are of the form
	$$
	U_{v^{\prime}, \epsilon} = \{w \in V \ | \ |\langle w, v^{\prime} \rangle| < \epsilon \}
	$$
	where $v^{\prime} \in V^{\prime},$ $\epsilon>0$ and $\langle \cdot, \cdot \rangle : V \times V^{\prime} \to \mathbb{C}$ denotes the pairing of $V$ and $V^{\prime}.$ Also, for each $n \in \mathbb{N},$ a pairing of the matrix spaces $M_n(V)$ and $M_n(V^{\prime})$ can be defined by 
	\begin{equation*}
		\langle\!\langle B, B^{\prime} \rangle\!\rangle = \sum_{i, j} \langle B_{i, j}, B^{\prime}_{i, j} \rangle
	\end{equation*}
	for $B \in M_n(V)$ and $B^{\prime} \in M_n(V^{\prime}).$ A net in $M_n(V)$ converges weakly if and only if it converges entrywise (see, e.g., \cite[Section 2]{EW}) and open neighbourhoods of $0$ in the weak topology of $M_n(V)$ are of the form 
	$$
	U^n_{B^{\prime}, \epsilon} = \{B \in M_n(V) \ | \ |\langle B_{i, j}, B^{\prime}_{i, j}  \rangle| < \epsilon \ \text{ for }\  i, j=1, \ldots, n\},
	$$
	where $B^{\prime} \in M_n(V^{\prime})$ and $\epsilon>0.$

	As $0 \in \text{int}\,K_1,$ there is a $v^{\prime} \in V^{\prime}$ and $\epsilon>0$ such that $U_{v^{\prime}, \epsilon} \subseteq K_1.$ Let  $w \in U_{v^{\prime}, \epsilon}=-U_{v^{\prime}, \epsilon}.$ 
	Matrix convexity of \textbf{$K$} and the unitary similarity
	of $
		\begin{pmatrix}
			0 & 1 \\
			1 & 0
		\end{pmatrix}$
and
$\begin{pmatrix}
		1 & 0 \\
		0 & -1
	\end{pmatrix}
	$,
	as well as
of $
		\begin{pmatrix}
			0 & 1 \\
			-1 & 0
		\end{pmatrix}$
and
$\begin{pmatrix}
		i & 0 \\
		0 & -i
	\end{pmatrix}
	$
	imply
	$$
	\begin{pmatrix}
		0 & w \\
		w & 0
	\end{pmatrix},
	\begin{pmatrix}
		0 & w \\
		-w & 0
	\end{pmatrix} \in K_2.
	$$
	Whence
	$$
	\begin{pmatrix}
		0 & w \\
		0 & 0
	\end{pmatrix} =
\frac{1}{2}\bigg(\begin{pmatrix}
	0 & w \\
	w & 0
\end{pmatrix} +
\begin{pmatrix}
	0 & w \\
	-w & 0
\end{pmatrix}\bigg) \in K_2.
	$$
It is easy to see that also $w\, E_{i,j}\in K_n,$ where $(E_{i, j})_{i, j}$ denote the standard $n\times n$ matrix units.
	
Now take $B^{\prime} = \big(\frac{v^\prime}{n^2}\big)_{i, j} \in M_n(V^{\prime})$ and observe that for any $B \in U^n_{B^{\prime}, \epsilon},$  
$$n^2 B_{i, j} E_{i,j}\in K_n$$ for all $i, j.$ Then again by (matrix) convexity,
$$
B = \frac{1}{n^2} \sum_{i, j} n^2 B_{i, j} E_{i, j} \in K_n.
$$
We deduce that $\oplus_n 0 \in U^n_{B^{\prime}, \epsilon} \subseteq K_n.$
\end{dokaz}

\subsection{Free spectrahedra}\label{subsec21} We now describe an important class of matrix convex sets arising from spectrahedra.
Taking $V = \mathbb{C}^g$ for some $g \in \mathbb{N}$ we get $M_n(V) \cong \mathbb{M}_n^g.$ For $k \in \mathbb{N}$ denote by $\mathbb{S}^g_k$ the space of $g$-tuples of complex self-adjoint $k \times k$ matrices. For $A = (A_0, \ldots, A_g) \in \mathbb{S}^g_k,$ the corresponding linear matrix-valued polynomial
$$
L_A = A_0 + \sum_{i=1}^g A_i x_i
$$ 
in the noncommuting variables $x_1, \ldots, x_g$ is called a \textbf{linear pencil}. It can be evaluated at a point $x \in \R^g,$ producing a Linear Matrix Inequality  $L_A(x) \succeq 0$ with the solution set $\{x \in \R^g \ | \ L_A(x) \succeq 0\}$ called a \textbf{spectrahedron} (see, e.g., \cite{RG, HV, HM12}). Similarly, $L$ is evaluated at a tuple $X \in \mathbb{S}^g_n$ as 
$$
L_A(X) = A_0 \otimes \mathbb{I}_n + \sum_{i=1}^g A_i \otimes X_i,
$$
where $\otimes$ denotes the Kronecker (tensor) product.
Then the matricial solution set $\mathcal{D}_{A} = \big(\mathcal{D}_{A}(n)\big)_n,$ where
$$
\mathcal{D}_{A}(n) = \{X \in \mathbb{S}^g_n \ | \ L_A(X) \succeq 0 \}
$$
is referred to as a \textbf{free spectrahedron} (see, e.g., \cite{HKM12, HKM13, FNT, EH, Kr}) and is easily seen to be matrix convex.

We will often assume that $\mathcal{D}_{A}(1)$ has nonempty interior. This may be done without loss of generality as we now explain. Since $\mathcal{D}_{A}(1)$ is a finite-dimensional convex set, it has nonempty relative interior, i.e., nonempty interior in the relative topology of its affine span. This means it is contained in a proper affine subspace of $\mathbb{R}^g,$ i.e., 
\begin{equation}\label{eq1}
\varphi |_{\mathcal{D}_{A}(1)} = a
\end{equation} 
for some functional $\varphi : \mathbb{R}^g \to \mathbb{R}$ and $a \in \mathbb{R},$
which implies 
$$
	\varphi_n |_{\mathcal{D}_{A}(n)} = a \otimes \mathbb{I}_n
$$
for all $n \in \mathbb{N}$ (see \cite[Corollary 3.6]{HKM16}). So if $\mathcal{D}_{A}(1)$ has no interior points, we can use the relations given by \eqref{eq1} to express some of the variables $x_i$ in terms of the others and thus reduce dimensions. 

The assumption int\,$\mathcal{D}_{A}(1) \neq \emptyset$ in turn implies that $L$ can be assumed to be \textbf{monic}, meaning $A_0 = \mathbb{I}_k$ (see \cite[Proposition 2.1]{HKM13}).
Then by the Effros-Winkler matricial Hahn-Banach separation theorem \cite{EW,HM12}, a spectrahedron is the matrix analogue of an affine half-space with the corresponding affine hyperplane being
\[
\partial \mathcal{D}_{A}(n) = \{X \in \mathbb{S}^g_n \ | \ L_A(X) \succeq 0 \text{ singular}\}.\qedhere
\]

\section{Matrix exposed points}\label{sec2}

In this section we introduce and study matrix exposed points in matrix convex sets, a notion originating in \cite{Kr}. We generalise the notion to arbitrary infinite-dimensional vector spaces while at the same time streamlining many of the arguments from \cite{Kr}. The main results are Theorem \ref{tr111} asserting the connection between matrix exposed points and matrix extreme points and Theorem \ref{th213} giving a matrix analogue of the classical Straszewicz theorem, more precisely, its generalisation due to Klee stating that every compact convex set in a normed space is the closed convex hull of its exposed points. 

\subsection{Definition and basic properties}\label{subsec31} The definition of a matrix exposed point aims to simulate properties of an exposed point in the classical sense. Throughout we assume that \textbf{$K$} is a matrix convex set in a dual vector space $V,$ endowed with the corresponding weak topology with respect to which the involved linear maps are assumed continuous.
We will also assume that $K_1$ has more than one point and (after a translation if needed) that $0 \in K_1.$ 

\begin{definicija}\label{def511}
	Let $\textbf{$K$} = (K_n)_{n \in \mathbb{N}}$ be a matrix convex set in a dual vector space $V$. An element $A \in K_n$ is called a \textbf{matrix exposed point} of \textbf{$K$} if there exist a continuous linear map $\Phi : V \to \mathbb{M}_n$ and a self-adjoint matrix $\alpha \in \mathbb{M}_n$ such that the following conditions hold:\looseness=-1
	\begin{enumerate}[(a)]
		\item \label{p1} for all positive integers $r$ and $B \in K_r$ we have $\Phi_r(B) \preceq \alpha \otimes \mathbb{I}_r;$ 
		
		\item \label{p2} $\{ B \in K_n \ | \  \alpha \otimes \mathbb{I}_n - \Phi_n(B) \succeq 0 \text{ singular}\} = \{U^\ast A U \ | \ U \in \mathbb{M}_n \text{ unitary}\}.$
	\end{enumerate}
	We say that a pair $(\Phi, \alpha)$  \textbf{matricially exposes} the point $A$ and denote the set of all matrix exposed points of \textbf{$K$} by mexp\,$\textbf{$K$}$.
\end{definicija}

\begin{opomba} \label{op12}
	(a) For any linear map $\Phi : V \to \mathbb{M}_n$ the corresponding family $(\Phi_r|_{K_r}: K_r \to M_r(\mathbb{M}_n))_{r\in \mathbb{N}}$ of restricted canonical ampliations defines a matrix affine map, i.e.,~it satisfies:
	$$
	\Phi_r\bigg(\sum_{i=1}^k \gamma_i^\ast A_i \gamma_i\bigg) =
	\sum_{i=1}^k(\gamma_i^\ast \otimes \mathbb{I}_r)\, \Phi_{r_i}(A_i)  (\gamma_i \otimes \mathbb{I}_r)
	$$
	for all $k$-tuples $(A_i)_{i=1}^k$ and $(\gamma_i)_{i=1}^k$ such that $A_i \in K_{r_i}$ and $\gamma_i \in \mathbb{M}_{r_i,r}$ for $i = 1,\ldots, k$ with the property $\sum_{i=1}^k \gamma_i^\ast \gamma_i = \mathbb{I}_r$.
	
	(b) \label{op12b} Recall that an ordinary exposed point $A \in K$ can be weakly separated from the other points of a convex set $K$ by an affine hyperplane. In other words, for the functional $\varphi:V\to\mathbb{C}$ and $a,$ which determine the hyperplane, we have that $A$ is the only point of $K$ in the kernel of the map $a - \varphi,$ while for all $x$ in $K$ we have $\varphi(x) \leq a.$ 
	On the other hand, if a linear map $\Phi$, matrix $\alpha \in \mathbb{M}_n$ and  $A \in K_n$ are as in Definition \ref{def511}, then we have for any unitary matrix $U \in \mathbb{M}_n,$
	\begin{align*}
		\alpha \otimes \mathbb{I}_n - \Phi_n(U^\ast A U) &= (U^\ast \otimes \mathbb{I}_n) \big(  \alpha \otimes \mathbb{I}_n - \Phi_n(A)\big) (U \otimes \mathbb{I}_n).
	\end{align*}
	Note that if the matrix $\alpha \otimes \mathbb{I}_n - \Phi_n(A)$ is singular, then so is $\alpha \otimes \mathbb{I}_n - \Phi_n(U^\ast A U).$ Condition (b) of Definition \ref{def511} additionally demands for the points of the unitary orbit of $A$ to be exactly the ones from $K_n$ in the kernel of the map $\alpha \otimes \mathbb{I}_n - \Phi_n$. We conclude that if $A$ is matrix exposed, then so is any point from its unitary orbit (being exposed by the same pair $(\Phi, \alpha)$ as $A$).

	(c) \label{c-op12} From Definition \ref{def511} we see that for $r < n$ and $B \in K_r$ the strict inequality $\Phi_r(B) \prec \alpha \otimes \mathbb{I}_r$ holds as we now explain. If $r < n$ and $B \in K_r$ are such that $\alpha \otimes \mathbb{I}_r - \Phi_r(B)$ is singular (while also positive semidefinite), then for 
	any $C \in K_{n-r},$
	$$
	 \alpha \otimes \mathbb{I}_n - \Phi_n(B \oplus C) = \big(\alpha \otimes \mathbb{I}_r - \Phi_r(B)\big) \oplus \big(\alpha \otimes \mathbb{I}_{n-r} - \Phi_{n-r}(C)\big),
	$$
	and from the singularity of $\alpha \otimes \mathbb{I}_r - \Phi_r(B),$ the singularity of $\alpha \otimes \mathbb{I}_n - \Phi_n(B \oplus C)$ follows. But then for any choice of $C \in K_{n-r}$ for which $B \oplus C$ is not unitarily equivalent to $A,$ the last statement contradicts part (b) of Definition \ref{def511}. 
	To prove the existence of such a $C$ we 
	proceed as follows.
	As explained above, without loss of generality $K_1$ contains a nonzero $v \in V$ as well as $0 \in V.$ In particular, by convexity, $w_t = \oplus_{i=1}^{n-r} t v \in  K_{n-r}$ for every $t \in [0,1]$.
	
	If $C \in K_{n-r}$ with the required properties does not exist, then for any $w \in K_{n-r},$ the direct sum $B \oplus w$ is unitarily equivalent to $A,$ hence all the $B \oplus w_t$ for $t \in [0,1]$ are unitarily equivalent. Now any functional $\varphi$ on $V$ with $\varphi(v)\neq 0$ gives rise to a nonzero diagonal matrix $\varphi_{n-r}(w_1) = \oplus_{i=1}^{n-r} \varphi(v)$ and thus yields a continuous family of complex matrices
	$$
	\varphi_n(B \oplus w_t) = 
	\begin{pmatrix}
		\varphi_r(B) & 0 \\
		0 & t\, \varphi_{n-r}(w_1)
	\end{pmatrix} \in \mathbb{M}_n,
	$$ 
	which are all unitarily equivalent. But this is a contradiction by a simple eigenvalue count.

	Note that the strict inequality $\Phi_r(B) \prec \alpha \otimes \mathbb{I}_r$ for all $B \in K_r$ does not hold  when $r \geq n$, since the singularity of $\alpha \otimes \mathbb{I}_n - \Phi_n(A)$ implies the singularity of $\alpha \otimes \mathbb{I}_{n+s} - \Phi_{n+s}(A\oplus C)$ for any $s \in \mathbb{N}$ and $C \in K_s$.
\end{opomba}

\begin{trditev}\label{prop33}
	Let $\textbf{$K$} = (K_n)_{n \in \mathbb{N}}$ be a matrix convex set. Then the matrix exposed points in $K_1$ coincide with the ordinary exposed points of $K_1$.
\end{trditev}

\begin{dokaz}
	First assume $A$ is a matrix exposed point in $K_1.$
	Conditions (a) and (b) in Definition \ref{def511} for the case $n=1$ imply the existence of a continuous linear functional $\varphi : V \to \mathbb{C}$ and a real number $\alpha$ such that:
	\begin{align}\label{13-1}
		\begin{split}
			\varphi(B) &= \varphi_1(B) \leq \alpha \ \text{ for all } B \in K_1,\\
			\varphi^{-1}(\alpha) &= \{A\}.
		\end{split}
	\end{align}
	So $A$ is ordinary exposed in $K_1$.
	
	For the converse, assume $A$ is an ordinary exposed point in $K_1$ and $\varphi : V \to \mathbb{C}$ a continuous linear functional satisfying (\ref{13-1}). We need to prove that for any positive integer $r$ and $B \in K_r$ the ampliation $\varphi_r$ satisfies $\varphi_r(B) \preceq \alpha \otimes \mathbb{I}_r = \alpha \mathbb{I}_r.$ Assume otherwise. Then there is an $r$ > $1$ and $B \in K_r,$ for which $\varphi_r(B) \npreceq \alpha \mathbb{I}_r,$ meaning that for some unit vector $x \in \mathbb{C}^r$ we have
	\begin{equation}\label{en1}
		x^\ast \big(\varphi_r(B) - \alpha \mathbb{I}_r\big) x =
		\big\langle \big(\varphi_r(B) - \alpha \mathbb{I}_r\big) x, x\big\rangle \in \mathbb{C}\backslash \{t \in \R \ | \ t \leq 0\}.
	\end{equation}
	Since $(\varphi_n)_{n \in \mathbb{N}}$ is matrix affine, we also have $x^\ast \varphi_r(B) x = \varphi(x^\ast B x) \leq \alpha,$ from which we deduce
	$$
	x^\ast \big(\varphi_r(B) - \alpha \mathbb{I}_r\big) x =
	x^\ast \varphi_r(B) x - \alpha \leq 0,
	$$
	contradicting (\ref{en1}). We conclude that the pair ($\varphi,$ $\alpha$) satisfies the conditions in Definition \ref{def511}, i.e.,~it matricially exposes $A$ in $\textbf{$K$}.$
\end{dokaz}

\subsection{Interplay between matrix extreme points and matrix exposed points} 
While in the classical theory it is straightforward to see that any exposed point is extreme, we need some preliminary observations (given in Proposition \ref{lema14}) to prove the matrix analogue of this claim in Proposition \ref{tr15}. Theorem \ref{tr111} then asserts when the converse holds: a matrix extreme point that is ordinary exposed is in fact matrix exposed.

\begin{trditev}\label{lema14}
	Let $A \in K_n$ be a matrix exposed point with an exposing pair \text{\textnormal{(}}$\Phi, \alpha$\text{\textnormal{)}}. Then the following statements hold.
	\begin{enumerate}[\text\textnormal{{(a)}}]
		\item For any nonzero $x = \sum_{j=1}^n x_j \otimes e_j \in \mathbb{C}^n \otimes \mathbb{C}^n$ in the kernel of $\alpha \otimes \mathbb{I}_n - \Phi_n(A),$ the components $x_1, \ldots, x_n$ form a basis of $\mathbb{C}^n.$ 
	\end{enumerate}
	\begin{enumerate}[\text\textnormal{{(b)}}]
		\item The kernel of $\alpha \otimes \mathbb{I}_n - \Phi_n(A)$ is one-dimensional.
	\end{enumerate}
\end{trditev}

\begin{dokaz}
	(a) Let $x = \sum_{j=1}^n x_j \otimes e_j \in \mathbb{C}^n \otimes \mathbb{C}^n$ be an arbitrary nonzero vector in the kernel of $\alpha \otimes \mathbb{I}_n - \Phi_n(A)$ and suppose that the matrix $(x_1, \ldots, x_n) \in \mathbb{M}_n$ is singular, say of rank $r < n$. Without loss of generality assume its rank is achieved in the first $r$ columns. Let $P \in  \mathbb{M}_{r, n}$ be the projection of $\mathbb{C}^n$ onto span$\{x_1, \ldots, x_r\},$ so that $(P \otimes \mathbb{I}_n)x = \sum_{j=1}^r x_j \otimes e_j.$ Then by denoting $x^{(r)} = \sum_{j=1}^r x_j \otimes e_j,$ we can express:
	\begin{align*}
		(\alpha \otimes \mathbb{I}_r- \Phi_r(P A P^\ast))x^{(r)}  &= (P \otimes \mathbb{I}_n)(\alpha \otimes \mathbb{I}_n - \Phi_n(A))(P^\ast \otimes \mathbb{I}_n)x^{(r)}\\ &=
		(P \otimes \mathbb{I}_n)(\alpha \otimes \mathbb{I}_n- \Phi_n(A))x = 0.
	\end{align*}
	We deduce that $\alpha \otimes \mathbb{I}_r - \Phi_r(P A P^\ast)$ is singular,  which by part (c)  of Remark \ref{op12} contradicts the fact that $A$ is matrix exposed.
	
	(b) Suppose there are two linearly independent vectors $x = \sum_{j=1}^n x_j \otimes e_j$ and $y = \sum_{j=1}^n y_j \otimes e_j$ in the kernel of $\alpha \otimes \mathbb{I}_n - \Phi_n(A).$ If $P \in \mathbb{M}_{r, n}$ is a projection of rank $r$ then $P A P^\ast \in K_r,$ and we have $\mathbb{I}_r \otimes \alpha - \Phi_r(P A P^\ast) = (P \otimes \mathbb{I}_n)(\alpha \otimes \mathbb{I}_n - \Phi_n(A))(P^\ast \otimes \mathbb{I}_n).$
	Assume $r < n.$ If there exists a complex number $\lambda$ for which the linear combination $\lambda x + y$ lies in the image of $P^\ast \otimes \mathbb{I}_n,$ i.e.~$\lambda x + y = (P^\ast \otimes \mathbb{I}_n)z$ for some $z \in \mathbb{C}^r \otimes \mathbb{C}^n,$ then the vector $z$ lies in the kernel of $\mathbb{I}_r \otimes \alpha - \Phi_r(P A P^\ast).$ But this contradicts $A$ being matrix exposed.
	
	To find such a $\lambda$ and projection $P$ consider the matrices $M = (x_1, \ldots, x_n)$ and $N = (y_1, \ldots, y_n).$ They are both invertible by part (a), whence
	$$
	\det(\lambda M + N) = \det(M^{-1}) \det(\lambda \mathbb{I}_n + M^{-1}N).
	$$
	Since the invertible matrix $M^{-1}N$ has a nonzero eigenvalue, there is a $\lambda$ for which $\lambda M + N$ is singular. Taking $P$ to be the projection onto span$\{\lambda x_1 + y_1, \ldots, \lambda x_n + y_n\}$ then finishes the proof. \qedhere 
\end{dokaz}

\begin{trditev}\label{tr15}
	Let $\textbf{$K$} = (K_n)_{n \in \mathbb{N}}$ be a matrix convex set. Then every matrix exposed point of \textbf{$K$} is matrix extreme. 
\end{trditev}

\begin{dokaz}
	Let $A \in K_n$ be matrix exposed and ($\Phi,$ $\alpha$) the corresponding exposing pair. Suppose we can express $A$ as a proper matrix convex combination:
	\begin{equation}\label{c1}
		A = \sum_{i=1}^k V_i^\ast A_i V_i
	\end{equation}
	for $k$-tuples $(A_i)_{i=1}^k$ and $(V_i)_{i=1}^k,$ where $A_i \in K_{n_i}$ and the matrices $V_i \in \mathbb{M}_{n_i,n}$ are surjective (implying $n_i \leq n$) with the property $\sum_{i=1}^k V_i^\ast V_i = \mathbb{I}_n$. By assumption we have $\mathbb{I}_{n_i} \otimes \alpha - \Phi_{n_i}(A_i) \succeq 0$ for $i=1, \ldots, k,$ and 
	\begin{align}\label{ie-1.3}
		\alpha \otimes \mathbb{I}_n - \Phi_n(A) &= \sum_{i=1}^k (V_i^\ast \otimes \mathbb{I}_n)\,\big(\mathbb{I}_{n_i} \otimes \alpha - \Phi_{n_i}(A_i)\big)\, (V_i\otimes \mathbb{I}_n).
	\end{align} 
	
	Suppose one of the $A_i$ (without loss of generality $A_1$) is not unitarily equivalent to $A$ and hence satisfies $\mathbb{I}_{n_1} \otimes \alpha - \Phi_{n_1}(A_1) \succ 0 .$ We will prove that this implies $V_1 = 0$.
	First notice that for $i=1, \ldots, k$ we have $(V_i^\ast \otimes \mathbb{I}_n)\,\big(\mathbb{I}_{n_i} \otimes \alpha - \Phi_{n_i}(A_i)\big)\, (V_i\otimes \mathbb{I}_n) \succeq 0.$ By Proposition \ref{lema14}, there is $x = \sum_{j=1}^n x_j \otimes e_j \in \mathbb{C}^n \otimes \mathbb{C}^n$ from the kernel of $\alpha \otimes \mathbb{I}_n - \Phi_n(A)$ such that its components  $x_1, \ldots, x_n$ span $\mathbb{C}^n.$ From (\ref{ie-1.3}) we can deduce (using the positive semidefiniteness of the summands on the right-hand side) that $x$ lies in the intersection of the kernels of $(V_i^\ast \otimes \mathbb{I}_n)\,\big(\mathbb{I}_{n_i} \otimes \alpha - \Phi_{n_i}(A_i)\big)\, (V_i\otimes \mathbb{I}_n)$ for $i=1, \ldots, k.$ In particular, we have $(V_1^\ast \otimes \mathbb{I}_n)\,\big(\mathbb{I}_{n_1} \otimes \alpha - \Phi_{n_1}(A_1)\big)\, (V_1\otimes \mathbb{I}_n)x = 0.$ Now the positive definiteness of the middle factor and injectivity of $(V_1^\ast \otimes \mathbb{I}_n)$ imply that $x$ lies in the kernel of $V_1 \otimes \mathbb{I}_n,$ i.e.,
	$$
	(V_1 \otimes \mathbb{I}_n)x = (V_1 \otimes \mathbb{I}_n) \bigg( \sum_{j=1}^n x_j \otimes e_j\bigg) = \sum_{j=1}^n V_1 x_j \otimes e_j = 0.
	$$
	So $V_1 x_j = 0$ for $j=1,\ldots, n$ and hence $V_1 = 0.$
\end{dokaz}

\begin{primer}\label{ex36}
	We give an example of a matrix convex set \textbf{$K$} with a matrix extreme point that is not matrix exposed. It is a free spectrahedrop (i.e.,~a coordinate projection of a free spectrahedron), whose ground level component $K_1$ is the set in Figure \ref{izpneext}. This is a convex set with an extreme point that is not exposed.
	\begin{figure}[h!]
		\begin{center}
			\begin{tikzpicture}[scale=2]
				\begin{scope}
					\clip (-2,-1) rectangle (2,2);
					\pgfsetstrokecolor{orange};
					\draw[domain=-0.9:1.2] plot({\x},{\x^3}) node [right, color=orange]{$x_2=x_1^3$};
					\draw (-2,1) -- (2,1);
					\draw (-1,-0.75) -- (-1,1.7);
					\filldraw[fill=yellow!40!white, domain=0:1] plot({\x},{\x^3}) -- (-1,1) -- (-1,0) -- (0,0);
				\end{scope}
				\draw[->] (-2,0) -- (2,0) node[right]{$x_1$};
				\draw[->] (0,-0.75) -- (0,1.8) node[above]{$x_2$};
				\fill[color=red] (0,0) circle (1pt)node[below right]{$0$};
				\fill (-1,0) circle (1pt)node[below left]{$-1$};
				\fill (1,0) circle (1pt)node[below]{$1$};
				\fill (0, 1) circle (1pt)node[above right]{$1$};
			\end{tikzpicture}
			\caption{The origin of the coordinate system is an extreme point, which is not exposed.}
			\label{izpneext}
		\end{center}
	\end{figure}
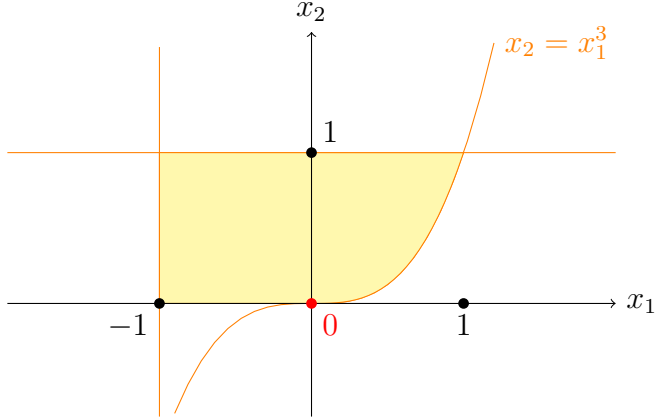
	
	As proved in \protect{\cite[Example 3.7]{NPS}}, the intersections of the depicted set with the first and second quadrant admit a so-called \textbf{exact Lasserre relaxation}, this being a sufficient condition for them and also their convex hull, i.e.,~their union, to be described by a spectrahedrop. So there is a linear pencil 
	$$L = C + \sum_{i=1}^2 A_i x_i + \sum_{j=1}^g B_j y_j$$ 
	in variables $(x_1, x_2, y_1, \ldots, y_g)$ such that
	$$
	K_1 = \{x \in \R^2 \ | \ \exists y \in \R^g \colon L(x, y) \succeq 0\}.
	$$
	Then
	\[
	\textbf{$K$} = \bigcup_{n \in \mathbb{N}}\{X \in \mathbb{S}^2_n \ | \ \exists Y \in  \mathbb{S}^g_n \colon L(X, Y) \succeq 0\}
	\]
	is a free spectrahedrop with a matrix extreme point that is not matrix exposed.
\end{primer}

\begin{primer}
	In analogy to the equality between extreme and exposed points of polyhedra in Euclidean spaces $\R^n,$ every matrix extreme point of a free spectrahedron is matrix exposed as showed in \cite[Corollary 6.21]{Kr} (cf.~Subsection \ref{sub331}).
\end{primer}

We now state a partial converse to Proposition \ref{tr15}, which will be proved after two technical lemmas.

\begin{izrek} \label{tr111}
	Let $\textbf{$K$} = (K_n)_{n \in \mathbb{N}}$ be a matrix convex set. Then:
	\begin{enumerate}[\text\textnormal{{(a)}}]
		\item Every matrix exposed point $A$ in $K_n$ is ordinary exposed in $K_n$.
	\end{enumerate}
	\begin{enumerate}[\text\textnormal{{(b)}}]\label{tr111bb}
		\item Every point $A$ in $K_n,$ which is both exposed and matrix extreme, is a matrix exposed point of $\textbf{$K$}$.
	\end{enumerate}
\end{izrek}

Denote by $S(\mathbb{M}_n) = \{p: \mathbb{M}_n \to \mathbb{C} \ | \ p \text{ unital positive}\}$ the state space of $\mathbb{M}_n$. 
Each positive functional $p \in S(\mathbb{M}_n)$ is of the form $p(\alpha) = \text{tr}\,(\gamma \alpha)$ for a fixed positive semidefinite matrix $\gamma \in \mathbb{M}_n$ with trace $1$ and so $S(\mathbb{M}_n)$ is a compact convex subset of $\mathbb{M}_n^*.$
The next lemma concerns the set we just introduced; it is a strict-positivity analogue of \cite[Lemma 5.2]{EW} along the lines of \protect{\cite[Lemma 2.16]{Kr}}.  In fact, \cite[Lemma 5.2]{EW} is a key result leading to the Effros-Winkler matricial Hahn-Banach separation theorem \cite{EW}, but the matricial separation there originates from the ability to separate a closed convex set from an outer point. On the other hand, to obtain the desired matricial separation in Theorem \ref{tr111}, we will implicitly use the idea that any (not necessarily closed) convex set can be separated from an outer point by a functional with values in an ordered extension field of $\R$ (cf.~\protect{\cite[Theorem 2.1]{NT}}). This, together with a finite intersection property \cite[Theorem 2.7.2]{BCR} motivates the introduction of real closed fields in the next lemma. For a real closed field $\mathcal{R}$  we will denote $>_{\mathcal{R}}$ the order relation on $\mathcal{R}.$ For $x, y >_{\mathcal{R}} 0$ we write $x \gg y$ if $x >_{\mathcal{R}} ny$ for all $n \in \mathbb{N}.$ For $x \in \mathcal{R}$ with an $n \in \mathbb{N}$ such that $-n <_{\mathcal{R}} x <_{\mathcal{R}} n$ we denote by st$(x) \in \R$ the standard part of $x$ (for more about real closed fields see \cite{BCR}).

\begin{lema}\label{lema1}
	Let $\mathcal{C}$ be a (convex) cone of continuous real affine functions on the state space $S(\mathbb{M}_n) \subseteq \mathbb{M}_n^\ast$ such that for every $f \in \mathcal{C}$ there is a state $p_f \in S(\mathbb{M}_n)$ with $f(p_f) > 0.$ Then there exists a real closed field $\mathcal{R}$ containing $\mathbb{R},$ and a unital positive $\R$-linear functional $p_0 : \mathbb{M}_n \to \mathcal{R}[\mathfrak{i}]$ satisfying $f_{\mathcal{R}}(p_0) >_{\mathcal{R}} 0$ for all $f \in \mathcal{C}\backslash \{0\}$ \text{\textnormal{(}}here $\mathcal{R}[\mathfrak{i}]$ stands for the algebraic closure of $\mathcal{R}$ and $f_{\mathcal{R}}$ is the unique extension of $f$ to an $\mathcal{R}$-linear map $M_n(\mathcal{R}[\mathfrak{i}])^\ast \rightarrow \mathcal{R}[\mathfrak{i}]$\text{\textnormal{)}}.
\end{lema}

\begin{dokaz}
	Denote the set of unital positive $\R$-linear functionals $p : \mathbb{M}_n \to \mathcal{R}[\mathfrak{i}]$ by $S_{\mathcal{R}}(\mathbb{M}_n)$ and
	for a given function $f \in \mathcal{C}\backslash \{0\}$ consider the following family 
	\begin{align*}\{f > 0\} := \{(\mathcal{R}, p) \ | \ \mathcal{R} \text{ real closed field over }\R,\ p \in S_{\mathcal{R}}(\mathbb{M}_n),\  f_{\mathcal{R}}(p) > 0\},
	\end{align*} 
	which is a type as we now explain (see \cite{Ho} as a reference for model theory). By assumption we have $\{f > 0\} \neq \emptyset$ for any $f \in \mathcal{C}\backslash \{0\}$ and we need to prove
	\begin{align*}
		\bigcap_{f \in \mathcal{C}\backslash \{0\}} \{f > 0\}\neq \emptyset.
	\end{align*}
	If we show that all finite intersections of the sets $\{f > 0\}$ are nonempty, then by a compactness argument as in \protect{\cite[Theorem 2.7.2]{BCR}}, there is a real closed field $\mathcal{R}$ over $\R,$ and a state $p \in S_{\mathcal{R}}(\mathbb{M}_n)$ such that for all $f \in \mathcal{C}\backslash \{0\}$ we have $f(p) > 0.$ So suppose there is an $n \in \mathbb{N}$ and functions $f_1, \ldots, f_n \in \mathcal{C}$ such that:
	$$
	\bigcap_{i=1}^n \{f_i > 0\} = \emptyset.
	$$
	
	Define the map $\theta: S(\mathbb{M}_n) \to \R^n$ with $\theta(p) = \big(f_1(p), \ldots, f_n(p)\big).$ This is clearly a continuous affine map, from which we see that $\theta(S(\mathbb{M}_n))$ is a compact convex subset of $\R^n.$ On the other hand we have by assumption that:
	$$
	\theta(S(\mathbb{M}_n)) \cap \R^n_+ = \emptyset,
	$$ 
	where $\R_+^n := [0,\infty)^n$. By a geometric version of the Hahn-Banach theorem (see e.g.~\protect{\cite[Section III.1]{Ba}}), the sets $\theta(S(\mathbb{M}_n))$ and $\R_+^n$ can be strictly separated, i.e.,~there is a linear function $g(x_1, \ldots, x_n) = c_1 x_1 + \cdots + c_n x_n$ on $\mathbb{R}^n$ and a real number $b$ such that 
	$g(y) \geq b$ for all $y \in \R_+^n$ and $g(z) < b$ for all $z \in \theta(S(\mathbb{M}_n))$. Moreover, since $\R_+^n$ is a cone, we have $b = 0.$

	For every standard unit vector $e_i \in \R^n$ we have $g(e_i) = c_i \geq 0.$ 
	So the function $f\circ \theta = c_1 f_1 + \cdots + c_n f_n$ is a conic combination of $f_1, \ldots, f_n$ and hence an element of $\mathcal{C}.$ But $f$ satisfies $\{p \in S(\mathbb{M}_n) \ | \ f(p)  > 0\} = \emptyset,$ which contradicts the assumption of the lemma.
\end{dokaz}

\begin{lema}\label{lema2}
	Let $\textbf{$K$} = (K_n)_{n \in \mathbb{N}}$ a matrix convex set for which $0 \in K_1.$ 
	Suppose there exist a linear functional $\varphi: M_n(V) \to \mathbb{C}$ and a real number $a > 0$ such that $\text{\textnormal{Re}}\,\varphi|_{K_n} < a.$ Then there is a state  $p : \mathbb{M}_n \to \mathbb{C}$ such that
	$$
	\text{\textnormal{Re}}\,\varphi(\alpha^\ast B \alpha) < a\, p(\alpha^\ast \alpha)
	$$
	for all $B \in K_r$, nonzero matrices $\alpha \in \mathbb{M}_{r, n}$ and positive integers $r$.
\end{lema}

\begin{dokaz}
	Let $\mathcal{C}$ be the set of all continuous affine functions on $S(\mathbb{M}_n)$ of the form:
	$$
	f_{v, \alpha}(p) = a\,p(\alpha^\ast \alpha) - \text{\textnormal{Re}}\,\varphi(\alpha^\ast v \alpha)
	$$
	for a matrix $\alpha \in \mathbb{M}_{r, n},$ $v \in K_r$ and $r \in \mathbb{N}.$ The set $\mathcal{C}$ is a cone as we can express
	\begin{align*}
		f_{v, \alpha} + f_{w, \beta} = f_{x, \gamma}, \quad c f_{v, \alpha} = f_{v, \sqrt{c}\alpha}
	\end{align*}
	for any real number $c \geq 0,$ where the matrix $\gamma$ is defined to be $\gamma^\ast = (\alpha \  \beta) \in \mathbb{M}_{r, 2n}$ and $x=v \oplus w.$ 
	
	Let us prove that for any $f = f_{v, \alpha} \in \mathcal{C}$ there is a state $p_f \in S(\mathbb{M}_n)$ with $f(p_f) > 0.$ Suppose $\alpha \neq 0$ and let $p_f$ be the state on $\mathbb{M}_n,$ for which $p_f(\alpha^\ast \alpha) = \|\alpha\|^2.$
	Then for the matrix $\beta = \frac{\alpha}{\|\alpha\|}$ the point $\beta^\ast v \beta$ lies in $K_n$ and by considering the assumption $\text{\textnormal{Re}}\,\varphi|_{K_n} < a,$ we have that
	$$
	\text{\textnormal{Re}}\,\varphi(\alpha^\ast v \alpha) = \|\alpha\|^2 \text{\textnormal{Re}}\,\varphi(\beta^\ast v \beta) < a\, p_f(\alpha^\ast \alpha).
	$$
	By Lemma \ref{lema1} there exist a real closed field $\mathcal{R}$ over $\R$ and a unital positive functional $p : \mathbb{M}_n \to \mathcal{R}[\mathfrak{i}]$ with $f_{v, \alpha}(p) >_{\mathcal{R}} 0$ for all $f_{v, \alpha} \in \mathcal{C}$ with $\alpha \neq 0.$
	
	Now using the procedure described in \protect{\cite[Corollary 2.17]{Kr}} we construct from $p$ a state on $\mathbb{M}_n$ with required properties. Since $p$ is positive, it is of the form $p(A) = \text{tr}\,(CA)$ for a positive semidefinite matrix $C \in M_n(\mathcal{R})$ with tr\,$C = 1.$ Let $D\in M_n(\mathcal{R})$ be its positive square root, i.e.,~$C = D^2$ with $D$ being positive semidefinite, and define $q \in S_{\mathcal{R}}(\mathbb{M}_n)$ by $q(A) = \text{tr}\,(DA)$. By \protect{\cite[Lemma 2.14]{Kr}} there exist $r \in \mathbb{N}$ and $\lambda_1 \gg...\gg\lambda_r>0$ in $\mathcal{R}$ along with $D_1, \ldots, D_r \in S_{\mathcal{R}}(\mathbb{M}_n)$ such that $D = \sum_{i=1}^r \lambda_j D_j$ (here $\lambda_1 D_1 = \text{st}\,D$ is the standard part of $D$). In the case $r=1,$ the state $q$ is already all we need; so assume $r \geq 2$ and let $q_j(A) = \text{tr}(D_j^2 A)$ for $j=1, \ldots, r.$ 
	
	Since $q(\mathbb{I}_n) = 1,$ we can take $\lambda_1 = 1.$ The goal is to prove we can replace the matrix $C,$ which defines $p,$ with $E^2 + D_r^2,$ where $E = \sum_{j=1}^{r-1} \lambda_j D_j.$ Note that for every $f \in \mathcal{C}$ we have st$f(p) = f(q_1) \geq 0$ and $f(q_r) \geq 0.$ Let $q_E(A) = \text{tr}(E^2 A)$ and suppose there was an $f \in \mathcal{C},$ for which $f(q_E + q_r) \leq_{\mathcal{R}} 0.$ Then $f(q_1) + f(q_r) = \text{st}\,f(q_E + q_r) \leq 0$ and hence $f(q_r) = 0.$ We deduce $f(q_E + q_r) = f(p) >_{\mathcal{R}} 0,$ which contradicts our assumption. Hence we have $f(q_E + q_r) >_{\mathcal{R}} 0$  for all $f \in \mathcal{C}.$
	
	Now continue the above procedure, i.e.,~in the next step replace the matrix $E^2 + D_r^2$ with $(\sum_{j=1}^{r-2} \lambda_j D_j)^2 + D_{r-1}^2 + D_r^2$ etc.~By induction we conclude that the state $q_F$ with the corresponding matrix $F=\sum_{j=1}^r D_r^2$ has the desired properties. 
\end{dokaz}

\begin{dokaz}[Proof of Theorem \ref{tr111}]
	
	(a) Let $A \in K_n$ be matrix exposed with the pair $\Phi:V \to \mathbb{M}_n$ and $\alpha \in \mathbb{M}_n$ exposing it as in Definition \ref{def511}. By Proposition \ref{lema14}, there is a vector $x = \sum_{j=1}^n x_j \otimes e_j \in \mathbb{C}^n \otimes \mathbb{C}^n$ spanning the kernel of $\alpha \otimes \mathbb{I}_n - \Phi_n(A)$ and whose components $x_1, \ldots, x_n$ are linearly independent. Define the functional $\varphi: M_n(V) \to \mathbb{C}$ by:
	$$
	\varphi(B) = x^\ast \Phi_n(B) x,
	$$
	and $a = x^\ast (\alpha \otimes \mathbb{I}_n) x \in \R.$
	Since the pair $(\Phi, \alpha)$ matricially exposes $A,$ we have
	$$
	a - \varphi(B) = x^\ast\big(\alpha \otimes \mathbb{I}_n - \Phi_n(B)\big) x \geq 0
	$$ 
	for all $B \in K_n$ and $a > \varphi(B)$ for each $B$ not in the unitary orbit of $A.$
	
	Let us check that $a > \varphi(B)$ also for any $B \in K_n\backslash\{A\},$ which is unitarily equivalent to $A.$ 
	Let $U \in \mathbb{M}_n$ be a unitary matrix such that $a = \varphi(U^\ast A U).$ Then
	$$
	\quad \quad \quad 0 = x^\ast\big(\alpha \otimes \mathbb{I}_n  - \Phi_n(U^\ast A U)\big) x = \big((U \otimes \mathbb{I}_n)x\big)^\ast \big(  \alpha \otimes \mathbb{I}_n - \Phi_n(A) \big)(U \otimes \mathbb{I}_n)x.
	$$
	The above together with the positive semidefiniteness of $\alpha \otimes \mathbb{I}_n - \Phi_n(A)$ implies that $(U \otimes \mathbb{I}_n)x$ lies in the one-dimensional kernel of $\alpha \otimes \mathbb{I}_n - \Phi_n(A)$ and is hence a multiple of $x.$ But then the unitary $U$ is a scalar multiple of the identity and so $U^\ast A U = A.$ We conclude that the pair $(\varphi, a)$ exposes $A$ in $K_n.$

	(b) We may assume $A \in K_n$ is a nonzero exposed and matrix extreme point and that $0 \in K_1$ as explained in the introductory section. 
	It is clear from the definition of a matrix extreme point that $A$ is not contained in the matrix convex set $\textbf{$L$}:=$ mconv$(K_n\backslash\{U^\ast A U \ | \ U \in \mathbb{U}_n\}).$ Since $A$ is nonzero, we have $0 \in L_1$ and hence \textbf{$L$} is closed under conjugation by contractions by Proposition \ref{pr15}.
	By assumption there is a continuous linear functional $\varphi: M_n(V) \to \mathbb{R}$ and real number $a$ with $\varphi(A) = a$ and $\varphi(B) < a$ for all $B \in K_n \backslash \{A\}.$ Since $0 \in K_1,$ we have $0 = \varphi(0) \leq a.$

	As in the proof of the Effros-Winkler matricial Hahn-Banach theorem in \protect{\cite[Theorem 5.4]{EW}} we divide our reasoning in three parts. First we gather together the key tools given by the previous lemmas, then we construct a candidate for the exposing pair. Finally we prove it does satisfy the desired separating conditions.
	
	By Lemma \ref{lema2}, there is a state $p : \mathbb{M}_n \to \mathbb{C}$ with: 
	\begin{equation} \label{eqstr}
		\text{Re}\,\varphi(\alpha^\ast B \alpha) < a\,p(\alpha^\ast \alpha)
	\end{equation}
	for any $B \in L_r,$ nonzero matrix  $\alpha \in \mathbb{M}_{r, n}$ and positive integer $r.$ Since $0 \in K_1,$ we have $ 0 < a\,p(\alpha^\ast \alpha)$ for any nonzero $\alpha \in \mathbb{M}_{r, n}$ so that $p$ is a faithful state.
	The condition $0 \in K_1$ also implies \textbf{$K$} is closed under conjugation by contractions, so we see by the properties of $\varphi$ that
	\begin{equation}\label{lastG}
		\text{Re}\,\varphi(\alpha^\ast B\alpha) \leq a\, p(\alpha^\ast \alpha)
	\end{equation}
	holds for any positive integer $r$, $B \in K_r$ and contraction $\alpha \in \mathbb{M}_{r, n}$ such that $p(\alpha^\ast \alpha) = 1.$

	By the GNS construction the map $p$ is determined by a representation $\pi : \mathbb{M}_n \to \mathcal{B}(\mathcal{H})$, where $\mathcal{H}$ is a finite-dimensional Hilbert space, and a cyclic and separating vector $x \in \mathcal{H}$ so that we have for arbitrary $\gamma \in \mathbb{M}_n$ the expression:
	$$
	p(\gamma) = \langle  \pi(\gamma)\,x, x \rangle.
	$$

	We now proceed to construct of the candidate  $\Phi: V \to \mathbb{M}_n$ for the map that  matricially  exposes $A.$ 
	To a row matrix $\alpha = [\alpha_1, \ldots, \alpha_n]\in \mathbb{M}_{1, n}$ assign the matrix $\tilde{\alpha} \in \mathbb{M}_n$, defined by
	$$
	\tilde{\alpha} = 
	\begin{pmatrix}
		\alpha_1 & \alpha_2 & \cdots & \alpha_n \\
		0 & 0 & \cdots & 0 \\
		\vdots  & \vdots  & \ddots & \vdots  \\
		0 & 0 & \cdots & 0
	\end{pmatrix}.
	$$
	Denote by $\widetilde{\mathbb{M}}_{1, n}$ the vector space consisting of matrices of such form and let $\mathcal{H}_0 := \pi(\widetilde{\mathbb{M}}_{1, n})\,x$. Since $x$ is separating, the space $\mathcal{H}_0$ is clearly an $n$-dimensional subspace of $\mathcal{H}$. We define on $\mathcal{H}_0$ a family of sesquilinear forms:
	$$
	\Psi_v(\pi(\tilde{\alpha})x, \pi(\tilde{\beta})x) =
	\varphi(\alpha^\ast v \beta)
	$$
	indexed by vectors $v$ from $V.$ For each $v \in V$ the form $\Psi_v$ is well-defined. Indeed, suppose we have $\pi(\tilde{\alpha_1})x = \pi(\tilde{\alpha_2})x$ for matrices
	$\alpha_1$ and $\alpha_2$ from $\mathbb{M}_{1, n}$. Linearity of $\pi$ gives $\pi(\tilde{\alpha_1}- \tilde{\alpha_1})x = 0$ and hence: 
	$$
	\big\langle \pi(\tilde{\alpha_1}- \tilde{\alpha_2})x, x\big\rangle = p(\alpha_1 - \alpha_2) = 0.
	$$
	Since $p$ is a faithful state, we have $\alpha_1 = \alpha_2,$ which implies that $\Psi_v$ is well-defined.
	
	Each sesquilinear form $\Psi_v$ on the finite-dimensional space $\mathcal{H}_0$ over $\mathbb{R}$ is uniquely determined by the linear map
	$\Phi(v) : \mathcal{H}_0 \to \mathcal{H}_0$ as follows:
	$$
	\varphi(\alpha^\ast v \beta) = \big\langle \Phi(v)\pi(\tilde{\alpha})x, \pi(\tilde{\beta})x \big\rangle.
	$$
	We thus get a map $\Phi : V \to \mathcal{B}(\mathcal{H}_0),$ which is both linear and weakly continuous. After choosing an orthonormal basis for $\mathcal{H}_0$, we identify $\mathcal{H}_0$ with $\mathbb{C}^n$ and the bounded operators $\mathcal{B}(\mathcal{H}_0)$ on it with $\mathbb{M}_n.$ 
	
	Now by letting $(e_i)_{i=1}^n$ be the standard basis of $\mathbb{C}^n$ and denoting $f_i := e_i^\ast$ for $i=1, \ldots, n$, we can express any matrix $B = (B_{i, j})_{i, j} \in M_n(V)$ (using the bimodule action on $V$) as a combination of the form:
	$$
	B = \sum_{i, j} e_i B_{i, j} f_j.
	$$
	Hence
	$$
	\varphi(B) = \sum_{i, j} \varphi(e_i B_{i, j} f_j) = 
	\sum_{i, j}\big\langle \Phi(B_{i, j})\,\pi(\tilde{f_j})x, \pi(\tilde{f_i})x \big\rangle = \big\langle \Phi_n(B) \eta_0, \eta_0\big\rangle,
	$$
	where
	$$
	\eta_0 = 
	\begin{pmatrix}
		\pi(\tilde{f_1})x \\
		\vdots \\
		\pi(\tilde{f_n})x
	\end{pmatrix}
	$$
	is a vector from $\mathcal{H}_0^n$ satisfying
	$$
	\|\eta_0\|^2 = \sum_{i=1}^n \|  \pi(\tilde{f_i})x \|^2 =
	\sum_{i=1}^n  p(f_i^\ast f_i) =  p\bigg(\sum_{i=1}^n f_i^\ast f_i\bigg) =  p(\mathbb{I}_n) = 1.
	$$
	
	In the last part of the proof we argue that $\Phi$ is the desired map. For this we need to check that for all $B \in K_r$ and $r \in \mathbb{N}$ the condition $ \Phi_r(B)\preceq a \mathbb{I}_{n} \otimes \mathbb{I}_{r} = a \mathbb{I}_{n \cdot r}$ holds or equivalently, 
	\begin{equation}\label{pog1}
		\text{Re}\,\big\langle \Phi_r(B) \eta, \eta\big\rangle = \big\langle \text{Re}\,\Phi_r(B) \eta, \eta\big\rangle \leq a\, \langle \eta, \eta \rangle
	\end{equation}
	holds for every vector $\eta \in (\mathbb{C}^n)^r.$ Since $x$ is cyclic, we can write any $\eta \in (\mathbb{C}^n)^r$ as
	\begin{align}\label{eta1}
		\eta = 
		\begin{pmatrix}
			\pi(\tilde{\alpha_1})x \\
			\vdots \\
			\pi(\tilde{\alpha_r})x
		\end{pmatrix},
	\end{align}
	where $\alpha_i \in \mathbb{M}_{1, n}$ for $i=1, \ldots, r$. In addition, we can express the norm of $\eta$ through the values of $p$ by
	\begin{equation}\label{eta}
		\|\eta\|^2 = \sum_{i=1}^n \|  \pi(\tilde{\alpha_i})x \|^2 = 
		\sum_{i=1}^n p(\alpha_i^\ast \alpha_i) = p(\alpha^\ast \alpha), \text{ where }
		\alpha = 
		\begin{pmatrix}
			\alpha_1 \\
			\vdots \\
			\alpha_r
		\end{pmatrix}\in \mathbb{M}_{r, n}.
	\end{equation}
	If $\eta$ is a unit vector, then for every $z \in \mathbb{C}^n$
	we have by the Cauchy-Schwarz inequality that:
	\begin{equation}\label{eq-110}
		\|\alpha z\|^2 = \sum_{i=1}^n |\alpha_i z|^2 \leq \|\eta \| \|z\| = 1
	\end{equation}
	showing $\alpha$ is a contraction. 
	
	We can finally check the validity of condition \eqref{pog1} for any unit vector $\eta \in (\mathbb{C}^n)^r$ using the property \eqref{lastG} of $\varphi$ and the just established connection \eqref{eta}:
	\begin{align*} 
		\big\langle \text{Re}\,\Phi_r(B) \eta, \eta\big\rangle &= 
		\sum_{i, j}\big\langle \text{Re}\, \Phi(B_{i, j})\, \pi(\tilde{\alpha_j})x,  \pi(\tilde{\alpha_i})x \big\rangle \\
		&= \sum_{i, j}\text{Re}\, \varphi(\alpha^\ast_i B_{i, j} \alpha_j) \\
		&= \text{Re}\,\varphi(\alpha^\ast B \alpha)\\
		&\leq a\, p(\alpha^\ast \alpha) \\
		&= a\,\|\eta\|^2 = a.
	\end{align*}
	Since $\varphi$ weakly separates the point $A$ from the set \textbf{$K$}, we have: 
	$$
	\big\langle\Phi_n(A) \eta_0, \eta_0\big\rangle
	= \varphi(A) = a,
	$$
	so the matrix $a\,\mathbb{I}_{n^2} - \Phi_n(A)$ is singular. Then by part (a) of Remark \ref{op12},  $a\,\mathbb{I}_{n^2} - \Phi_n(U^\ast A U)$ is singular for any $U \in \mathbb{U}_n.$
	To finish the proof we argue that for any $B \in K_n$ not in the unitary orbit of $A,$ the matrix $a\,\mathbb{I}_{n^2} - \Phi_r(B)$ is not singular. Indeed, if a unit vector $\eta \in \mathbb{C}^{n^2}$ of the form \eqref{eta1} satisfies $\big\langle\text{Re}\,\Phi_n(B) \eta, \eta\big\rangle = \text{Re}\,\varphi(\alpha^\ast B \alpha) =  a\,p(\alpha^\ast \alpha) = a,$ then $\alpha$ is a contraction by \eqref{eq-110}. Hence $\alpha^\ast B \alpha \in L_n$ as \textbf{$L$} is closed under conjugation of its elements by contractions.
	But then we have by the strong separation \eqref{eqstr} that $\text{Re}\,\varphi(\alpha^\ast B \alpha) <  a\,p(\alpha^\ast \alpha)$, which is a contradiction. We conclude that the pair $(\Phi, a\,\mathbb{I}_n)$ matricially exposes $A.$\qedhere
\end{dokaz}

Combining the techniques used in the proofs of Lemma \ref{lema2} and Theorem \ref{tr111} we obtain an Effros-Winkler type weak Hahn-Banach separation theorem for (not necessarily closed) matrix convex sets analogous to \protect{\cite[Corollary 2.17]{Kr}}.\looseness=-1

\begin{posledica}[Weak separation theorem for matrix convex sets]\label{th310}
	Let \textbf{$K$} be a matrix convex set in a dual space $V$ with $0 \in K_1$ and $A \notin K_n.$ Suppose there is a continuous linear functional $\varphi: M_n(V) \to \mathbb{C}$ and real number $a > 0$ such that $\text{\textnormal{Re}}\,\varphi|_{K_n} < a$ and $\varphi(A) = a.$ Then there exists a continuous linear map $\Phi : V \to \mathbb{M}_n$ such that
	$$
	\mathbb{I}_n \otimes \mathbb{I}_r  - \text{\textnormal{Re}}\, \Phi_r(B) \succ 0
	$$
	for every positive integer $r$ and $B \in K_r,$ but
	$$
	\text{\textnormal{ker}}\,\big(\mathbb{I}_n \otimes \mathbb{I}_n  - \text{\textnormal{Re}}\, \Phi_n(A)\big) \neq \{0\}.
	$$
\end{posledica}

\subsection{Straszewicz-Klee theorem for matrix convex sets}\label{subsec33}
A classical result on exposed points due to Straszewicz \cite[Section II.2]{Ba} states that the exposed points of a finite-dimensional compact convex set $K$ form a dense subset of the extreme points and hence their closed convex hull equals $K$.
This section extends the Straszewicz theorem, more precisely, its generalisation for normed spaces \cite{K} due to Klee, to the matrix convex setting.

We proceed by following the idea in \cite{HL} associating to a matrix convex set $\textbf{$K$} = (K_r)_{r \in \mathbb{N}}$ in the space $V$ a family of convex sets $\{\Gamma_n(\textbf{$K$})\}_{n \in \mathbb{N}}$ given by
\begin{equation}\label{eq311}
\Gamma_n(\textbf{$K$}) = \{ (\gamma^\ast \gamma, \gamma^\ast A \gamma) \ | \ \gamma \in \mathbb{M}_{k, n}, \text{tr}(\gamma^\ast \gamma) = 1, k \in \mathbb{N}, A \in K_k \} \subseteq \mathbb{M}_n \times M_n(V).
\end{equation}
The set $\Gamma_n(\textbf{$K$})$ is indeed convex as we can express
$$
t \gamma^\ast A \gamma + (1-t)\delta^\ast B \delta = 
\bigg(t^{1/2}\gamma^\ast \ \  (1-t)^{1/2}\delta^\ast\bigg) 
\begin{pmatrix}
A & 0\\
0 & B \\
\end{pmatrix}
\begin{pmatrix}
t^{1/2}\gamma\\
(1-t)^{1/2}\delta\\
\end{pmatrix}
$$
for elements $(\gamma^\ast \gamma, \gamma^\ast A \gamma)$ and $(\delta^\ast \delta, \delta^\ast B \delta)$ from $\Gamma_n(\textbf{$K$}),$ where $A\in K_r,$ $B \in K_s$ and $\gamma \in \mathbb{M}_{r, n},$ $\delta \in \mathbb{M}_{s, n}$ are matrices satisfying $\text{tr}(\gamma^\ast \gamma) = \text{tr}(\delta^\ast \delta) = 1,$ and arbitrary real number $t$ in $[0, 1].$
Since \textbf{$K$} is closed under direct sums, 
and we have
$$
\text{tr}\bigg(
\begin{pmatrix}
t^{1/2}\gamma\\
(1-t)^{1/2}\delta\\
\end{pmatrix}^\ast
\begin{pmatrix}
t^{1/2}\gamma\\
(1-t)^{1/2}\delta\\
\end{pmatrix}
\bigg)
 =
t\,\text{tr}(\gamma^\ast \gamma) + (1-t)\,\text{tr}(\delta^\ast \delta) = 1,
$$
the convex combination $t \big(\gamma^\ast \gamma, \gamma^\ast A \gamma\big) + (1-t)\big(\delta^\ast \delta, \delta^\ast B \delta\big)$ lies in $\Gamma_n(\textbf{$K$}).$

Additionally we can assume for any element $(\gamma^\ast \gamma, \gamma^\ast A \gamma)$ from $\Gamma_n(\textbf{$K$})$, where $A$ is in $K_r$, that the matrix $\gamma \in \mathbb{M}_{r, n}$ is surjective (and hence $r \leq n$) as we now explain. Let $\gamma \in \mathbb{M}_{r, n}$ be an arbitrary matrix of rank $s \in \mathbb{N}$ with the property tr$(\gamma^\ast \gamma) = 1.$  Let $\xi \in \mathbb{M}_{r, s}$ be an isometry from $\mathbb{C}^s$ to the range of $\gamma.$ Then
\begin{equation}\label{izo}
\gamma^\ast A \gamma = (\xi \gamma)^\ast (\xi A \xi^\ast) (\xi \gamma),
\end{equation}
where $\xi A \xi^\ast$ lies in $K_s$ and the matrix $\xi \gamma$ is surjective with tr$((\xi \gamma)^\ast (\xi \gamma)) = 1$.

The following variant of the Douglas Lemma \cite{D} will be often used to establish well-definedness of objects and maps in the remainder of this section.

\begin{lema}\label{lema18}
Let $\gamma \in \mathbb{M}_{r, n}$ and $\delta \in \mathbb{M}_{s, n}$ be surjective matrices. Then $\gamma^\ast \gamma = \delta^\ast \delta$ if and only if $r=s$ and there exists a unitary matrix $U \in \mathbb{M}_r$ such that $\gamma = U \delta.$
\end{lema}

We now state a correspondence between matrix exposed points of \textbf{$K$} at level $n$ and ordinary exposed points of $\Gamma_n(\textbf{$K$})$ (cf.~\cite[Proposition 2.14]{HL}).

\begin{trditev} \label{tr19}
	Let $\textbf{$K$} = (K_m)_{m \in \mathbb{N}}$ be a matrix convex set and $A \in K_r$. 
	\begin{enumerate}[\text\textnormal{{(a)}}]
		\item Let $\gamma \in \mathbb{M}_{r, n}$ be a surjective matrix with $\text{\textnormal{tr}}(\gamma^\ast \gamma)=1$ such that the point
		$(\gamma^\ast \gamma, \gamma^\ast A \gamma)$ is exposed in $\Gamma_n(\textbf{$K$}).$ Then $A$ is a matrix exposed point of \textbf{$K$}.
\end{enumerate}	
	\begin{enumerate}[\text\textnormal{{(b)}}]
		\item If $A$ is matrix exposed in \textbf{$K$}, then for any invertible $\gamma \in \mathbb{M}_r$ with $\text{\textnormal{tr}}(\gamma^\ast \gamma)=1,$ the point
		$(\gamma^\ast \gamma, \gamma^\ast A \gamma)$ is exposed in $\Gamma_r(\textbf{$K$}).$
	\end{enumerate}
\end{trditev}

\begin{dokaz} 
	To prove (a) suppose that $(\gamma^\ast \gamma, \gamma^\ast A \gamma)$ is an exposed point of $\Gamma_n(\textbf{$K$}),$ where $A\in K_r$ and $\gamma \in \mathbb{M}_{r,n}$ is a surjective matrix satisfying $\text{\textnormal{tr}}(\gamma^\ast \gamma)=1.$ By assumption there exists a continuous linear functional $\varphi: \mathbb{M}_n \times M_n(V) \to \mathbb{C}$ and a real number $a$ such that $\varphi(\gamma^\ast \gamma, \gamma^\ast A \gamma) = a$ and $\varphi(\delta, C) > a$ for all $(\delta, C) \in \Gamma_n(\textbf{$K$}) \backslash \{(\gamma^\ast \gamma, \gamma^\ast A \gamma)\}.$ 
	
	Note that $\varphi$ can be written as  $\varphi = \varphi_1 + \varphi_2,$ where $\varphi_1 : \mathbb{M}_n \to \mathbb{C}$ and $\varphi_2 : M_n(V) \to \mathbb{C}.$ By the Riesz representation theorem there is a matrix $\beta \in \mathbb{M}_n$ such that for all $\delta \in \mathbb{M}_n$ we have $\varphi_1(\delta) = \text{tr}(\beta \delta).$ It is easy to check that $\beta$ needs to be self-adjoint and that every matrix $\eta = (\eta_1, \ldots, \eta_n) \in \mathbb{M}_{m, n}$ and its corresponding vectorization 
	\begin{equation}\label{eq-111}
	v_{\eta}= \sum_{i=1}^n e^{\prime}_i \otimes \eta_i \in \mathbb{C}^n \otimes \mathbb{C}^m
	\end{equation} satisfy
	$$
	\varphi_1(\eta^\ast \eta) = \text{tr}(\beta \eta^\ast \eta) = \text{tr}(\eta \beta \eta^\ast) =v_{\eta}^\ast(\beta \otimes \mathbb{I}_m) v_{\eta}.
	$$
	By the canonical correspondence between linear functionals on $M_n(V)$ and linear maps $V \to \mathbb{M}_n,$ the functional $\varphi_2$ gives rise to a map $\Psi: V \to \mathbb{M}_n.$ For every $B \in M_m(V),$ matrix $\eta = (\eta_1, \ldots, \eta_n) \in \mathbb{M}_{m, n}$ and its vectorization 
	$
	v_{\eta}
	$
	as in \eqref{eq-111}
	we have
	$$
	\varphi_2(\eta^\ast B \eta) = e^*\Psi_n(\eta^\ast B \eta)e =  v_{\eta}^\ast \Psi_m(B) v_{\eta},
	$$
	where $e=\sum_{i=1}^n e^\prime_i \otimes e^\prime_i \in \mathbb{C}^n \otimes \mathbb{C}^n$ and $e^\prime_i$ are standard basis vectors of $\mathbb{C}^n.$
	
	So for every positive integer $m,$ element $B \in K_m$ and matrix $\eta \in \mathbb{M}_{m, n}$ with the property tr$(\eta^\ast \eta) = 1,$ 
	\begin{align}\label{al11}
	v_{\eta}^\ast \big( (\beta - a \mathbb{I}_n)  \otimes \mathbb{I}_m - \Psi_m(B)\big)  v_{\eta} = \varphi(\eta^\ast \eta, \eta^\ast B \eta) - a \geq 0.
	\end{align}
	Since tr$(\eta^\ast \eta) = 1$ if and only if $v_{\eta} \in \mathbb{C}^{mn}$ is a unit vector, the above implies that the matrix $(\beta - a \mathbb{I}_n) \otimes \mathbb{I}_m - \Psi_m(B)$ is positive semidefinite. We also have:
	\begin{align} \label{al12}
	v_{\gamma}^\ast \big( (\beta - a \mathbb{I}_n) \otimes \mathbb{I}_r - \Psi_r(A)\big)  v_{\gamma} = \varphi(\gamma^\ast \gamma, \gamma^\ast A \gamma) - a= 0,
	\end{align}
	which gives that $(\beta - a \mathbb{I}_n) \otimes \mathbb{I}_r - \Psi_r(A)$ is singular (as it is positive semi-definite).
	
	From Lemma \ref{lema18}, the points $(\eta^\ast \eta, \eta^\ast B \eta)$ and $(\gamma^\ast \gamma, \gamma^\ast A \gamma)$ from $\Gamma_n(\textbf{$K$})$ are equal if and only if there is a unitary matrix $U \in \mathbb{M}_r$ such that $\eta = U \gamma$ and $A = U^\ast B U.$ Hence:
	\begin{align}\label{al18}
	\varphi(\eta^\ast \eta, \eta^\ast B \eta) - a &= \varphi(\gamma^\ast \gamma, \gamma^\ast A \gamma) - a \nonumber\\ 
	&= v_{\gamma}^\ast \big( (\beta - a \mathbb{I}_n) \otimes \mathbb{I}_r - \Psi_r(A)\big)  v_{\gamma} \nonumber\\
	&= v_{\gamma}^\ast \big( (\beta - a \mathbb{I}_n) \otimes \mathbb{I}_r - \Psi_r(U^\ast B U)\big)  v_{\gamma}\\
	&= ((\mathbb{I}_n \otimes U)v_{\gamma})^\ast \big( (\beta - a \mathbb{I}_n) \otimes \mathbb{I}_r - \Psi_r( B )\big)  (\mathbb{I}_n \otimes U)v_{\gamma}\nonumber\\
	&=  v_{\eta}^\ast \big( (\beta - a \mathbb{I}_n) \otimes \mathbb{I}_r - \Psi_r(B)\big)  v_{\eta} \nonumber.
	\end{align}
	Using the properties of $\varphi$ we deduce that for $B \in K_m$ the matrix $(\beta - a \mathbb{I}_n) \otimes \mathbb{I}_m - \Psi_m(B)$  is singular if and only if $m = r$ and $B$ is unitarily equivalent to $A.$
	
	It remains to correct the target space of $\Psi$ and the size of $\beta$, i.e., we need a pair $(\Phi,\alpha)$ with $\Phi: V \to \mathbb{M}_r$ and $\alpha \in \mathbb{M}_r$ that matricially exposes $A.$
	Let $\delta \in \mathbb{M}_{r,n}$ be any surjective matrix such that the range of $\delta^\ast \otimes \mathbb{I}_r$ contains $v_\gamma$ (that is in the kernel of $(\beta - a \mathbb{I}_n) \otimes \mathbb{I}_r - \Psi_r(A)$). Then the compression $(\Phi, \alpha) = \delta (\Psi, \beta - a \mathbb{I}_n) \delta^\ast$ defines an exposing pair for $A;$ the positive semi-definiteness in all the points of \textbf{$K$} is clear, definiteness outside of the unitary orbit of $A$ follows from the injectivity of $\delta^*,$ and singularity at the unitary conjugates of $A$ holds by the choice of $\delta$ (cf.~equations \eqref{al11} -- \eqref{al18}). In fact, we can take $\delta = \gamma.$ To see that, write $v_\gamma = \sum_{i=1}^r \tilde{\gamma}_i \otimes e_i,$ where $\tilde{\gamma}_i$ is the $i$-th row of $\gamma$ and $e_i$ are standard basis vectors of $\mathbb{C}^r.$ Then for any $y = \sum_{i=1}^r y_i \otimes e_i \in \text{ker}(\gamma \otimes \mathbb{I}_r),$
	$$
	\langle v_\gamma, y \rangle = \sum_{i=1}^r \langle \tilde{\gamma}_i, y_i \rangle = 0.
	$$
	Hence, $v_\gamma$ lies in the orthogonal complement of ker\,$(\gamma \otimes \mathbb{I}_r)$, i.e., in the range of $\gamma^* \otimes \mathbb{I}_r.$

	To prove (b) assume that $A \in K_r$ is matrix exposed and $\Phi:V \to \mathbb{M}_r$ together with $\alpha \in \mathbb{M}_r$ are as in Definition \ref{def511}.
	Let $\gamma \in \mathbb{M}_{r}$ be an invertible matrix with $\text{\textnormal{tr}}(\gamma^\ast \gamma)=1.$ We claim that $(\gamma^\ast \gamma, \gamma^\ast A \gamma)$ is exposed in $\Gamma_r(\textbf{$K$}).$ 
	
	By Proposition \ref{lema14}, there is a vector $x = \sum_{j=1}^r x_j \otimes e_j \in \mathbb{C}^r \otimes \mathbb{C}^r$ that spans the kernel of $\alpha \otimes \mathbb{I}_r - \Phi_r(A)$ and whose components $x_1, \ldots, x_r$ are linearly independent.
	The pair $(\Phi, \alpha)$ produces another pair $(\Psi, \beta) = \gamma^* (\Phi, \alpha)\gamma$ and denoting by $\delta$ the inverse of $\gamma$, we can express 
	$$(\Phi, \alpha) = \delta^* (\Psi, \beta)\delta.$$ 
	Now define a functional $\varphi: \mathbb{M}_r \times M_r(V) \to \mathbb{C},$ which is for $C \in M_r(V)$ and $\mu \in \mathbb{M}_r$ given by
	$$
	\varphi(\mu, C) = y^*\big(\beta \otimes \mu - \Psi_r(C)\big)y,
	$$
	where $y = (\delta \otimes \delta)x \in \mathbb{C}^r \otimes \mathbb{C}^r.$ 
	The calculations in the above proof of (a) together with the fact that $\Psi$ is a matrix affine map show that for any positive integer $m,$ element $B \in K_m$ and surjective matrix $\eta \in \mathbb{M}_{m, r}$ with tr$(\eta^\ast \eta) = 1,$ 
	\begin{align*}
		\varphi(\eta^\ast \eta, \eta^\ast B \eta)& =\big((\mathbb{I}_r \otimes \eta)y\big)^\ast \big(\beta \otimes \mathbb{I}_m - \Psi_m(B)\big)  (\mathbb{I}_r \otimes \eta)y 
		\geq 0.
	\end{align*}
	Next we show that if $(\eta^\ast \eta, \eta^\ast B \eta)$ is different from $(\gamma^\ast \gamma, \gamma^\ast A \gamma),$ then $\varphi(\eta^\ast \eta, \eta^\ast B \eta) > 0.$
	If $B$ is not unitarily equivalent to $A,$ the positive definiteness of $\beta \otimes \mathbb{I}_m - \Psi_m(B)$ (following from the positive definiteness of $\alpha \otimes \mathbb{I}_m - \Phi_m(B)$ and the invertibility of $\gamma$) implies $(\mathbb{I}_r \otimes \eta)y=  0$ or $\varphi(\eta^\ast \eta, \eta^\ast B \eta) > 0.$ But $(\mathbb{I}_r \otimes \eta)y = (\delta \otimes \eta \delta)x \neq 0$ for nonzero $\eta$ since span$\{x_1, \ldots, x_r\} = \mathbb{C}^r$ and $\delta$ is invertible.
	On the other hand, if $B = U^* A U$ for some unitary $U \in \mathbb{M}_r$ and  $\varphi(\eta^\ast \eta, \eta^\ast B \eta) =0,$ then 
	\begin{align*}
		0 = \varphi(\eta^\ast \eta, \eta^\ast B \eta) &= ((\mathbb{I}_r \otimes \eta)y)^\ast \big(\beta \otimes \mathbb{I}_r  - \Psi_r(B)\big)  (\mathbb{I}_r \otimes \eta)y \\
		&= ((\mathbb{I}_r \otimes U \eta)y)^\ast \big(\beta \otimes \mathbb{I}_r  - \Psi_r(A)\big)  (\mathbb{I}_r \otimes U \eta)y \\
		&= ((\gamma \otimes U \eta)y)^\ast \big(\alpha \otimes \mathbb{I}_r  - \Phi_r(A)\big)  (\gamma \otimes U \eta)y \\
		&=((\mathbb{I}_r \otimes U \eta \delta)x)^\ast \big(\alpha \otimes \mathbb{I}_r  - \Phi_r(A)\big)  (\mathbb{I}_r \otimes U \eta\delta)x.
	\end{align*} 
	By part (b) of Proposition \ref{lema14}, the kernel of $\alpha \otimes \mathbb{I}_r  - \Phi_r(A)$ is spanned by $x$, hence $(\mathbb{I}_r \otimes U \eta\delta)x = \lambda x$
	for some nonzero $\lambda \in \mathbb{C}.$ Since the components of $x$ form a basis of $\mathbb{C}^n,$ this implies $U \eta\delta = \lambda \mathbb{I}_r,$ i.e., $U \eta = \lambda\gamma.$ A comparison of Hilbert-Schmidt norms now yields $|\lambda|=1,$ which in turn implies $(\eta^\ast \eta, \eta^\ast B \eta) = (\gamma^\ast \gamma, \gamma^\ast A \gamma),$ a contradiction.

	It remains to observe that
	\begin{align*}
	\varphi(\gamma^\ast \gamma, \gamma^\ast A \gamma) &= ((\mathbb{I}_r \otimes \gamma)y)^\ast \big(\beta \otimes \mathbb{I}_r  - \Psi_r(A)\big)  (\mathbb{I}_r \otimes \gamma)y \\
	&= ((\delta \otimes \mathbb{I}_r)x)^\ast \big(\beta \otimes \mathbb{I}_r  - \Psi_r(A)\big)  (\delta \otimes \mathbb{I}_r)x \\
	&=x^*\big(\delta^* \beta \delta \otimes \mathbb{I}_r  - (\delta^*\Psi_r\delta)(A) \big)x \\
	&= x^*\big(\alpha \otimes \mathbb{I}_r  - \Phi_r(A)\big)x=0,
	\end{align*}
	which proves that the pair $(-\varphi, 0)$ exposes the point $(\gamma^\ast \gamma, \gamma^\ast A \gamma)$ in $\Gamma_n(\textbf{$K$}).$
\end{dokaz}

Using the results established above we can prove, applying similar techniques as in \protect{\cite[Theorem 4.3]{WW}}, a version of the Krein-Milman matricial theorem featuring matrix exposed points.

\begin{izrek}[\textbf{The Straszewicz-Klee theorem for matrix convex sets}]\label{th213}
	Let \textbf{$K$} be a compact matrix convex set in a normed vector space $V.$ Then $\text{\textnormal{mexp}}\,\textbf{$K$} \neq \emptyset$ and
	$$
	\textbf{$K$} = \overline{\text{\textnormal{mconv}}}\,(\text{\textnormal{mexp}}\,\textbf{$K$}).
	$$
\end{izrek}

\begin{dokaz}
	
Let $\textbf{$K$}=(K_n)_{n \in \mathbb{N}}$ be a compact matrix convex set in a normed space $V$. By the classical Straszewicz-Klee theorem (see \cite[Section II.2]{Ba} and \cite{K}), the set of exposed points of $K_1$ is non-empty, so by Proposition \ref{prop33}, also $\text{\textnormal{mexp}}\,\textbf{$K$} \neq \emptyset.$ The inclusion $\overline{\text{\textnormal{mconv}}}(\text{\textnormal{mexp}}\textbf{$K$}) \subseteq \textbf{$K$}$ holds by definition, so we only need to prove $\textbf{$K$} \subseteq \overline{\text{\textnormal{mconv}}}(\text{\textnormal{mexp}}\textbf{$K$})$. We can assume $0 \in \overline{\text{\textnormal{mconv}}}(\text{\textnormal{mexp}}\textbf{$K$}),$ otherwise translate each $K_n$ by an element $a \otimes \mathbb{I}_n$ for some fixed $a \in K_1$ (as translations preserve matrix exposed points).

Suppose there is $A$ in $K_n\backslash(\overline{\text{\textnormal{mconv}}}(\text{\textnormal{mexp}}\textbf{$K$}))_n.$ By the matricial Hahn-Banach separation theorem \protect{\cite[Theorem 5.4]{EW}} there exists a continuous linear map $\Phi : V \to \mathbb{M}_n$ such that:\looseness=-1
\begin{equation}\label{pr}
	\text{\textnormal{Re}}\, \Phi_r(B) \preceq \mathbb{I}_{n\cdot r}
\end{equation}
for every positive integer $r$ and $B \in (\overline{\text{\textnormal{mconv}}}(\text{\textnormal{mexp}}\textbf{$K$}))_r,$ but: 
\begin{equation}\label{dr}
	\text{\textnormal{Re}}\,\Phi_n(A) \npreceq \mathbb{I}_{n^2}.
\end{equation}

The idea in the following is to reduce the problem to a point where the classical  Straszewicz theorem \cite[Section II.2]{Ba}, more precisely, its generalisation for normed spaces due to Klee \cite{K}, applies. For that we use the fact that the map $\Phi$ associates to every pair of matrices $\gamma, \delta \in \mathbb{M}_{r, n}$
a continuous linear functional $\varphi_2: M_n(V) \to \mathbb{C}$ with:\looseness=-1
$$
\varphi_2(\gamma^\ast B \delta) = v_{\gamma}^\ast\Phi_r(B)v_{\delta}
$$
for $B \in M_r(V)$. Here $v_{\gamma}$ and $v_{\delta}$ from $(\mathbb{C}^n)^r$ are the vectorizations of $\gamma$ and $\delta$ as in \eqref{eq-111}.

Let $(\delta, C)$ be an exposed point of $\Gamma_n(\textbf{$K$})$. By Proposition \ref{tr19} it can be expressed as $(\delta, C) = (\eta^\ast \eta, \eta^\ast B \eta)$ for some matrix exposed point $B \in K_r,$ matrix $\eta \in \mathbb{M}_{r, n}$ with tr$(\eta^\ast \eta)=1$ (hence its vectorization $v_{\eta}$ satisfies $v_{\eta}^\ast v_{\eta} = 1$) and positive integer $r \leq n.$ Accordingly we extend $\varphi_2$ to a functional $\varphi: \mathbb{M}_n \times M_n(V) \to \mathbb{C}$ by
$$
\varphi(\eta^\ast \eta, \eta^\ast B \eta) = v_{\eta}^\ast (\mathbb{I}_{n\cdot r} + \Phi_r(B)) v_{\eta}.
$$ 
Inequality \eqref{pr} now gives
\begin{align*}
	\text{\textnormal{Re}}\,\varphi(\eta^\ast \eta, \eta^\ast B \eta) &= 
	\text{\textnormal{Re}}\,\big(v_{\eta}^\ast (\mathbb{I}_{n\cdot r} + \Phi_r(B)) v_{\eta} \big)\\
	&= 1 + v_{\eta}^\ast \,\text{\textnormal{Re}}\,\Phi_r(B) v_{\eta}\\
	&\leq 1 + v_{\eta}^\ast \,\mathbb{I}_{n\cdot r} v_{\eta} = 2.
\end{align*}
Since the above holds for every $(\delta, C)$ in $\Gamma_n(\textbf{$K$})$ and the compactness of $\textbf{$K$}$ implies the compactness of $\Gamma_n(\textbf{$K$})$, we have by the Klee generalisation of the Straszewicz theorem \cite{K} that
\begin{equation}\label{eq-117}
	\text{\textnormal{Re}}\,\varphi(\epsilon, D) \leq 2
\end{equation}
for every $(\epsilon, D)  \in \Gamma_n(\textbf{$K$})$ (since $\varphi$ is linear and every such $(\epsilon, D)$ is a limit of a sequence of convex combinations of exposed points). Equation \eqref{eq-117} implies that for every positive integer $r,$ unit vector $v_{\eta} \in \mathbb{C}^{r\cdot n}$ and $B \in K_r,$
$$
1 + v_{\eta}^\ast \,\text{\textnormal{Re}}\,\Phi_r(B) v_{\eta} = \text{\textnormal{Re}}\,\varphi(\eta^\ast \eta, \eta^\ast B \eta) \leq 2 = 1 + v_{\eta}^\ast \,\mathbb{I}_{n\cdot r} v_{\eta},
$$
which in turn implies $\text{\textnormal{Re}}\,\Phi_r(B) \preceq \mathbb{I}_{n\cdot r}$. But this contradicts condition \eqref{dr}, hence $\overline{\text{\textnormal{mconv}}}(\text{\textnormal{mext}}\textbf{$K$}) = \textbf{$K$}.$
\end{dokaz}

\subsection{Matrix exposed points of matrix state spaces}\label{subsec34}

For an operator system $\mathcal{R}$ let $\text{CP}_n(\mathcal{R})$ denote the set of all completely positive (cp) maps and $\text{UCP}_n(\mathcal{R})$ the set of all unital completely positive (ucp) maps from $\mathcal{R}$ to $\mathbb{M}_n$ respectively. We identify both sets in a canonical way with subsets of $M_n(\mathcal{R}^\ast).$ Then the family $\textbf{CP}(\mathcal{R})= (\text{CP}_n(\mathcal{R}))_{n \in \mathbb{N}}$ is a weak$^\ast$ closed matrix convex cone in $\mathcal{R}^\ast$ meaning it is closed under formation of the following combinations:
\begin{equation}\label{eq12}
\sum_{i=1}^{k} V_i^\ast \psi_i V_i   \in M_n(\mathcal{R}^\ast)
\end{equation}
for $k \in \mathbb{N},$ elements $\psi_i \in\text{CP}_{n_i}(\mathcal{R})$ and matrices $V_i \in \mathbb{M}_{n_i, n}.$
On the other hand, the family $\textbf{\text{UCP}}(\mathcal{R})= (\text{UCP}_n(\mathcal{R}))_{n \in \mathbb{N}},$ usually referred to as the \textbf{matrix state space} of the operator system $\mathcal{R},$ is a weak$^\ast$ compact matrix convex set in $\mathcal{R}^\ast,$ i.e.,~closed under matrix convex combinations \eqref{eq12}, where $\sum_{i=1}^{k} V_i^\ast V_i = \mathbb{I}_n$.
Moreover, by \protect{\cite[Proposition 1.2]{F}} we see that $\textbf{CP}(\mathcal{R})$ is the matrix conic hull of $\textbf{\text{UCP}}(\mathcal{R})$ as we have 
$$
\text{CP}_n(\mathcal{R}) = \{\gamma^\ast \varphi \gamma \ | \ \gamma \in \mathbb{M}_n,\  \varphi \in \text{UCP}_n(\mathcal{R})\}.
$$

For every positive integer $n$ the component $\text{CP}_n(\mathcal{R})$ is a convex cone and the set $\text{UCP}_n(\mathcal{R})$ is convex. 
In this section we describe the relationship between matrix exposed points of $\textbf{\text{UCP}}(\mathcal{R})$ and the exposed rays of $\textbf{CP}(\mathcal{R})$ while keeping in mind the analogous connection between matrix extreme points and extremal rays (see \cite{F}). We begin by recalling some of the known results.

A cp map $\varphi : \mathcal{R} \to \mathbb{M}_n$ is called \textbf{pure} if for every cp map $\psi : \mathcal{R} \to \mathbb{M}_n,$ for which $\varphi - \psi$ is completely positive, we have that $\psi = t \varphi$ for some real number $t \in [0, 1]$. 
It turns out (see \protect{\cite[Section II.8]{Ba}}) that pure cp maps in $\text{CP}_n(\mathcal{R})$ determine the extremal rays of this cone, i.e.,~the extremal rays are exactly rays of the form $\{t \varphi \ | \ t \geq 0\}$ for some pure cp map $\varphi.$ The following proposition explains the interplay between pure cp maps and pure states, i.e.,~pure ucp maps.

\begin{trditev}[\protect{\cite[Theorem 2.2]{F}}]\label{tr115}
	Let $\psi \in\text{\textnormal{CP}}_n(\mathcal{R})$ be pure. Then there exists a positive integer $k \leq n$, a pure state $\varphi \in \text{\textnormal{UCP}}_k(\mathcal{R})$ and a matrix $\gamma \in \mathbb{M}_{k, n}$ such that $\psi = \gamma^\ast \varphi \gamma.$
\end{trditev}

If the matrix $\gamma$ in this proposition is invertible, then the reverse implication holds, i.e., for a pure state $\varphi \in \text{UCP}_k(\mathcal{R})$ the map $\psi = \gamma^\ast \varphi \gamma \in\text{CP}_n(\mathcal{R})$ is pure ($\psi$ is clearly completely positive even with no condition on $\gamma$). Indeed, let $\theta \in\text{CP}_n(\mathcal{R})$ satisfy $\psi - \theta \in\text{CP}_n(\mathcal{R}).$ Then the map $ (\gamma^\ast)^{-1}\theta\gamma^{-1}$ lies in CP$_k(\mathcal{R})$ and: 
$$
(\gamma^\ast)^{-1}(\psi - \theta)\gamma^{-1} = \varphi - (\gamma^\ast)^{-1}\theta\gamma^{-1}  \in\text{CP}_k(\mathcal{R}).
$$
Since $\varphi$ is pure, there is a $t \in [0, 1]$ such that $(\gamma^\ast)^{-1}\theta\gamma^{-1} = t \varphi$ and hence $\theta = t \psi.$

To conclude this short summary we note the fact that the pure states on the operator system $\mathcal{R}$ correspond to matrix extreme points in $\textbf{\text{UCP}}(\mathcal{R})$ (see \cite{F}). We now turn our attention to the case of matrix exposed points.

\begin{trditev}\label{tr116}
	Suppose for a completely positive map $\psi \in \text{\textnormal{CP}}_n(\mathcal{R})$ the corresponding ray $\{t\psi \ | \ t \geq 0\}$ is exposed in $\textnormal{CP}_n(\mathcal{R}).$ Then there exists a positive integer $k \leq n$, a state $\varphi \in \text{\textnormal{UCP}}_k(\mathcal{R})$ and a matrix $\gamma \in \mathbb{M}_{k, n}$ such that $\psi = \gamma^\ast \varphi \gamma,$ and $\varphi$ is a matrix exposed point in  $\text{\textnormal{\textbf{UCP}}}(\mathcal{R}).$ 
\end{trditev}

\begin{dokaz}
	Since $\psi \in \text{\textnormal{CP}}_n(\mathcal{R})$ determines the exposed ray $\{t\psi \ | \ t \geq 0\}$ in $\textnormal{CP}_n(\mathcal{R}),$ there is a linear functional $F: M_n(\mathcal{R}^\ast) \to \mathbb{C}$ satisfying $F(t \psi) = 0$ for all $t \geq 0$ and $F(\theta) > 0$ for all $\theta \in\text{CP}_n(\mathcal{R}) \backslash \{t\psi \ | \ t \geq 0\}.$ By Proposition \ref{tr115}, there is a positive integer $k \leq n$, a pure state $\varphi \in \text{UCP}_k(\mathcal{R})$ and $\gamma \in \mathbb{M}_{k, n}$ such that $\psi = \gamma^\ast \varphi \gamma;$ so $\varphi$ is matrix extreme in $\textbf{\text{UCP}}(\mathcal{R}).$ If we prove that $\varphi$ is an exposed point in $\text{UCP}_k(\mathcal{R}),$ the by part (b) of Theorem \ref{tr111} it is also matrix exposed. 
	To prove $\varphi$ is exposed define the functional $G: M_k(\mathcal{R}^\ast) \to \mathbb{C}$ by
	$$
	G(\theta) = F(\gamma^\ast \theta \gamma)
	$$
	for $\theta \in M_k(\mathcal{R}^\ast).$ By the properties of $F$ we have $G(\theta) \geq 0$ for all $\theta \in \text{UCP}_k(\mathcal{R})$ (since $\gamma^\ast \theta \gamma \in \text{CP}_n(\mathcal{R})$) and $G(\varphi) = 0.$ As $\varphi$ is the only unital element of the ray $\{t\psi \ | \ t \geq 0\},$ we have $G(\theta) > 0$ for all $\theta \in  \text{UCP}_k(\mathcal{R}) \backslash \{\varphi\}.$  Hence the pair $(G, 0)$ exposes $\varphi.$
\end{dokaz}

\begin{opomba}
	To give a partial converse to Proposition \ref{tr116}, analogous to that of Proposition \ref{tr115}, let $\gamma \in \mathbb{M}_n$ in the above notation be invertible and let $\varphi \in \text{\textnormal{UCP}}_n(\mathcal{R})$ be matrix exposed. Then $\varphi$ is both matrix extreme and exposed in the classical sense; so the remark after Proposition \ref{tr115} implies that the map $\gamma^\ast \varphi \gamma$ is pure completely positive. Suppose we are in the special case when the pair $(G, 0)$ for some functional  $G: M_n(\mathcal{R}^\ast) \to \mathbb{C}$  exposes $\varphi$ in $\text{\textnormal{UCP}}_n(\mathcal{R}).$ We now claim that
	$$
	F(\theta) = G\big((\gamma^\ast)^{-1} \theta \gamma^{-1}\big)
	$$
	defines a functional $F: M_n(\mathcal{R}^\ast) \to \mathbb{C}$ which exposes the ray $\{t\, \gamma^\ast \varphi \gamma\ | \ t \geq 0\}.$ Clearly, $F(t\, \gamma^\ast \varphi \gamma) = t\, G(\varphi) = 0.$ For $\theta \in \text{\textnormal{CP}}_n(\mathcal{R})$ denote $t_{\theta} = \|\theta(1)\|$ so that we have $F\big(\frac{\theta}{t_{\theta}}\big) = G\big((\gamma^\ast)^{-1} \frac{\theta}{t_{\theta}} \gamma^{-1}\big) \geq 0$  (since $\frac{\theta}{t_{\theta}} \in \text{\textnormal{UCP}}_n(\mathcal{R}$)) implying $F(\theta) \geq 0.$ On the other hand, we have
	$F\big(\frac{\theta}{t_{\theta}}\big)  = G\big((\gamma^\ast)^{-1} \frac{\theta}{t_{\theta}} \gamma^{-1}\big) = 0$  (so $F(\theta) = 0$) if and only if $(\gamma^\ast)^{-1} \theta \gamma^{-1} = t_{\theta}\, \varphi,$ whence $\theta = t_{\theta} \, \gamma^\ast \varphi \gamma.$

\end{opomba}

\begin{primer}\label{ex:cstarexposed} Let $\mathcal{A}$ be a separable unital $C^\ast$-algebra. Then the matrix extreme points, i.e., pure matrix states, of the matrix state space $\textbf{\text{UCP}}(\mathcal{A})$ are matrix exposed. Indeed, by part (b) of Theorem \ref{tr111}, it is enough see that every ordinary extreme point $\varphi$ of $\text{CP}_n(\mathcal{A})$ is ordinary exposed. By the canonical correspondence between linear maps $\mathcal{A} \to \mathbb{M}_n$ and linear functionals on $M_n(\mathcal{A})$ (see \cite[Chapter 6]{Pa}), such a $\varphi$ determines a state $\tilde{\varphi}: M_n(\mathcal{A}) \to \mathbb{C}$ given by
	$$
	\tilde{\varphi}(X) = \frac{1}{n} \langle \varphi_n(X)e, e \rangle,
	$$
	where $e= e_1 \oplus \cdots \oplus e_n$ and the $e_i$ are standard basis vectors of $\mathbb{C}^n.$ The above correspondence is in fact a linear bijection and as such it preserves extreme and exposed points. Hence, $\varphi$ being extreme in $\text{CP}_n(\mathcal{A})$ implies that $\tilde{\varphi}$ is a pure state of the $C^*$-algebra $M_n(\mathcal{A}).$ But then by \cite[Corollary 3.55]{AS}, $\tilde{\varphi}$ is exposed in $M_n(\mathcal{A})^*$, which in turn implies that $\varphi$ is an exposed point of $\text{CP}_n(\mathcal{A}).$ \qedhere

\end{primer}

\section{Matrix faces and matrix exposed faces}\label{sec3}

We proceed by discussing several possible notions of a face and an exposed face of a matrix convex set. The main distinction is whether one considers subsets of a single component $K_n$ for some $n \in \mathbb{N}$ or multicomponent subsets of a matrix convex set \textbf{$K$}.
All the presented definitions aim to extend the concepts of a matrix extreme point or a matrix exposed point. We also explain how their properties resemble those of the (exposed) faces in the classical sense and investigate the interplay between the notions of a matrix face and a matrix exposed face. 

\subsection{Fixed-level matrix faces} In this section we present three aspirant definitions of a non-graded matrix face and explore the resemblance of their properties with the classical theory:
Proposition \ref{tr34} states that the $C^\ast$-extreme points of a matrix face of a matrix convex set \textbf{$K$} are matrix extreme in \textbf{$K$} and in Theorem \ref{tr211} we prove that any matrix face that is ordinary exposed is in fact a matrix exposed face. As a corollary of the latter we observe that every \ti or \tii of a free spectrahedron is matrix exposed. Although it is not clear whether they are ordinary faces, \tiii are included in the list as the most straightforward generalisation of a matrix extreme point.

\begin{definicija}\label{def41}
	Let $\textbf{$K$} = (K_r)_{r \in \mathbb{N}}$ be a matrix convex set in the space $V$ and $F$ a convex subset of $K_n$.
	
	(a) Then $F$ is a \textbf{\ti} if for every tuple of points $A_1, \ldots, A_k$ from \textbf{$K$} with $A_i \in K_{n_i}$ and every tuple of surjective matrices $\gamma_i \in \mathbb{M}_{n_i, n}$ satisfying $\sum_{i=1}^k \gamma_i^\ast \gamma_i = \mathbb{I}_n,$ the condition
	\begin{equation} \label{eq21}
		\sum_{i=1}^k \gamma_i^\ast A_i \gamma_i \in F,
	\end{equation}
	implies $n_i=n$ and $A_i \in F$ for $i=1, \ldots, k.$
	
	(b) If $F$ is a $C^\ast$-convex \ti, then it is a \textbf{\tii}.
	
	(c) The set $F$ is a \textbf{\tiii} if for every tuple of points $A_1, \ldots, A_k$ from \textbf{$K$} with $A_i \in K_{n_i}$ and every tuple of surjective matrices $\gamma_i \in \mathbb{M}_{n_i, n}$ satisfying $\sum_{i=1}^k \gamma_i^\ast \gamma_i = \mathbb{I}_n,$ the condition
	\begin{equation*}
		\sum_{i=1}^k \gamma_i^\ast A_i \gamma_i \in F,
	\end{equation*}
	implies $n_i=n$ and each $A_i$ is unitarily equivalent to some element in $F.$ We will denote by $\mathcal{U}(F) = \{U^\ast A U \ | \ A \in F, \ U \in \mathbb{M}_n \text{ unitary}\}$
	the unitary orbit of $F.$
\end{definicija}

\begin{opomba}\label{op22}
		(a) It is clear from the definition that a \ti or \tii $F$ is itself a face in the classical sense and that the matrix convex set  mconv$(K_n\backslash F)$ is disjoint from $F.$ Similarly, in the case of a \tiii, the set mconv$(K_n\backslash \mathcal{U}(F))$ is disjoint from $F.$
		
		(b) Since in dimension $1$ unitary equivalence implies equality and $C^\ast$-convexity implies classical convexity, the matrix faces of all types in $K_1$ coincide with its faces in the classical sense. Also, $K_1$ itself is a matrix face, however $K_n$ for $n > 1$ is never a matrix face of any type. Indeed, for any $v \in K_1$ we have
		$$
		\oplus_n v = \sum_{i=1}^n e_i v e_i^\ast \in K_n, 
		$$
		where $(e_i)_i$ is the standard basis of $\mathbb{C}^n$ satisfying $\sum_{i=1}^n e_i e_i^\ast = \mathbb{I}_n.$
		
		In general, a \ti (\tii) $F \subseteq K_n$ (or its unitary orbit in the case of a \tiii) does not contain any reducible elements of the form $A \oplus B \in M_{r+s} (V)$ for $r, s < n$ as we can express
		$$
		A \oplus B 
		=
		\begin{pmatrix} 
			\mathbb{I}_r \\
			0 
		\end{pmatrix}
		A
		\begin{pmatrix} 
			\mathbb{I}_r & 0 \\
		\end{pmatrix}
		+ 
		\begin{pmatrix} 
			0 \\
			\mathbb{I}_s
		\end{pmatrix}
		B
		\begin{pmatrix} 
			0 & \mathbb{I}_s
		\end{pmatrix}.
		$$
		This gives the intuition that just as points of a matrix convex set are only sporadically matrix extreme, there are in general very few (boundary) points that are contained in a matrix face.
		
		(c) A \ti $F$ is closed under unitary conjugation. Indeed, for $A \in F,$ $B \in K_n$ and a unitary matrix $U \in \mathbb{M}_n,$ by \eqref{eq21} the condition $U^\ast B U = A \in F$ implies $B \in F.$
		
		(d) If $A \in K_n$ is matrix extreme, then its unitary orbit $\mathcal{U}(A) = \{U^\ast A U \ | \ U \in \mathbb{M}_n \text{ unitary} \}$, though in general not convex, satisfies the \ti condition \eqref{eq21}. Indeed, if for a tuple of points $A_1, \ldots, A_k$ from \textbf{$K$} and a tuple of surjective matrices $\gamma_i \in \mathbb{M}_{n_i, n}$ with $\sum_{i=1}^k \gamma_i^\ast \gamma_i = \mathbb{I}_n$ we have
		$$
		U^\ast A U = \sum_{i=1}^k \gamma_i^\ast A_i \gamma_i,
		$$
		where $U \in \mathbb{M}_n$ is a unitary matrix, then $A = \sum_{i=1}^k U\gamma_i^\ast A_i \gamma_iU^\ast$ and $A$ being matrix extreme forces all $A_i$ to be unitarily equivalent to $A.$
		
		(e) A singleton  $F = \{A\}$ is a \tiii if and only if $A$ is a matrix extreme point.
		
		(f) To get rid of the surjectivity of the matrices $\gamma_i$ in \eqref{eq21}, an equivalent condition for a convex set $F \subseteq K_n$ to be a \ti demands that for every tuple of points $A_1, \ldots, A_k$, where each $A_i$ belongs to $K_{n_i}$ for some $n_i \leq n,$ and every tuple of nonzero matrices $\gamma_i \in \mathbb{M}_{n_i, n}$ satisfying $\sum_{i=1}^k \gamma_i^\ast \gamma_i = \mathbb{I}_n$ the condition
		\begin{equation} \label{eq22}
			\sum_{i=1}^k \gamma_i^\ast A_i \gamma_i \in F,
		\end{equation}
		implies $n_i=n$ and $A_i \in F$ for $i=1, \ldots, k.$ It is easy to check that (\ref{eq21}) and (\ref{eq22}) are equivalent. The analogous definitions for $C^\ast$-faces and weak matrix faces can be formulated similarly.
\end{opomba}

The next example will show that as in the case of matrix extreme points, some of the components $K_n$ need not contain any matrix faces. 
\begin{primer}[all types]\label{ex33}
	Given a pair of real numbers $a$ and $b$ with $a < b$ define the corresponding \textbf{matrix interval} $[a \mathbb{I}, b \mathbb{I}]$ := $([a \mathbb{I}_n, b \mathbb{I}_n])_{n \in \mathbb{N}},$ where
	$$
	[a \mathbb{I}_n, b \mathbb{I}_n] := \{\alpha \in \mathbb{M}_n\ | \ a \mathbb{I}_n \preceq \alpha \preceq b \mathbb{I}_n \}.
	$$
	A simple argument in \protect{\cite[Example 2.2]{WW}} shows the only matrix extreme points of $[a \mathbb{I}, b \mathbb{I}]$ are the numbers $a, b \in [a \mathbb{I}, b \mathbb{I}]_1 = [a, b].$ Each $A \in [a \mathbb{I}, b \mathbb{I}]_n$ can be expressed as a matrix convex combination of $a$ and $b.$ So none of the $[a \mathbb{I}, b \mathbb{I}]_n$ for $n > 1$ contains any matrix face.
\end{primer}

In general, a compact matrix convex set \textbf{$K$} over a finite-dimensional space is the (already closed) matrix convex hull of its matrix extreme points by a version of the matricial Krein-Milman theorem \protect{\cite[Theorem 4.3]{WW}} (see also \protect{\cite[Theorem 2.9]{HL}}), and every point $A \in K_n$ can be expressed as a matrix convex combination of matrix extreme points from the sets $K_1, \ldots, K_n.$ 
Hence by part (f) of Remark \ref{op22}, a component $K_n$ which contains a matrix face of any type must also contain a matrix extreme point. Moreover, the elements of a \ti or a \tii can only be described by matrix convex combinations of matrix extreme points of \textbf{$K$}, which lie in $F.$ Similarly, only matrix convex combinations of matrix extreme points that lie in $\mathcal{U}(F)$ can describe points in a \tiii $F$.

An important property of classical faces is that their extreme points are also extreme in the whole set and the next proposition gives a matrix version of it. While the notion of a $C^\ast$-extreme point in a (not necessarily $C^\ast$-convex) convex set might not be natural, this hereditary property of $C^\ast$-extreme points holds for all types of matrix faces.

\begin{trditev}[all types]\label{tr34}
	Let \textbf{$K$} be a matrix convex set and $F \subseteq K_n$ a matrix face. Every $C^\ast$-extreme point of $F$ is a matrix extreme point of \textbf{$K$}.
\end{trditev}

\begin{dokaz}
	First suppose that $F \subseteq K_n$ is a \ti or a \tii and $A$ is a $C^\ast$-extreme point in $F.$ Suppose we can express $A$ as
	\begin{equation}\label{c2}
		A = \sum_{i=1}^k \gamma_i^\ast A_i \gamma_i \in F
	\end{equation}
	for $k$-tuples $(A_i)_{i=1}^k$ and $(\gamma_i)_{i=1}^k$, where each $A_i$ is in $K_{n_i}$ and $\gamma_i \in \mathbb{M}_{n_i,n}$ are surjective matrices with $\sum_{i=1}^k \gamma_i^\ast \gamma_i = \mathbb{I}_n$. Since $F$ is a matrix face, we deduce $n_i=n$ and $A_i \in F$ for $i=1, \ldots, k.$ But then $A$ being $C^\ast$-extreme implies that all the $A_i$ are unitarily equivalent to $A.$
	
	If $F$ is a \tiii, the condition \eqref{c2} implies $n_i=n$ and $A_i \in \mathcal{U}(F)$ for $i=1, \ldots, k.$ So each $A_i$ is of the form $U_i^\ast B_i U_i$ for some unitary $U_i \in \mathbb{M}_n$ and $B_i \in F.$ Then as $A$ is $C^\ast$-extreme, all the $B_i$ and therefore all the $A_i,$ too, are unitarily equivalent to $A.$
\end{dokaz}

\subsection{Matrix exposed faces}\label{subsec42}

This section gives three possible generalisations of the concept of a matrix exposed point to a set, namely a matrix exposed face. 

\begin{definicija} \label{def21}
	Let $\textbf{$K$} = (K_r)_{r \in \mathbb{N}}$ be a matrix convex set in a dual vector space $V$ and $F$ a convex subset of $K_n$. 
	
	(a) Then $F$ is a \textbf{\eti} if there exists a continuous linear map $\Phi:V \to \mathbb{M}_n$ and a self-adjoint matrix $\alpha \in \mathbb{M}_n$ satisfying the following conditions:
	\begin{enumerate}[(i), leftmargin=2cm]
		\item \label{24-a} for every positive integer $m$ and $B \in K_m$ we have $\Phi_m(B) \preceq \alpha \otimes \mathbb{I}_m;$ 
		
		\item for any $m < n$ and $B \in K_m$ we have $\Phi_m(B) \prec \alpha \otimes \mathbb{I}_m;$
		
		\item \label{24-b} $\{ B \in K_n \ | \  \alpha \otimes \mathbb{I}_n - \Phi_n(B) \succeq 0 \text{ is singular}\} = F.$
	\end{enumerate}

	(b) If $F$ is a $C^\ast$-convex \eti, then it is a \textbf{\etii}.
	
	(c) We call $F$ a \textbf{\etiii} if there exists a continuous linear map $\Phi:V \to \mathbb{M}_n$ and a self-adjoint matrix $\alpha \in \mathbb{M}_n$ satisfying the following conditions:
	\begin{enumerate}[(i), leftmargin=2cm]
		\item  for every positive integer $m$ and $B \in K_m$ we have $\Phi_m(B) \preceq \alpha \otimes \mathbb{I}_m;$
		
		\item for any $m < n$ and $B \in K_m$ we have $\Phi_m(B) \prec \alpha \otimes \mathbb{I}_m;$ 
		
		\item $\{ B \in K_n \ | \  \alpha \otimes \mathbb{I}_n - \Phi_n(B) \succeq 0 \text{ is singular}\} = \mathcal{U}(F).$
	\end{enumerate}
\end{definicija}

We call $(\Phi, \alpha)$ in the notation above an \textbf{exposing pair} of the matrix exposed face $F.$\looseness=-1
%
%
%


\begin{opomba}\label{op25}
		(a) As in the case of matrix faces of all types, the matrix exposed faces in $K_1$ of all types coincide with its ordinary faces.
		Indeed, it is clear that every matrix exposed face in $K_1$ is ordinary exposed. For the converse we only need to observe that if $F \subseteq K_1$ is ordinary exposed with an exposing pair $(\varphi, \alpha),$ then the condition $\varphi|_{K_1} \leq \alpha$ implies
		 $$\varphi_m(B) \preceq \alpha \otimes \mathbb{I}_m$$ 
		 for every positive integer $m$ and $B \in K_m.$ Indeed, if 
		 $\varphi_m(B) \npreceq \alpha \otimes \mathbb{I}_m$
		 for some $m \in \mathbb{N}$ and $B \in K_m,$ then there is a unit vector $\xi \in \mathbb{C}^{nm}$ such that 
		 $$
		 0 > \xi^\ast(\alpha \otimes \mathbb{I}_m - \varphi_m(B))\xi = \alpha - \varphi(\xi^\ast B \xi)
		 $$
		But as $\xi^\ast B \xi \in K_1,$ this contradicts the condition $\varphi|_{K_1} \leq \alpha.$
	
		(b) For any matrix affine map $\Phi,$ unitary $U \in \mathbb{M}_n$ and $B \in K_n,$  the matrix 
		$$\alpha \otimes \mathbb{I}_n - \Phi_n(U^\ast B U) = (\mathbb{I}_n \otimes U^\ast)\big( \alpha \otimes \mathbb{I}_n - \Phi_n(B) \big)(\mathbb{I}_n \otimes U)$$
		is singular if and only if $\alpha \otimes \mathbb{I}_n - \Phi_n(B)$ is singular. Hence, the condition (iii) in part (a) of the above definition implies that a \eti is closed under conjugation by unitaries.
		
		(c) Let $A \in K_n$ be a matrix exposed point with the pair $(\Phi, \alpha)$ exposing it. Then by definition and part (c) of Remark \ref{op12}, the same pair  $(\Phi, \alpha)$ satisfies the conditions in part (a) of Definition \ref{def21} for exposing the (in general non-convex) unitary orbit of $A$. 
		
		(d) A singleton $F = \{A\}$ is a \etiii if and only if $A$ is a matrix exposed point.
		
		(e) For a matrix exposed face $F \subseteq K_n$ of any type the intersection: 
		$$
		\mathcal{N} = \bigcap_{A \in F} \ker(\alpha \otimes \mathbb{I}_n - \Phi_n(A))
		$$
		is always nontrivial. If we restrict the search to a unit vector in $\mathcal{N},$ then by the finite intersection property it suffices to prove that for every finite selection of $A_1, \ldots, A_k \in F$ the intersection
		\begin{equation}\label{eq24}
			\bigcap_{i=1}^k \ker\big(\alpha \otimes \mathbb{I}_n - \Phi_n(A_i)\big) \cap S^{n^2}
		\end{equation}
		is nonempty. Here $S^{n^2}$ denotes the unit sphere in $\mathbb{C}^n \otimes \mathbb{C}^n.$ So suppose $A_1, \ldots, A_k \in F$ are such that the intersection (\ref{eq24}) is empty and consider the convex combination $A = \frac{1}{k}\big(A_1 + \cdots + A_k\big) \in F.$ For this point the matrix $\alpha \otimes \mathbb{I}_n - \Phi_n(A)$ is nonsingular, which implies by Definition \ref{def21} that $F$ is not a matrix face. 
\end{opomba}

\begin{trditev}[all types]\label{prop27}
	Let $F \subseteq K_n$ be a matrix exposed face of a matrix convex set \textbf{$K$} and $\Phi:V \to \mathbb{M}_n$ together with $\alpha \in \mathbb{M}_n$ an exposing pair. Then for every nonzero $x = \sum_{i=1}^n x_i \otimes e_i \in \mathbb{C}^n \otimes \mathbb{C}^n$ in $\mathcal{N},$ its components $x_1, \ldots, x_n$ span $\mathbb{C}^n.$ Moreover, $\mathcal{N}$ is one-dimensional.
\end{trditev}	

\begin{dokaz}
	Let $P$ be the projection onto span$\{x_1, \ldots, x_n\}.$ Then as in the proof of part (a) of Proposition \ref{lema14}, $\alpha \otimes \mathbb{I}_r - \Phi_r(PAP^*)$ is singular for any $A \in F,$ but $PAP^* \in K_r$ for some $r <n.$ The proof of part (b) of Proposition \ref{lema14} then shows that $\mathcal{N}$ is one-dimensional.
\end{dokaz}	

\begin{trditev} \label{tr39}
	Let $\textbf{$K$} = (K_n)_{n \in \mathbb{N}}$ be a matrix convex set. Then every \eti $F \subseteq K_n$ is an exposed face of $K_n$.
\end{trditev}

\begin{dokaz}
	Let $F \subseteq K_n$ be a matrix exposed face and $\Phi:V \to \mathbb{M}_n$ together with $\alpha \in \mathbb{M}_n$ an exposing pair. By part (e) of Remark \ref{op25}, there is a nonzero $x \in \mathcal{N}.$ Now define the functional $\varphi: M_n(V) \to \mathbb{C}$ by
	\begin{equation}\label{eq-29}
	\varphi(B) = x^\ast \Phi_n(B) x,
	\end{equation}
	and the real number $a = x^\ast (\alpha \otimes \mathbb{I}_n) x \in \R.$
	Since the pair $(\Phi, \alpha)$ matricially exposes $F,$ we have by the choice of $x$ that
	$$
	a - \varphi(B) = x^\ast\big(\alpha \otimes \mathbb{I}_n - \Phi_n(B)\big) x = 0
	$$
	for all $B \in F,$ but also
	$$
	a - \varphi(B) = x^\ast\big(\alpha \otimes \mathbb{I}_n - \Phi_n(B)\big) x \geq 0
	$$ 
	for all $B \in K_n.$ Moreover, if $B \in K_n$ satisfies $\varphi(B)=a,$ then $(\alpha \otimes \mathbb{I}_n - \Phi_n(B)) x = 0.$ Hence $\alpha \otimes \mathbb{I}_n - \Phi_n(B)$ is singular, so $B \in F$. We conclude that $\varphi(B)=a$ if and only if $B \in F,$ which shows that the pair $(\varphi, a)$ exposes $F$ in $K_n.$
\end{dokaz}

\begin{primer}\label{ex410}
	Suppose $L$ and $M$ are two linear pencils such that $\mathcal{D}_M \subseteq \mathcal{D}_L$ and $\mathcal{D}_M(n) \cap \partial\mathcal{D}_L(n) \neq \emptyset$ for some $n \in \mathbb{N}.$ Let 
	$$n_{\min}:= \min\{n \in \mathbb{N} \ | \ \mathcal{D}_M(n) \cap \partial\mathcal{D}_L(n) \neq \emptyset\}.$$
	If $F= \mathcal{D}_M(n_{\min}) \cap \partial\mathcal{D}_L(n_{\min})$ is convex, then it is clearly a \eti in $\mathcal{D}_M.$ Indeed, if $L = A_0 + \sum_{i=1}^g A_i x_i,$ then the pair $(L - A_0, A_0)$ matricially exposes $F.$  
\end{primer}

\subsection{Interplay between matrix faces and matrix exposed faces} \label{subsec43} 

The aim of this section is to show that a \ti (\tii) is exposed if and only if it is a \eti (\etii). It is not clear if a \etiii is an exposed face, however the ``only if" part of the claim still holds in this case, i.e., an ordinary exposed \etiii is in fact weakly matrix exposed (see Theorem \ref{tr211}).

\begin{opomba}\label{op211} 
 We emphasize that the zero map and the zero matrix always define an exposing pair (of any type) for $K_n,$ even though $K_n$ for $n > 1$ is never a matrix face.
\end{opomba}

\begin{trditev}[all types] \label{tr412}
	Let \textbf{$K$} be a matrix convex set and $F \subsetneq K_n$ a matrix exposed face. Then $F$ is a matrix face.
\end{trditev}

\begin{dokaz} The proof is analogous to that of Proposition \ref{tr15} regarding matrix exposed points. 
	Let $F \subseteq K_n$ be a matrix exposed face of any type with $\Phi:V \to \mathbb{M}_n$ and $\alpha \in \mathbb{M}_n$ an exposing pair. Further, assume we have
	\begin{equation*}
		\sum_{i=1}^k \gamma_i^\ast A_i \gamma_i = A  \in F
	\end{equation*}
	for $k$-tuples $(A_i)_{i=1}^k$ and $(\gamma_i)_{i=1}^k$, where each $A_i$ is in $K_{n_i}$ and $\gamma_i \in \mathbb{M}_{n_i,n}$ are surjective matrices (so $n_i \leq n$) with $\sum_{i=1}^k \gamma_i^\ast \gamma_i = \mathbb{I}_n$. By assumption we have $\mathbb{I}_{n_i} \otimes \alpha - \Phi_{n_i}(A_i) \succeq 0$ for every $i=1, \ldots, k.$ Also
	\begin{align}\label{ie1.3}
		\alpha \otimes \mathbb{I}_n - \Phi_n(A) &= \sum_{i=1}^k (\gamma_i^\ast \otimes \mathbb{I}_n)\,\big(\mathbb{I}_{n_i} \otimes \alpha - \Phi_{n_i}(A_i)\big)\, (\gamma_i\otimes \mathbb{I}_n),
	\end{align} 
	where clearly $(\gamma_i^\ast \otimes \mathbb{I}_n)\,\big(\mathbb{I}_{n_i} \otimes \alpha - \Phi_{n_i}(A_i)\big)\, (\gamma_i\otimes \mathbb{I}_n) \succeq 0$ for $i=1, \ldots, k.$
	Suppose any of the points $A_i,$ say $A_1,$ is not in $F$ (not in $\mathcal{U}(F)$ for the case of a \etiii) and hence satisfies $\mathbb{I}_{n_i} \otimes \alpha - \Phi_{n_i}(A_i) \succ 0$ by definition. We will prove that the latter implies $\gamma_1 = 0$.
	By Proposition \ref{prop27}, there is an $x = \sum_{i=1}^n x_i \otimes e_i \in \mathcal{N} = \bigcap_{A \in F} \ker(\alpha \otimes \mathbb{I}_n - \Phi_n(A))$ such that span$\{x_1, \ldots, x_n\} = \mathbb{C}^n.$ Equality (\ref{ie1.3}) together with the positive semi-definiteness of the right-hand side summands imply that $x$ lies in the intersection of the kernels of the matrices $(\gamma_i^\ast \otimes \mathbb{I}_n)\,\big(\mathbb{I}_{n_i} \otimes \alpha - \Phi_{n_i}(A_i)\big)\, (\gamma_i\otimes \mathbb{I}_n)$ for $i=1, \ldots, k.$ In particular, we have $(\gamma_1^\ast \otimes \mathbb{I}_n)\,\big(\mathbb{I}_{n_1} \otimes \alpha - \Phi_{n_1}(A_1)\big)\, (\gamma_1\otimes \mathbb{I}_n)x = 0,$ which, together with the positive definiteness of the middle factor and injectivity of $(\gamma_1^\ast \otimes \mathbb{I}_n),$ implies that $x$ lies in the kernel of $\gamma_1 \otimes \mathbb{I}_n,$ i.e.,~we have:
	$$
	0 = (\gamma_1 \otimes \mathbb{I}_n)x = (\gamma_1 \otimes \mathbb{I}_n) \bigg( \sum_{j=1}^n x_j \otimes e_j\bigg) = \sum_{j=1}^n \gamma_1 x_j \otimes e_j.
	$$
	So $\gamma_1 x_j^i = 0$ for $i=1, \ldots, d$ and $j=1,\ldots, n,$ whence $\gamma_1 = 0.$
\end{dokaz}

\begin{izrek}[all types] \label{tr211} Let \textbf{$K$} be a matrix convex set and $F \subseteq K_n$ a matrix face that is also an exposed face.  Then $F$ is a matrix exposed face.
\end{izrek}

\begin{dokaz}
	First, suppose $F \subseteq K_n$ is a \ti or a \tii, hence it is closed under conjugation by unitaries, and suppose $F$ is also exposed in the classical sense. By part (a) of Remark \ref{op25} we may assume $n > 1.$ Moreover, we may assume $0_n \notin F$ and $0 \in  K_1.$ Indeed, otherwise replace \textbf{$K$} by $-v \,+\, $\textbf{$K$} for any vector $v \in K_1.$
	Note that $\mathbb{I}_n \otimes v = \oplus_n v \notin F$ as $F$ does not contain reducible elements by part (b) of Remark \ref{op22}. By assumption, we also have $F \subsetneq K_n$ since $K_n$ is never a matrix face.
	
	By considering the matrix convex set $\textbf{$L$} =$ mconv$(K_n\backslash F),$ which is disjoint from $F,$ we construct, analogously to the proof of part (b) of Theorem \ref{tr111}, from the pair $(\varphi, a)$ exposing $F,$ a continuous linear map $\Phi: V \to \mathbb{M}_n$ and a self-adjoint matrix $\alpha = a \mathbb{I}_n \in \mathbb{M}_n$ that satisfy $\Phi_m(B) \preceq \alpha \otimes \mathbb{I}_m$ for all $m \in \mathbb{N}$ and $B \in K_m.$ Also, the following relation holds:\looseness=-1
	$$
	\big\langle\Phi_n(A) \eta_0, \eta_0\big\rangle
	= \varphi(A) = a
	$$
	for some unit vector $\eta_0 \in (\mathbb{C}^n)^n$ and all $A \in F.$
	So for all $A \in F$ the corresponding matrix $a\,\mathbb{I}_{n^2} - \Phi_n(A)$ is singular.

	We now show that for a point $B \in K_n \backslash F$ the matrix $a\,\mathbb{I}_{n^2} - \Phi_r(B)$ is nonsingular. So suppose some unit vector $\eta \in (\mathbb{C}^n)^n$ of the form \eqref{eta1} satisfies $\big\langle\text{Re}\,\Phi_n(B) \eta, \eta\big\rangle = \text{Re}\,\varphi(\alpha^\ast B \alpha) =  a\,p(\alpha^\ast \alpha) = a.$ By the same reasoning as in \eqref{eq-110} we see that $\alpha$ is a contraction. As we assumed $0 \in L_1,$ the point $\alpha^\ast B \alpha$ belongs to $L_n,$ since \textbf{$L$} is closed under conjugation by contractions.
	Because of the strong separation of $F$ from $L_n$ determined by $(\varphi, a),$ we get $\text{Re}\,\varphi(\alpha^\ast B \alpha) <  a\,p(\alpha^\ast \alpha)$, which is a contradiction.
	The above reasoning can be easily adapted for case of a \tiii by considering the matrix convex set mconv$(K_n\backslash \mathcal{U}(F))$ instead of mconv$(K_n\backslash F).$

Finally, we claim that for any $m < n$ and $B \in K_m$ the strict inequality $\Phi_m(B) \prec \alpha \otimes \mathbb{I}_m$ holds. So suppose $m < n$ and $B \in K_m$ are such that $\alpha \otimes \mathbb{I}_m - \Phi_m(B)$ is singular (while also positive semidefinite). Then by part (b) of Remark \ref{op22}, $B \oplus C \notin F$ for any $C \in K_{n-m}.$ By the singularity of 
	$$
	\alpha \otimes \mathbb{I}_n - \Phi_n(B \oplus C) = \big(\alpha \otimes \mathbb{I}_m - \Phi_m(B)\big) \oplus \big(\alpha \otimes \mathbb{I}_{n-m} - \Phi_{n-m}(C)\big),
	$$ 
	this contradicts $F$ being a matrix exposed face (cf.~part (c) of Remark \ref{op12}). The claim also holds if $F$ is a \etiii, where in the above reasoning we use the fact that $\mathcal{U}(F)$ does not contain any reducible elements. Having proven the last claim, we conclude that the pair $(\Phi, a\,\mathbb{I}_n)$ matricially exposes $F.$
\end{dokaz}

Recall how a matrix convex set $\textbf{$K$}$ determines a family of convex sets $\{\Gamma_n(\textbf{$K$})\}_{n \in \mathbb{N}}$ given by \eqref{eq311} in Section \ref{subsec33}.
We proceed by showing that for each $n \in \mathbb{N},$ there is a connection between matrix exposed faces of $\textbf{$K$}$ and exposed faces of $\Gamma_n(\textbf{$K$})$ analogous to the one regarding points in Proposition \ref{tr19}. We state it as Proposition \ref{tr416}, whose proof mimics the one of Proposition \ref{tr19}.

\begin{trditev}\label{tr416}
	Let $\textbf{$K$} = (K_n)_{n \in \mathbb{N}}$ be a matrix convex set and $F \subseteq K_r.$ Then the following holds:
	\begin{enumerate}[\text\textnormal{{(a)}}]
		\item Let $\gamma \in \mathbb{M}_{r, n}$ be a surjective matrix with $\text{\textnormal{tr}}(\gamma^\ast \gamma)=1$ such that the set 
		\begin{equation}\label{eq47}
			F_{\gamma} := \{(\gamma^\ast \gamma, \gamma^\ast A \gamma) \ | \ A \in F\}
		\end{equation}
		is an exposed face of $\Gamma_n(\textbf{$K$}).$
		Then $F$ is a \etiii of \textbf{$K$}.
	\end{enumerate}
	\begin{enumerate}[\text\textnormal{{(b)}}]
		\item If $F$ is a \etiii of \textbf{$K$}, then for any invertible $\gamma \in \mathbb{M}_{r}$ with $\text{\textnormal{tr}}(\gamma^\ast \gamma)=1,$ the set $F_{\gamma}$ in \eqref{eq47}
		is an exposed face of $\Gamma_r(\textbf{$K$}).$
	\end{enumerate}
\end{trditev}

\begin{dokaz}
	To prove (a) suppose $F_{\gamma}$ as above is an exposed face of $\Gamma_n(\textbf{$K$})$ for some surjective $\gamma \in \mathbb{M}_{r,n}$ with $\text{\textnormal{tr}}(\gamma^\ast \gamma)=1.$ By assumption there exists a continuous linear functional $\varphi: \mathbb{M}_n \times M_n(V) \to \mathbb{C}$ and a real number $a$ such that $\varphi(\gamma^\ast \gamma, \gamma^\ast A \gamma) = a$ and $\varphi(\delta, C) > a$ for all $(\delta, C) \in \Gamma_n(\textbf{$K$}) \backslash F_{\gamma}.$ 
	
	Now decompose $\varphi$ into $\varphi = \varphi_1 + \varphi_2,$ where $\varphi_1 : \mathbb{M}_n \to \mathbb{C}$ and $\varphi_2 : M_n(V) \to \mathbb{C}.$ By the Riesz representation theorem there is a self-adjoint matrix $\beta \in \mathbb{M}_n$ such that $\varphi_1(\delta) = \text{tr}(\beta \delta).$ For any matrix $\eta = (\eta_1, \ldots, \eta_n) \in \mathbb{M}_{m, n}$ and its vectorization 
	\begin{equation}\label{eq-312}
		v_{\eta}= \sum_{i=1}^n e^{\prime}_i \otimes \eta_i \in \mathbb{C}^n \otimes \mathbb{C}^m
	\end{equation}
	we have
	$$
	\varphi_1(\eta^\ast \eta) = \text{tr}(\beta \eta^\ast \eta) = \text{tr}(\eta \beta \eta^\ast) =v_{\eta}^\ast(\beta \otimes \mathbb{I}_m) v_{\eta}.
	$$
	By the canonical correspondence between linear functionals on $M_n(V)$ and linear maps $V \to \mathbb{M}_n,$ the functional $\varphi_2$ gives rise to a map $\Psi: V \to \mathbb{M}_n.$ For every $B \in M_m(V)$ and matrix $\eta = (\eta_1, \ldots, \eta_n) \in \mathbb{M}_{m, n}$ with its vectorization $v_{\eta}$ as in \eqref{eq-312} we have
	$$
	\varphi_2(\eta^\ast B \eta) = v_{\eta}^\ast \Psi_m(B) v_{\eta}.
	$$
	
	So for every positive integer $m,$ element $B \in K_m$ and matrix $\eta \in \mathbb{M}_{m, n}$ with the property tr$(\eta^\ast \eta) = 1,$ 
	\begin{align}\label{eq49}
		v_{\eta}^\ast \big( (\beta - a \mathbb{I}_n)  \otimes \mathbb{I}_m - \Psi_m(B)\big)  v_{\eta} = \varphi(\eta^\ast \eta, \eta^\ast B \eta) - a \geq 0.
	\end{align}
	Since tr$(\eta^\ast \eta) = 1$ if and only if $v_{\eta} \in \mathbb{C}^{mn}$ is a unit vector, the above implies that the matrix $(\beta - a \mathbb{I}_n) \otimes \mathbb{I}_m - \Psi_m(B)$ is positive semidefinite. For $A \in F$ we also have
	\begin{align} \label{eq410}
		v_{\gamma}^\ast \big( (\beta - a \mathbb{I}_n) \otimes \mathbb{I}_r - \Psi_r(A)\big)  v_{\gamma} = \varphi(\gamma^\ast \gamma, \gamma^\ast A \gamma) - a= 0.
	\end{align}
	The positive semi-definiteness of $(\beta - a \mathbb{I}_n) \otimes \mathbb{I}_r - \Psi_r(A)$ implies $\big( (\beta - a \mathbb{I}_n) \otimes \mathbb{I}_r - \Psi_r(A)\big)  v_{\gamma} = 0).$ So $(\beta - a \mathbb{I}_n) \otimes \mathbb{I}_r - \Psi_r(A)$ is singular for every $A \in F$. 
	
	Again, by the Douglas Lemma \ref{lema18}, $(\eta^\ast \eta, \eta^\ast B \eta) \in \Gamma_n(\textbf{$K$})$ with $B \in K_m$ and $\eta \in \mathbb{M}_{m, n}$ is of the form $(\gamma^\ast \gamma, \gamma^\ast A \gamma)$ for some $A \in F$ if and only if $r=m$ and there is a unitary matrix $U \in \mathbb{M}_r$ such that $\eta = U \gamma$ and $A = U^\ast B U.$ By the calculation (\ref{al18}) we have
	\begin{align*}
		\varphi(\eta^\ast \eta, \eta^\ast B \eta) - a =  v_{\eta}^\ast \big( (\beta - a \mathbb{I}_n) \otimes \mathbb{I}_r - \Psi_r(B)\big)  v_{\eta}. \nonumber
	\end{align*}
	Using the properties of $\varphi$ we deduce that for $B \in K_m$ the matrix $(\beta - a \mathbb{I}_n) \otimes \mathbb{I}_m - \Psi_m(B)$  is singular if and only if $m = r$ and $B\in \mathcal{U}(F).$
	
	Finally, we define the map $\Phi: V \to \mathbb{M}_r$ by $\Phi = \gamma \Psi \gamma^\ast$ and the self-adjoint matrix $\alpha \in \mathbb{M}_r$ by $\alpha = \gamma(\beta - a \mathbb{I}_r) \gamma^\ast.$ Now the same reasoning as in the proof of part (a) of Proposition \ref{tr19} shows that $v_\gamma$ is in the range of $(\gamma^* \otimes \mathbb{I}_r)$ so that $\alpha \otimes \mathbb{I}_r - \Phi_r(A)$ is singular for all $A\in F.$ Since the positive semi-definiteness in all the points of \textbf{$K$} is clear and the definiteness outside of $\mathcal{U}(F)$ follows from the injectivity of $\gamma^*,$ the pair $(\Phi, \alpha)$ matricially exposes $F$ and $F$ is a \etiii.
	
	To prove (b) assume $F \subseteq K_r$ is a \etiii and let $\Phi:V \to \mathbb{M}_r$ together with $\alpha \in \mathbb{M}_r$ be an exposing pair.
	We claim that for any invertible $\gamma \in \mathbb{M}_{r}$ with $\text{\textnormal{tr}}(\gamma^\ast \gamma)=1,$ $F_{\gamma} = \{(\gamma^\ast \gamma, \gamma^\ast A \gamma) \ | \ A \in F\}$ is an exposed face of $\Gamma_r(\textbf{$K$}).$ 
	
	By Proposition \ref{prop27}, there is an $x = \sum_{i=1}^r x_i \otimes e_i \in \mathbb{C}^r \otimes \mathbb{C}^r$ in $\mathcal{N} = \bigcap_{A \in F} \ker(\alpha \otimes \mathbb{I}_r - \Phi_r(A))$ such that span$\{x_1, \ldots, x_r\} = \mathbb{C}^r.$ Let $(\Psi, \beta) = \delta(\Phi, \alpha)\delta^*$ and $y = (\delta \otimes \delta)x,$ where $\delta \in \mathbb{M}_r$ denotes the inverse of $\gamma,$ and define a functional $\varphi: \mathbb{M}_r \times M_r(V) \to \mathbb{C},$ 
	$$
	\varphi(\mu, C) =  y^\ast\big(\beta \otimes \mu -\Psi_r(C)\big) y
	$$
	for $C \in M_r(V)$ and $\mu \in \mathbb{M}_r.$ 
	Since $\Psi$ is matrix affine, we see that for any positive integer $n,$ element $B \in K_n$ and surjective matrix $\eta \in \mathbb{M}_{n, r}$ with tr$(\eta^\ast \eta) = 1,$ 
	$$
	\varphi(\eta^\ast \eta, \eta^\ast B \eta) = 
	\big((\eta \otimes \mathbb{I}_r)y\big)^\ast \big( \beta \otimes \mathbb{I}_n - \Psi_n(B)\big)  (\eta \otimes \mathbb{I}_r)y \geq 0.
	$$
	By Proposition \ref{prop27}, the same reasoning as in the proof of part (b) of Proposition \ref{tr19} shows that if $(\eta^\ast \eta, \eta^\ast B \eta) \in \Gamma_n(\textbf{$K$}) \backslash F_{\gamma},$ then $\varphi(\eta^\ast \eta, \eta^\ast B \eta)>0.$ On the other hand, for $(\gamma^\ast \gamma, \gamma^\ast A \gamma) \in F_{\gamma}$ we have
	\begin{align*}
		\varphi(\gamma^\ast \gamma, \gamma^\ast A \gamma) &= ((\mathbb{I}_r \otimes \gamma)y)^\ast \big(\beta \otimes \mathbb{I}_r  - \Psi_r(A)\big)  (\mathbb{I}_r \otimes \gamma)y \\
		&= ((\delta \otimes \mathbb{I}_r)x)^\ast \big(\beta \otimes \mathbb{I}_r  - \Psi_r(A)\big)  (\delta \otimes \mathbb{I}_r)x \\
		&= x^*\big(\alpha \otimes \mathbb{I}_r  - \Phi_r(A)\big)x=0.
	\end{align*}
	Hence the pair $(-\varphi, 0)$ exposes $F_{\gamma}$ in $\Gamma_n(\textbf{$K$}).$	
\end{dokaz}

Next we state a corollary of Proposition \ref{tr416} regarding matrix exposed faces. The  proof is omitted as it is similar to the one of Proposition \ref{tr416}. 

\begin{posledica} \label{tr315}
	Let $\textbf{$K$} = (K_n)_{n \in \mathbb{N}}$ be a matrix convex set and $F$ a subset $K_r$ that is closed under conjugation by unitaries. Then the following holds:
	\begin{enumerate}[\text\textnormal{{(a)}}]
		\item Let $\gamma \in \mathbb{M}_{r, n}$ be a surjective matrix with $\text{\textnormal{tr}}(\gamma^\ast \gamma)=1$ such that the set $F_{\gamma}$ in \eqref{eq47}
		is an exposed face of $\Gamma_n(\textbf{$K$}).$
		Then $F$ is a \eti of \textbf{$K$}.
	\end{enumerate}
	\begin{enumerate}[\text\textnormal{{(b)}}]
		\item If $F$ is a \eti of \textbf{$K$}, then for any invertible $\gamma \in \mathbb{M}_{r}$ with $\text{\textnormal{tr}}(\gamma^\ast \gamma)=1,$ the set $F_{\gamma}$ in \eqref{eq47}
		is an exposed face of $\Gamma_r(\textbf{$K$}).$
	\end{enumerate}
\end{posledica}

\begin{primer}\label{ex416}
	 Let \textbf{$K$} be a compact matrix convex set in the space $V.$ As a corollary of Proposition \ref{tr416} we now give an example of a \etiii of \textbf{$K$}. Since \textbf{$K$} is compact, $\Gamma_n(\textbf{$K$})$ is a compact subset of the finite-dimensional space $\mathbb{M}_n \times M_n(V)$ and so for any fixed $r < n,$ the minimum
	 $$
	 m := \min\{\|\gamma^\ast \gamma\| \ | \ \gamma \in \mathbb{M}_{r, n}, \  \text{tr}\,(\gamma^\ast \gamma) = 1, \text{ such that } \exists A \in \mathbb{M}_r \text{ with } (\gamma^\ast \gamma, \gamma^\ast A \gamma) \in \Gamma_n(\textbf{$K$}) \}
	 $$ 
	 is attained in some $\gamma \in \mathbb{M}_{r, n}$ with tr$(\gamma^\ast \gamma) = 1$. Then
	 $$F_{\gamma} = \{(\delta^\ast \delta, \delta^\ast A \delta) \ | \ A \in K_r,\  \delta \in \mathbb{M}_{r, n},\ \text{tr}\,(\delta^\ast \delta) = 1, \ \delta^\ast \delta = \gamma^\ast \gamma\}$$ 
	is the intersection of $\Gamma_n(\textbf{$K$})$ with the affine plane $\{(\epsilon, B) \in \mathbb{M}_n \times M_n(V) \ | \ \epsilon = \gamma^\ast \gamma\}.$ Hence $F_{\gamma}$ is an exposed face of $\Gamma_n(\textbf{$K$})$ which by Proposition \ref{tr416} implies that
	 $$
	 F=\{A \in K_r \ | \ \exists \delta \in \mathbb{M}_{r, n},\ \text{tr}\,(\delta^\ast \delta) = 1, \  \delta^\ast \delta = \gamma^\ast \gamma \text{ and } (\gamma^\ast \gamma, \gamma^\ast A \gamma) \in F_{\gamma}\}
	 $$ 
	 is a \etiii of \textbf{$K$}.
\end{primer}

\subsubsection{Matrix faces and matrix exposed faces in free spectrahedra}\label{sub331}
		As an application of Proposition \ref{tr211} we can deduce that every \ti (\tii) of a free spectrahedron is a \eti (\etii). Indeed, by Proposition \ref{tr211} it is enough to prove every face of a free spectrahedron is exposed. For this we apply the result \protect{\cite[Corrolary 1]{RG}} by Ramana and Goldman stating that every face of a spectrahedron is exposed. But since $\mathcal{D}_L(n)$ can be considered a spectrahedron in the Euclidean space $\mathbb{S}_n^g$ for arbitrary $n \in \mathbb{N},$ every face of a free spectrahedron is exposed, too.
		
		 We conclude with an insight into the structure of (exposed) faces of a free spectrahedron. First recall that for a convex set $K \subseteq \R^n$ and point $x \in K$, there is a unique face $F_K(x)$ of $K$ which contains $x$ in its relative interior (see \cite[Section II.2]{Ba}). In general we have\looseness=-1 
		$$F_K(x) = \text{aff}(F_K(x)) \cap K,$$
		where $\text{aff}(F_K(x))$ denotes the affine span of $F_K(x).$ Now let $n \in \mathbb{N}$ and $A_0, \ldots, A_g$ be complex self-adjoint matrices of size $n \times n$ and let $L = A_0 + \sum_{i=1}^g A_i x_i$ be the corresponding linear pencil.
		The following theorem (cf.~\protect{\cite[Theorem 1]{RG}}) gives for a point $X$ in the free spectrahedron $\mathcal{D}_L$ a concrete description of the unique face $F_L(X)$ that contains $X$ in its relative interior. 
		
		\begin{izrek}
			Let $L$ be a linear pencil and $X \in \mathcal{D}_L(n).$ Then
			\begin{align*}
				F_L(X) &= \{Y \in \mathcal{D}_L(n) \ | \ \ker\,L(Y) \supseteq \ker\,L(X)\}  \\ 
				&= \{Y \in \mathbb{S}^g_n \ | \ \ker\,L(Y) \supseteq \ker\,L(X)\} \cap \mathcal{D}_L(n)\\ 
				&=\{Y \in  \mathcal{D}_L(n) \ | \ x^\ast L(Y)x = 0 \ \forall x \in \ker\,L(X)\}  \\
				&= \{Y \in \mathbb{S}^g_n \ | \ x^\ast L(Y)x = 0 \ \forall x \in \ker\,L(X)\}  \cap \mathcal{D}_L(n).
			\end{align*}
			Moreover, $\text{\textnormal{aff}}(F_L(X)) = \{Y \in \mathbb{S}^g_n \ | \ \ker\,L(Y) \supseteq \ker\,L(X)\}.$\qedhere
		\end{izrek}
 
We now present a sufficient condition to determine for which $X \in \mathcal{D}_L$ the face $F_L(X)$ is a matrix (exposed) face. The following is an easy corollary of Proposition \ref{tr416}, adapted to free spectrahedra.

\begin{posledica} 
	Let $L$ be a linear pencil and $F \subseteq \mathcal{D}_L(r).$ 
	Suppose $F$ is closed under conjugation by unitaries. If for some $m \in \mathbb{N}$ and surjective $\gamma \in \mathbb{M}_{r, m}$ with $\text{\textnormal{tr}}(\gamma^\ast \gamma)=1$ the set 
	$$F_{\gamma} = \{(\gamma^\ast \gamma, \gamma^\ast X \gamma) \ | \ X \in F\}$$ 
	is an exposed face of 
	\begin{align*}
		\Gamma_m(\mathcal{D}_L) &= \{(\delta^\ast \delta, \delta^\ast Y \delta) \ | \ \delta \in \mathbb{M}_{k, m} \text{\textnormal{ onto}}, \text{\textnormal{ tr}}\,(\delta^\ast \delta) = 1, k \in \mathbb{N}, Y \in \mathcal{D}_L(k)\} \\
		&= \{(\delta^\ast \delta, \delta^\ast Y \delta) \ | \  \delta \in \mathbb{M}_{k, m} \text{\textnormal{ onto}}, \text{\textnormal{ tr}}\,(\delta^\ast \delta) = 1, k \in \mathbb{N}, Y \in \mathbb{S}_n^g, \\ & \hspace{3.1cm} L(Y)\succeq 0\} \subseteq \mathbb{M}_m \times \mathbb{M}_m,
	\end{align*}
then $F$ is a \eti of $\mathcal{D}_L.$
\end{posledica}

	As it is not clear if weak matrix (exposed) faces are ordinary (exposed) faces, the observations from this subsection do not directly apply to weak matrix faces.

\section{Multilevel matrix faces}\label{sec5}

\subsection{Matrix multifaces}
In this section we discuss two notions of a multicomponent face of a matrix convex set. We show how the classical theory connecting (archimedean) faces of compact convex sets and (archimedean) order ideals of the corresponding function systems presented in \cite[Section II.5]{Al} has its noncommutative counterpart featuring matrix multifaces. It also gives rise to a family of examples along with a sufficient condition to deduce whether a point is contained in some matrix multiface.

A notion similar to that of a matrix multiface recently appeared under the name nc face in \cite{KKM}. 
While nc faces extend the concept of absolute extreme points (see \cite{EHKM}), our interest is in generalizing the properties of matrix extreme points.

\begin{definicija}\label{def51}
	Let $\textbf{$K$} = (K_r)_{r \in \mathbb{N}}$ be a matrix convex set in the space $V$ and $\textbf{$F$} = (F_r)_{r \in \mathbb{N}} \subseteq \textbf{$K$}$ a levelwise convex subset of \textbf{$K$}.  
	
	(a) Then $\textbf{$F$}$ is a \textbf{\tmi} if for every tuple of points $A_1, \ldots, A_k$ from \textbf{$K$} and every tuple of surjective matrices $\gamma_i \in \mathbb{M}_{n_i, n}$ satisfying $\sum_{i=1}^k \gamma_i^\ast \gamma_i = \mathbb{I}_n,$ the condition
	\begin{equation} \label{eq411}
		\sum_{i=1}^k \gamma_i^\ast A_i \gamma_i \in \textbf{$F$},
	\end{equation}
	implies $A_i \in \textbf{$F$}$ for $i=1, \ldots, k.$
	
	(b) If $\textbf{$F$}$ is a matrix convex \tmi, then it is a \textbf{\tmii}.
\end{definicija}

\begin{opomba}\label{op422a}
	(a) It is straightforward that each component $F_n$ of a matrix (convex) multiface \textbf{$F$} is an ordinary face.
	
	(b) For $n \in \mathbb{N}$ and $F \subseteq K_n$ denote by $\widehat{F}$ the subset of \textbf{$K$} with $n$-th component $F$ and the other components being empty. Then for every \ti $F \subseteq K_n,$ the corresponding multicomponent set $\widehat{F}$ is a \tmi. Moreover, the matrix multifaces $\textbf{$F$} = (F_r)_{r \in \mathbb{N}}$ with $F_r = \emptyset$ for $r > 1$ coincide with subsets of \textbf{$K$}, whose first components are ordinary faces. Also, $\widehat{K_1}$ itself is a \tmi and we see as in part (b) of Remark \ref{op22} that $\widehat{K_n}$ for $n > 1$ is never a matrix multiface. Moreover, as in part (c) of Remark \ref{op22}, every \tmi is closed under unitary conjugation.
	
	
	(c) As in part (c) of Definition \ref{def41} one might consider weak matrix multifaces  to obtain that for any weak matrix multiface \textbf{$F$}, where $F_n = \{A\}$ for some $n \in \mathbb{N}$ and $F_m = \emptyset$ whenever $m \neq n,$ the point $A$ is matrix extreme. 
\end{opomba}

\begin{primer}\label{ex418} We give an example of a matrix face in a matrix convex set, inspired by the theory connecting (archimedean) faces of compact convex sets and (archimedean) order ideals of the corresponding function systems presented in \cite[Section II.5]{Al}. 
	Let \textbf{$K$} be a compact matrix convex set and denote by 
	$$
	A(\textbf{$K$}) = \{\theta = (\theta_n: K_n \to \mathbb{M}_n)_{n \in \mathbb{N}} \ | \ \theta \text{ continuous matrix affine}\}
	$$ 
	its dual operator system. For a ucp map $\Phi$ on $A(\textbf{$K$})$ with kernel $J$ let 
	\begin{equation}\label{eq412}
		J_n^{\perp} := \{A \in K_n \ | \ \theta_n(A) = 0 \ \  \forall \theta \in J\}.
	\end{equation}
	A straightforward calculation shows that $\textbf{$J$}^{\perp}:= (J_n^{\perp})_{n \in \mathbb{N}}$ is a matrix convex subset of \textbf{$K$}.
	We will show that if $J$ is spanned by its positive elements, i.e., $J= J^+ - J^+,$ then $\textbf{$J$}^{\perp}$ is a \tmii. 
	
	First note that for any $r \in \mathbb{N},$ the ampliation $\Phi_r$ also has the kernel generated by its positive elements. Indeed, if 
	$$
	\Phi_r(A) = (\Phi(A_{i, j})) = 0
	$$ 
	for some $A = (A_{i, j}) \in M_r(A(\textbf{$K$})),$ then every $A_{i, j}$ lies in $J= J^+ - J^+.$ Now the claim follows, since we can write $A = \sum_{i, j}E_{i, j} \otimes  A_{i, j},$ each standard base matrix $E_{i, j}$ can be expressed as a (complex) linear combination of positive matrices, and the elements of the form $\alpha \otimes \theta$ for $\alpha \in \mathbb{M}_r^+$ and $\theta \in A(\textbf{$K$})^+$ lie in $M_r(A(\textbf{$K$}))^+.$

	Now to prove $\textbf{$J$}^{\perp}$ is a \tmii suppose
	$$
	A = \sum_{i=1}^k \gamma_i^\ast A_i \gamma_i \in J_n^{\perp}
	$$
	for $k$-tuples $(A_i)_{i=1}^k$ and $(\gamma_i)_{i=1}^k$, where $A_i \in K_{n_i}$ and $\gamma_i \in \mathbb{M}_{n_i,n}$ are onto with $\sum_{i=1}^k \gamma_i^\ast \gamma_i = \mathbb{I}_n.$ Then for any $\theta \in J,$
	$$
	0 = \theta_n(A) = \sum_{i=1}^k \gamma_i^\ast \theta_{n_i}(A_i) \gamma_i.
	$$
	If $\theta$ is a positive element in $A(\textbf{$K$}),$ the above implies $\theta_{n_i}(A_i) = 0$ for $i=1, \ldots, k$ and since $J$ is spanned by its positive elements, we have $\theta_{n_i}(A_i) = 0$ for all $i$ and $\theta \in J.$ Hence $A_i \in \textbf{$J$}^{\perp}$ for $i=1, \ldots, k$ so that $\textbf{$J$}^{\perp}$ is a matrix multiface.
	\end{primer}

\begin{opomba} \label{op422}
		 
	 (a) By inspection of the observations of Example \ref{ex418}, we deduce a sufficient condition for a point $X \in \textbf{$K$}$ to be contained in a matrix multiface. Identifying $X$ with the corresponding evaluation map $\Phi_X \in \textbf{UCP}(A(\textbf{$K$})),$ we see that if every $\theta \in A(\textbf{$K$})$ with $\theta(X)=0$ can be decomposed as $\theta = \theta_1 - \theta_2,$ where $\theta_1, \theta_2 \in A(\textbf{$K$})^+$ and $\theta_1(X)= \theta_2(X) = 0,$ then $(\ker\,\Phi_X)^{\perp}$ is a matrix (convex) multiface that contains $X.$
	 
	 (b) One might try to adapt Example \ref{ex418} to obtain a fixed-level matrix face by assuming that for a ucp map $\varphi$ on $A(\textbf{$K$})$ the corresponding kernel $J$ is spanned by its positive elements and
	 $$
	 n_{\text{min}} := \min\,\{n \in \mathbb{N} \ | \ J_n^{\perp} \neq \emptyset\} < \infty.
	 $$
	 Then the set $J_{\min}^{\perp} := J_{n_{\min}}^{\perp}$ satisfies the conditions of a \tii. However, if $n_{\text{min}} < \infty,$ then $n_{\text{min}}=1$ because of the connection
	 $$
	 \xi^\ast \theta_n(A) \xi = \theta_1(\xi^\ast A \xi)
	 $$
	 for any $A \in K_n,$ matrix affine map $\theta \in A(\textbf{$K$})$ and unit vector $\xi \in \mathbb{C}^n.$ So the presented construction only reproduces some of the faces of $K_1.$
\end{opomba}

\begin{definicija} Let \textbf{$K$} be a matrix convex set.
	A ucp map $\Phi: A(\textbf{$K$}) \to \mathbb{M}_r$ is \textbf{partially order reflecting} if it satisfies 
	\begin{equation}\label{eq413}
		\Phi_n\big(M_n(A(\textbf{$K$}))^+\big) = \Phi_n\big(M_n(A(\textbf{$K$}))\big)^+
	\end{equation}
	for all $n \in \mathbb{N}$, i.e., for every $n \in \mathbb{N}$ and $A \in M_n(A(\textbf{$K$}))$ with $\Phi_n(A) \succeq 0$ there exists a $B \succeq 0$ such that $\Phi_n(A) = \Phi_n(B).$
\end{definicija}

The following two propositions together give a noncommutative analogue of \cite[II.5.11 -- II.5.13]{Al}.

\begin{trditev}\label{tr420}
	Let \textbf{$K$} be a compact matrix convex set and $\Phi: A(\textbf{$K$}) \to \mathbb{M}_r$ a ucp map. Denote by $J$ the kernel of $\Phi$ and $\textbf{$J$}^{\perp}= (J_m^{\perp})_{m \in \mathbb{N}}$ as in \eqref{eq412}. Let $q: A(\textbf{$K$}) \to A(\textbf{$K$}) / J$ be the canonical quotient projection. Then the following are equivalent:
	\begin{enumerate}[\text\textnormal{{(a)}}]
	\item The kernel $J$ is spanned by its positive elements and the quotient matrix ordered space $A(\textbf{$K$}) / J$ with positive cones $( q(M_n(A(\textbf{$K$}))^+))_n$ and unit $q(1_\textbf{$K$})$
	 is an operator system.
	\end{enumerate}
	\begin{enumerate}[\text\textnormal{{(b)}}]
	\item For each $n \in \mathbb{N}$ and $\theta \in M_n(A(\textbf{$K$}))$ with $\theta|_{\textbf{$J$}^{\perp}} \succeq 0$ there is a positive element $\psi \in M_n(A(\textbf{$K$}))^+$ such that 
	$$\psi \succeq \theta \quad \text{ and } \quad \psi|_{\textbf{$J$}^{\perp}} = \theta|_{\textbf{$J$}^{\perp}}.
	$$ 
\end{enumerate}
\end{trditev}

\begin{dokaz} 
	(a) $\Rightarrow$ (b)  Since $A(\textbf{$K$}) / J$ is an operator system, it can be equipped with a norm (see \cite[Proposition 13.3]{Pa}) in which the positive cones $(A(\textbf{$K$}) / J)^+_n := q(M_n(A(\textbf{$K$}))^+)$ are closed for all $n.$ Moreover, $A(\textbf{$K$}) / J$ satisfies the compatibility condition on its positive cones \cite[Chapter 13]{Pa}, which states that for any $n \times m$ complex matrix $\alpha,$ 
	$$
	\alpha^* \,(A(\textbf{$K$}) / J)^+_n \,\alpha \subseteq (A(\textbf{$K$}) / J)^+_m.
	$$
	This immediately implies that $((A(\textbf{$K$}) / J)^+_n)_n$ is a (closed) matrix convex set.
	
	 Observe that for any $\Psi \in \text{UCP}(A(\textbf{$K$}) / J),$ which is continuous with respect to the above norm, the composition $\overline{\Psi} = \Psi \circ q$ is a ucp map on $A(\textbf{$K$}),$ which is a matrix state.
	By the categorical duality \cite[Proposition 3.5]{WW}, there is a point $X \in \textbf{$K$}$ such that $$\overline{\Psi}_n(\theta) = \theta(X)$$ 
	for every $n \in \mathbb{N}$ and $\theta \in M_n(A(\textbf{$K$})).$ Moreover, for any $\theta \in J$ we have 
	$$\theta(X) = \overline{\Psi}(\theta) = 0,$$
	 so $X \in \textbf{$J$}^{\perp}.$ 
	Now fix $n \in \mathbb{N}$ and let $\theta \in M_n(A(\textbf{$K$}))$ be such that $\theta|_{\textbf{$J$}^{\perp}} \succeq 0.$ 
	Then $\overline{\Psi}_n(\theta) = \theta(X) \succeq 0$ for all $\Psi \in \text{UCP}(A(\textbf{$K$}) / J),$ which implies that $q(\theta) \succeq 0.$
	Indeed, if $q(\theta) \notin (A(\textbf{$K$}) / J)^+_n,$ then by the matricial Hahn-Banach theorem \cite{EW}, there is a continuous linear map $F: A(\textbf{$K$}) / J \to \mathbb{M}_n$ and a selfadjoint $\alpha \in \mathbb{M}_n$ such that
	\begin{equation}\label{eq: pos cond}
	F_r(\theta^\prime) \succeq \alpha \otimes \mathbb{I}_r
	\end{equation}
	for all $r \in \mathbb{N}$ and $\theta^\prime \in (A(\textbf{$K$}) / J)^+_r,$ but
	$$
	F_n(q(\theta)) \nsucceq \alpha \otimes \mathbb{I}_n.
	$$
	Since $(A(\textbf{$K$}) / J)^+_n$ is a cone, we cane take $\alpha=0.$ But now condition \eqref{eq: pos cond} implies that $F(q(1_{\textbf{$K$}}))^{-1} F$ is a ucp map on $A(\textbf{$K$}) / J$ such that $F_n(q(\theta)) \nsucceq 0.$ This is a contradiction, whence $q(\theta) \succeq 0.$
	
	By the definition of the quotient ordering
	there is a $\psi_1 \in  M_n(A(\textbf{$K$}))^+$ such that $q(\psi_1) = q(\theta),$ i.e., 
	$$\theta - \psi_1 \in \ker \Phi_n.$$
	As noted in Example \ref{ex418}, the fact that $J$ is spanned by its positive elements implies $\ker \Phi_n$ is also spanned by its positive elements, meaning $\ker \Phi_n = (\ker \Phi_n)^+ - (\ker \Phi_n)^+.$ Hence there is a $\psi_2 \in (\ker \Phi_n)^+$ such that $\theta - \psi_1 \preceq \psi_2.$ Then 
	$$\psi = \psi_1 + \psi_2 \in M_n(A(\textbf{$K$}))^+ + (\ker \Phi_n)^+ \subseteq M_n(A(\textbf{$K$}))^+$$ 
	clearly satisfies $\psi \succeq \theta \text{ and } \psi|_{\textbf{$J$}^{\perp}} = 
	 \psi_1|_{\textbf{$J$}^{\perp}} = \theta|_{\textbf{$J$}^{\perp}}.$ 
	
	(b) $\Rightarrow$ (a) To see that $J$ is spanned by its positive elements let $\theta \in J.$ Then $\theta |_{\textbf{$J$}^{\perp}} =0$ and by (b), there is a $\psi \in A(\textbf{$K$})^+$ such that $\psi \succeq \theta  \text{ and }  \psi|_{\textbf{$J$}^{\perp}} = \theta|_{\textbf{$J$}^{\perp}} = 0.$ Hence 
	$$\theta = \psi - (\psi - \theta) \in J^+ - J^+.$$

	Since $A(\textbf{$K$}) / J$ is easily seen to be a matrix ordered $*$-vector space in the quotient ordering, it remains to prove that it has an Archimedean matrix order unit. 
	 We first show  that $q(1_{\textbf{$K$}})$ is an Archimedean order unit. So suppose $\theta \in A(\textbf{$K$})$
	satisfies 
	$$
	q(\theta) + \frac{1}{n}\, q(1_{\textbf{$K$}}) \succeq 0
	$$
	for all $n \in \mathbb{N}.$ Then for each $n$ there is a map $\psi_n \in J$ such that 
	$$
	\theta + \frac{1}{n}\, 1_{\textbf{$K$}} + \psi_n \succeq 0,
	$$
	which implies $\theta|_{\textbf{$J$}^{\perp}} + \frac{1}{n}1_{\textbf{$J$}^{\perp}} \succeq 0$ for all $n \in \mathbb{N}.$ Hence $\theta|_{\textbf{$J$}^{\perp}} \succeq 0$ and by (b), there is a $\psi \in A(\textbf{$K$})^+$ such that 
	$\psi|_{\textbf{$J$}^{\perp}} = \theta|_{\textbf{$J$}^{\perp}}.$ But then $q(\theta) = q(\psi) \succeq 0$ as we now explain. 
	
	By the categorical duality \cite[Proposition 3.5]{WW}, $\Phi$ is the evaluation at some point $X \in \textbf{$K$}.$ Moreover, it is easily seen that $X$ must be in $\textbf{$J$}^{\perp}.$ This implies that $\Phi(\theta) = \Phi(\psi),$ whence
 	$\theta - \psi \in \ker \Phi$ and $q(\theta) = q(\psi).$

By an easy adaptation of the above argument to the $n$-th ampliation $q_n$ of $q,$ we see that $\mathbb{I}_n \otimes q(1_{\textbf{$K$}})$ is an Archimedean order unit in $M_n(A(\textbf{$K$}) / J)$ for every $n \in \mathbb{N}.$ In other words, $q(1_{\textbf{$K$}})$ is an Archimedean matrix order unit. Hence $A(\textbf{$K$}) / J$ satisfies the Choi-Effros axioms of an operator system (see \cite[Chapter 13]{Pa}).
\end{dokaz}

The next result combines and summarizes the conclusions from Example \ref{ex418} and Proposition \ref{tr420} giving the connection between certain closed matrix convex multifaces and positively generated kernels of partially order reflecting ucp maps.

\begin{trditev} \label{tr423} Let \textbf{$K$} be a compact matrix convex set.
	\begin{enumerate}[\text\textnormal{{(a)}}]
	\item Let $n \in \mathbb{N}$ and $\Phi: A(\textbf{$K$}) \to \mathbb{M}_n$ be a partially order reflecting ucp map with kernel $J$ spanned by its positive elements. Then $\textbf{$J$}^{\perp} \subseteq\textbf{$K$}$ is a closed \tmii that satisfies both of the equivalent conditions of Proposition \ref{tr420}.
	\end{enumerate}
	\begin{enumerate}[\text\textnormal{{(b)}}]
	\item Suppose $\textbf{$F$}  \subseteq \textbf{$K$}$ is a closed \tmii that satisfies the condition \text{\textnormal{(b)}} of Proposition \ref{tr420} (with $\textbf{$J$}^{\perp}$ replaced by \textbf{$F$}).
	Then
	$$
	J := \{\theta \in A(\textbf{$K$}) \ | \ \theta|_{\textbf{$F$}} = 0\}
	$$
	is spanned by its positive elements and is the kernel of a ucp map $\Phi: A(\textbf{$K$}) \to \mathcal{R}$ that satisfies the partially order reflecting property \eqref{eq413}, where $\mathcal{R}$ is an operator system.  \looseness=-1
	\end{enumerate}
\end{trditev}

\begin{dokaz} Let $q: A(\textbf{$K$})  \to A(\textbf{$K$}) / J$ denote the canonical projection. It is clearly completely positive, i.e., $\mathbb{I}_n\otimes q$ is order preserving for all $n,$ but $q$ is also partially order reflecting as we now explain. Note that for each $n$ there is a linear isomorphism 
	$$M_n(A(\textbf{$K$}) / J) \cong M_n(A(\textbf{$K$})) / M_n(J).$$
	So for any $n$ and $\theta \in M_n(A(\textbf{$K$}))$ with 
	$
	q(\theta) \in \big( M_n(A(\textbf{$K$}) / J)\big)^+ 
	$ 
	there is a $\psi \in M_n(J)$ such that $\theta + \psi \succeq 0.$ Whence, $q(\theta + \psi) = q(\theta).$

	For (a) note that for every $n,$ the induced map 
	$$\overline{\Phi}_n :   M_n(A(\textbf{$K$}) / J)   \to \text{im}\,\Phi_n$$ is bijective, ucp and partially order reflecting, hence an order isomorphism. This implies  that $A(\textbf{$K$}) / J$ is an operator system.
	Part (a) now follows from Example \ref{ex418} and Proposition \ref{tr420}.
	
	For (b) we first prove that $\textbf{$J$}^\perp = \textbf{$F$}.$ Clearly, $\textbf{$F$} \subseteq \textbf{$J$}^\perp.$ So suppose we have $X \in \textbf{$J$}^\perp \backslash \textbf{$F$}.$ As \textbf{$F$} is closed, by the matricial Hahn-Banach theorem \cite{EW} there is a matrix affine map $\theta \in A(\textbf{$K$})$ such that $\theta|_{\textbf{$F$}} \succeq 0,$ but $\theta(X) \nsucceq 0.$ By assumption there is a positive element $\psi \in A(\textbf{$K$})^+$ with $\psi \succeq \theta \text{ and } \psi|_{\textbf{$F$}} = \theta|_{\textbf{$F$}}.$ So we have both $\psi(X) \succeq 0$ and $\psi -\theta \in J.$ But then by the choice of $X$, $(\psi -\theta)(X) = 0,$ which is a contradiction.
Hence $\textbf{$J$}^\perp = \textbf{$F$}.$

Now the proof of the (b) $\Rightarrow$ (a) implication in Proposition \ref{tr420} shows that $J$ is spanned by its positive elements and that the quotient $\mathcal{R}:=A(\textbf{$K$}) / J$ is an operator system. Hence $q: A(\textbf{$K$})  \to \mathcal{R}$ is the desired ucp map with the partial order reflection property \eqref{eq413} for all $n \in \mathbb{N}.$
\end{dokaz}

\begin{opomba} If \textbf{$K$} is a matrix convex set in a finite-dimensional space $V,$ then the operator system $\mathcal{R}$ in part (b) of Proposition \ref{tr423} is finite-dimensional and the ucp map $\Phi$ is in fact a matrix state on $A(\textbf{$K$}).$ So in this case part (b) gives a proper converse to part (a). 
\end{opomba}

\begin{primer}\label{ex519}
	Let $K_1 \subseteq \mathbb{R}^2$ be the triangle with a vertex $X$ in Figure \ref{fig2} and let \textbf{$K$} be the matrix convex hull of $K_1.$ It is routine to check that $\ker\,\Phi_X$ is spanned by its positive elements, hence by Example \ref{ex418}, $(\ker\,\Phi_X)^{\perp}$ is a \tmii containing $X.$
	
	 Further, $\Phi_X$ is partially order reflecting as we now explain.
	Let $n \in \mathbb{N}$ and $\theta \in M_n(A(\textbf{$K$}))$ be such that $$(\Phi_X)_n(\theta) = \theta(X) \succeq 0.$$
	Without loss of generality assume $\theta(X)$ is diagonal with diagonal entries $\lambda_1, \ldots, \lambda_n \geq 0.$ For each $i$ choose an affine function $f_i,$ which is positive on $K_1,$ such that $f_i(X)=\lambda_i$ and define
	$$
	\psi = f_1 \oplus \cdots \oplus f_n.
	$$
	Note that for any matrix affine map $\varphi = (\varphi_n)_n \in A(\textbf{$K$})$  and unit vector $\xi \in \mathbb{C}^n,$ the property 
	$$\xi^\ast \varphi_n(A) \xi = \varphi_1(\xi^\ast A \xi)$$ 
	implies $\varphi_1 \geq 0$ on $K_1$ if and only if $\varphi \in A(\textbf{$K$})^+.$
	It is now easy to see that $\psi \in M_n(A(\textbf{$K$}))^+$ and $\psi(X) = \theta(X).$ This shows that $\Phi_X$ is indeed partially order reflecting.
	
	After a rotation, the above reasoning applies to any of the three vertices of $K_1,$ so each of them defines a partially order reflecting evaluation map. Moreover, a similar argument shows the same for each of the vertices of a simplex $S$ in an Euclidean space $\mathbb{R}^n.$ So each of them lies in a matrix convex multiface of mconv$(S)$ satisfying both of the (equivalent) conditions of Proposition \ref{tr420}.
	\begin{figure}[h!]
		\begin{center}
			\begin{tikzpicture}
				\draw[->] (-2,0) -- (4,0) node[right]{$x_1$};
				\draw[->] (0,-1) -- (0,4) node[above]{$x_2$};
				\filldraw[fill=yellow!40!white] (0,0) -- (3,0) -- (0,3)-- cycle;
				\fill (0,3) circle (1.5pt)node[left]{$X$};
				\fill (0,0) circle (1pt)node[below left]{$0$};
			\end{tikzpicture}
			\caption{The level-one point $X$ determines a matrix exposed multiface of the matrix convex hull of the above triangle.}
			\label{fig2}
		\end{center}
	\end{figure}
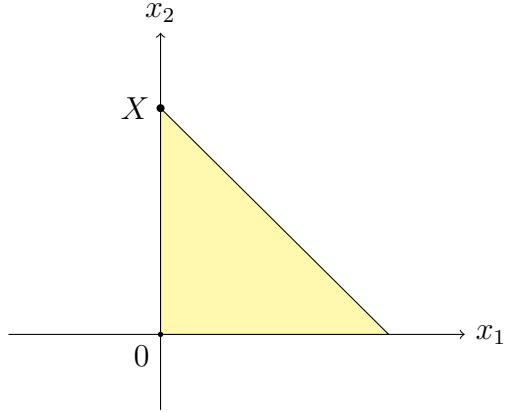
\end{primer}

With similar reasoning as in Proposition \ref{tr34} we can prove the following generalization of the extreme point preservation property of classical faces.

\begin{trditev} \label{tr58}
	Let \textbf{$K$} be a matrix convex set and $\textbf{$F$} \subseteq \textbf{$K$}$ a matrix \text{\textnormal{(}}convex\,\text{\textnormal{)}} multiface. Every matrix extreme point of \textbf{$F$} is a matrix extreme point of \textbf{$K$}.
\end{trditev}

\begin{primer}
	Note that a matrix multiface \textbf{$F$} of either type must contain all the matrix extreme points, whose proper matrix convex combinations describe the elements of \textbf{$F$}. In particular, the matrix interval $[a \mathbb{I}, b \mathbb{I}]$ := $([a \mathbb{I}_n, b \mathbb{I}_n])_{n \in \mathbb{N}}$ has very few matrix convex multifaces, i.e., $\emptyset,$ mconv$(\{a\}),$ mconv$(\{b\})$ and $[a \mathbb{I}, b \mathbb{I}],$
	and accordingly, its matrix multifaces that are not matrix convex are of the form $\widehat{F}$ in the notation from part (b) of Remark \ref{op422a}, where $F$ is a \ti. More generally, if  \textbf{$K$} is a matrix convex set and $x \in K_1$ is an extreme point, then mconv\,$\{x\}$ is a \tmii of \textbf{$K$} (see \cite[Example 9.9]{KKM}). 
\end{primer}

\subsection{Matrix exposed multifaces}\label{subsec52} This section introduces the exposed counterparts of the multilevel matrix faces and investigates their basic properties.

\begin{definicija}\label{def59}
	Let $\textbf{$K$} = (K_r)_{r \in \mathbb{N}}$ be a matrix convex set in a dual vector space $V$ and \textbf{$F$} a levelwise convex subset of \textbf{$K$}. 
	
	(a) Then \textbf{$F$} is a \textbf{\etmi} if there exists a positive integer $r,$ a continuous linear map $\Phi:V \to \mathbb{M}_r$ and a self-adjoint matrix $\alpha \in \mathbb{M}_r$ satisfying the following conditions:
	\begin{enumerate}[(i), leftmargin=2cm]
		\item for every positive integer $n$ and $B \in K_n$ we have $\Phi_n(B) \preceq \alpha \otimes \mathbb{I}_n;$ 
		
		\item for each $n \in \mathbb{N}$ we have $\{ B \in K_n \ | \  \alpha \otimes \mathbb{I}_n - \Phi_n(B) \succeq 0 \text{ is singular}\} = F_n.$
	\end{enumerate}
	
	(b) If \textbf{$F$} is a matrix convex \etmi, then it is a \textbf{\etmii.}
\end{definicija}

We call a pair $(\Phi, \alpha)$ in the notation above an \textbf{exposing pair} (of size $r$) for the matrix exposed face \textbf{$F$}. An exposing pair is \textbf{minimal} if there is no $s < r$ together with a linear map $\Psi:V \to \mathbb{M}_s$ and a self-adjoint matrix $\beta \in \mathbb{M}_s$ satisfying 

(i) $\Psi_m(B) \preceq \beta \otimes \mathbb{I}_m$ for all $B \in K_m$ and positive integers $m,$ 

(ii) $\{ B \in K_n \ | \  \alpha \otimes \mathbb{I}_n - \Phi_n(B) \succeq 0 \text{ is singular}\} = F_n$ for every $n \in \mathbb{N}.$

\begin{opomba}\label{op423}
	(a) Assume the notation from part (b) of Remark \ref{op422a}. Then for every \eti $F \subseteq K_n,$ the corresponding multilevel set $\widehat{F}$ is a \etmi. Moreover, the matrix exposed multifaces $\textbf{$F$} = (F_r)_{r \in \mathbb{N}}$ with $F_r = \emptyset$ for $r > 1$ coincide with subsets of \textbf{$K$}, whose first components are ordinary exposed faces.
	Also, as in part (b) of Remark \ref{op25}, a matrix exposed multiface is closed under conjugation by unitaries.

	(b) As in part (c) of Definition \ref{def21} one might consider weak matrix exposed multifaces to obtain that for any weak matrix exposed multiface \textbf{$F$}, where $F_n = \{A\}$ for some $n \in \mathbb{N}$ and $F_m = \emptyset$ whenever $m \neq n,$ the point $A$ is matrix exposed. 
	
	(c) We can show as in part (e) of Remark \ref{op25} that for a matrix exposed multiface \textbf{$F$} of any type the intersection 
	$$
	\mathcal{N}_n = \bigcap_{A \in F_n} \ker(\alpha \otimes \mathbb{I}_n - \Phi_n(A))
	$$
	is nontrivial for any positive integer $n$.
	
	(d) Note that the zero map and the zero matrix define an exposing pair for \textbf{$K$}.
	
	(e) Suppose $L$ and $M$ are two linear pencils such that $\mathcal{D}_M \subseteq \mathcal{D}_L$ and $\mathcal{D}_M(n) \cap \partial\mathcal{D}_L(n) \neq \emptyset$ for some $n \in \mathbb{N}.$ If $F= \mathcal{D}_M \cap \partial\mathcal{D}_L = (\mathcal{D}_M(m) \cap \partial\mathcal{D}_L(m))_m$ is levelwise convex, then it is  a matrix exposed multiface in $\mathcal{D}_M$ (cf.~Example \ref{ex410}).
	\end{opomba}

\begin{trditev}\label{pr427}
	Let $\textbf{$F$} \subseteq \textbf{$K$}$ be a matrix \text{\textnormal{(}}convex\,\text{\textnormal{)}} exposed multiface and $\Phi:V \to \mathbb{M}_r$ together with $\alpha \in \mathbb{M}_r$ a minimal exposing pair. Then for every $n \in \mathbb{N}$ and nonzero $x = \sum_{i=1}^n x_i \otimes e_i \in \mathbb{C}^r \otimes \mathbb{C}^n$ in $\mathcal{N}_n,$ the span of its components $x_1, \ldots, x_n$ is $m$-dimensional, where $m=\min (r,n).$
\end{trditev}

\begin{dokaz}
	Suppose $r \leq n$ (the other case is treated similarly). If the span $\mathcal{M}$ of $x_1, \ldots, x_n$ is of dimension $m<r,$ then the projection $P$ onto $\mathcal{M}$ gives rise to an exposing pair $(P\,\Phi P^*, P\alpha P^*)$ of size $m<r.$ But this contradicts the minimality of $(\Phi,\alpha).$
\end{dokaz}

\begin{trditev}\label{pr428}
	Let $\textbf{$K$} = (K_n)_{n \in \mathbb{N}}$ be a matrix convex set and $\textbf{$F$} \subseteq \textbf{$K$}$ a matrix \text{\textnormal{(}}convex\,\text{\textnormal{)}} exposed multiface. Then each $F_n$ is an ordinary exposed face of $K_n.$
\end{trditev}

\begin{dokaz}
	Let $\textbf{$F$} \subseteq \textbf{$K$}$ be a matrix exposed multiface and $\Phi:V \to \mathbb{M}_r$ together with $\alpha \in \mathbb{M}_r$ a minimal exposing pair. For fixed $n \in \mathbb{N},$ choose a nonzero $x_n$ in $\mathcal{N}_n$ as in Proposition \ref{pr427} and define the functional $\varphi: M_n(V) \to \mathbb{C}$ by
	\begin{equation*}
		\varphi_n(B) = x_n^\ast \Phi_n(B) x_n,
	\end{equation*}
	and the real number $a_n = x_n^\ast (\alpha \otimes \mathbb{I}_n) x_n.$ Now conclude as in the proof of Proposition \ref{tr39} that the pair $(\varphi_n, a_n)$ exposes $F_n$ in $K_n.$
\end{dokaz}

\subsubsection{Interplay between matrix multifaces and matrix exposed multifaces}\label{subsec53}

As an extension of Subsection \ref{subsec43}, here we give a few remarks on the connection between matrix multifaces and matrix exposed multifaces.

\begin{trditev} \label{tr513}
	Let \textbf{$K$} be a matrix convex set and $\textbf{$F$} \subsetneq \textbf{$K$}$ a matrix exposed multiface with minimal exposing pair of size $r.$ Denote by $\widetilde{\textbf{$F$}}$ the graded set with $\widetilde{F}_k = F_k$ for $k\leq r$ and $\widetilde{F}_k = \emptyset$ for $k>r.$ Then $\widetilde{\textbf{$F$}}$ is a matrix multiface.
\end{trditev}

\begin{dokaz}
	The proof is essentially same as that of Proposition \ref{tr412}, where the counterparts of the two key observations needed, namely part (b) of Remark \ref{op211} and Proposition \ref{prop27}, are covered by the definition of a matrix multiface and by Proposition \ref{pr427}, respectively.
\end{dokaz}

\begin{opomba}
	Note that if \textbf{$F$} is a \tmii and for some $n \in \mathbb{N}$, $F_n$ is a \etii, then $F_n=K_n.$
	Indeed, if $F_n$ is a \etii, which is a proper subset of $K_n,$ then it must not contain any reducible elements by part (b) of Remark \ref{op22}. But since \textbf{$F$} is matrix convex, $\oplus_n F_1 \subseteq F_n$ for all $n \in \mathbb{N}.$ So $F_n$ can only be  exposed by a map that is constant on the whole $K_n$, whence $F_n=K_n.$
\end{opomba}

%

\begin{trditev} \label{pr434} 
	Let \textbf{$K$} be a compact matrix convex set.
	Suppose $\textbf{$F$}  \subseteq \textbf{$K$}$ is a closed \tmii such that for each $n \in \mathbb{N}$ and $\theta \in M_n(A(\textbf{$K$}))$ with $\theta|_{\textbf{$F$}} \succeq 0$ there is a positive element $\psi \in M_n(A(\textbf{$K$}))^+$ with 
	$$\psi \succeq \theta \quad \text{ and } \quad \psi|_{\textbf{$F$}} = \theta|_{\textbf{$F$}}.
	$$ 
	Then for each $n \in \mathbb{N},$ $F_n$ is an exposed face of $K_n.$
\end{trditev}

\begin{dokaz} Fix $n \in \mathbb{N}$ and without loss of generality suppose $F_n \subsetneq K_n.$ 
	We first show that for every $X \in K_n \backslash F_n$ there is a continuous affine function $\varphi_X: M_n(V) \to \mathbb{C}$ such that
	\begin{equation*}
		\varphi_X|_{F_n} = 0 \quad \text{ and } \quad \varphi_X(X) > 0.
	\end{equation*}
	Indeed, since \textbf{$F$} is closed, by the matricial Hahn-Banach separation theorem \cite{EW}, there is a matrix affine map $\theta \in A(\textbf{$K$})$ with $\theta_1:V\to \mathbb{C}$ that satisfies 
	$$\theta|_{\textbf{$F$}} \succeq 0 \quad \text{ and } \quad \theta(X) \nsucceq 0.$$ 
 Hence there is a $y \in \mathbb{C}^n$ such that $y^* \theta(X) y <0.$
	By assumption, there is a $\psi \in A(\textbf{$K$})^+$ with $\psi - \theta \in A(\textbf{$K$})^+ \text{ and }  \psi|_{\textbf{$F$}} = \theta|_{\textbf{$F$}}.$ But the latter means we can take $\varphi_X$ to be $y^*(\psi - \theta)y.$ 
	
	Since $K_n$ is compact, there exist finitely many $X_1, \ldots, X_k \in K_n \backslash F_n$ such that for each $X \in K_n \backslash F_n$ there is an $1 \leq i \leq k$ such that
	$$\varphi_{X_i}|_{F_n} = 0 \quad \text{ and } \quad \varphi_{X_i}(X) > 0.$$ 
	It is then straightforward that $\varphi = \frac{1}{k}(\varphi_{X_1} + \cdots + \varphi_{X_k})$ exposes $F_n$ in $K_n.$ 
\end{dokaz}

\end{document}